\documentclass[runningheads, nospthms]{svjour2}
\journalname{Inventiones mathematicae}

\usepackage{amssymb, amsmath, amscd} 
\usepackage{mathptmx, mathrsfs}  

\usepackage{amsthm}

\def\RCS$#1: #2 ${\expandafter\def\csname RCS#1\endcsname{#2}}
\RCS$Revision: 1.65 $
\RCS$Date: 2007/11/06 18:26:57 $

\newcommand{\3}[1]{3\nobreakdash-\hspace{0pt}}
\newcommand{\1}[1]{1\nobreakdash-\hspace{0pt}}

\newcommand{\C}{{\mathbb C}} 
\newcommand{\Q}{{\mathbb Q}}
\newcommand{\F}{{\mathbb F}}
\newcommand{\R}{{\mathbb R}}
\newcommand{\Z}{{\mathbb Z}}
\newcommand{\T}{{\mathcal T}}

\newcommand{\M}{{\mathcal{M}}}
\newcommand{\A}{{\mathcal{A}}}
\newcommand{\E}{{\mathcal{E}}}
\newcommand{\maps}{\colon\thinspace} 
\DeclareMathOperator{\tr}{tr} 
\DeclareMathOperator{\Aut}{Aut}
\DeclareMathOperator{\Out}{Out}
\DeclareMathOperator{\Alt}{Alt}
\DeclareMathOperator{\vecspan}{span}
\DeclareMathOperator{\Sym}{Sym} 
\newcommand{\surjects}{\twoheadrightarrow}
\newcommand{\PSL}[2]{\mathrm{PSL}_{#1} #2}
\newcommand{\PGL}[2]{\mathrm{PGL}_{#1} #2}
\newcommand{\SL}[2]{\mathrm{SL}_{#1} #2}
\newcommand{\GL}[2]{\mathrm{GL}_{#1} #2}
\newcommand{\SU}[1]{\mathrm{SU}({#1})}
\newcommand{\Sp}[2]{\mathrm{Sp}_{#1} #2}
\newcommand{\abs}[1]{{\left| #1 \right|}}
\newcommand{\pair}[1]{\left\langle #1 \right\rangle}

 \newcommand{\spandef}[2]{{ \left\langle
      {#1} \ \left| \ {#2} \right. \right\rangle }}
\newcommand{\setdef}[2]{{ \left\{ {#1} \ \left| \ {#2} \right.
    \right\} }} \newcommand{\mtext}[1]{\quad\mbox{#1}\quad}
\newcommand{\num}[1]{\abs{#1}}
\DeclareMathOperator{\Gr}{Gr}

\newcommand{\comment}[1]{}

\renewcommand{\tilde}{\widetilde}

\newtheoremstyle{myplain}
  {}
  {}
  {\itshape}
  {}
  {\bfseries}
  {}
  {.5em}
  {}

\newtheoremstyle{myremark}
  {}
  {}
  {\itshape}
  {}
  {\bfseries}
  {}
  {.5em}
  {}

\makeatletter
\renewenvironment{proof}[1][\proofname]{\par
  \pushQED{\qed}%
  \normalfont \topsep6\p@\@plus6\p@\relax
  \trivlist
  \item[\hskip\labelsep
        \itshape
    #1]\ignorespaces
}{%
  \popQED\endtrivlist\@endpefalse
}
\makeatother

\theoremstyle{myplain}
\swapnumbers
\newtheorem{theorem}{Theorem}[section]
\newtheorem{conjecture}[theorem]{Conjecture}
\newtheorem{lemma}[theorem]{Lemma}
\newtheorem{corollary}[theorem]{Corollary}
\newtheorem{proposition}[theorem]{Proposition}
\newtheorem{claim}[theorem]{Claim}
\newtheorem{question}[theorem]{Question}

\newtheorem{virtual_haken}[theorem]{Virtual Haken Conjecture}

\newtheorem{simpleprobtheorem}{Theorem}

\newtheorem{theoremraq}{Theorem}

\newtheorem{theoremnohom}{Theorem}

\theoremstyle{myremark}
\newtheorem{remark}[theorem]{Remark}
\newtheorem{rmk}[theorem]{Remark}

\makeatletter
  \let\c@theorem=\c@subsection
  \let\c@table=\c@subsection
  \let\p@table=\p@subsection
  
  \let\cl@table=\cl@subsection
  \let\c@equation=\c@subsection
  \let\p@equation=\p@subsection
  
  \let\cl@equation=\cl@subsection
\makeatother

\begin{document}

\title{Finite covers of random 3-manifolds\thanks{Authors were partially supported by the US NSF and the Sloan Foundation.}}

\author{Nathan M.~Dunfield \and William P.~Thurston}

\institute{Mathematics 253-37, Caltech, Pasadena, CA 91125, USA.
 \email{dunfield@caltech.edu}
\and
Department of Mathematics, 310 Malott Hall , Cornell University, Ithaca, NY 14853-4201, USA. \email{wpt@math.cornell.edu}}

\date{Version: \RCSRevision, Compile: \today, Last commit:
    \RCSDate}

\maketitle

\begin{abstract} 
  A 3-manifold is Haken if it contains a topologically essential
  surface.  The Virtual Haken Conjecture posits that every irreducible
  3-manifold with infinite fundamental group has a finite cover which
  is Haken.  In this paper, we study \emph{random 3-manifolds} and
  their finite covers in an attempt to shed light on this difficult
  question.  In particular, we consider \emph{random Heegaard
    splittings} by gluing two handlebodies by the result of a random
  walk in the mapping class group of a surface.  For this model of
  random 3-manifold, we are able to compute the probabilities that the
  resulting manifolds have finite covers of particular kinds.  Our
  results contrast with the analogous probabilities for groups coming
  from random balanced presentations, giving quantitative theorems to
  the effect that 3-manifold groups have many more finite quotients
  than random groups.  The next natural question is whether these
  covers have positive betti number.  For abelian covers of a fixed
  type over 3-manifolds of Heegaard genus 2, we show that the probability
  of positive betti number is 0.
  
  In fact, many of these questions boil down to questions about the
  mapping class group.  We are led to consider the action of the mapping
  class group of a surface $\Sigma$ on the set of quotients $\pi_1(\Sigma) \to
  Q$.  If $Q$ is a simple group, we show that if the genus of $\Sigma$ is
  large, then this action is very mixing.  In particular, the action
  factors through the alternating group of each orbit.  This is
  analogous to Goldman's theorem that the action of the mapping class
  group on the $\SU{2}$ character variety is ergodic.
\end{abstract}

\setcounter{tocdepth}{2}
\tableofcontents

\section{Introduction}

Here, we study various notions of a random \3-manifold, and try to
understand the distribution of topological and group-theoretic
properties for such manifolds.  Our primary motivation is to try to
determine the underlying issues behind the Virtual Haken Conjecture
and related problems about properties of finite covers of
\3-manifolds.  While any hyperbolic \3-manifold has many finite
covers---its fundamental group is residually finite---we do not know
what most of the covering groups are, much less the properties that we
may reasonably expect these covers to have.  The reason that the
fundamental group of a hyperbolic \3-manifold $M$ is residually finite
is that it is a finitely generated group of matrices, and this gives
many quotients of $\pi_1(M)$ of the form $\PSL{2}{\F}$, where $\F$ is a
finite field.  Lubotzky has shown \cite{Lubotzky95} that such
quotients have measure zero among all quotients of $\pi_1(M)$, but his
proof provides little insight into what those other quotients might
be. For instance, is it reasonable to expect $\pi_1(M)$ to have a
quotient which is an alternating group $A_n$?  Should there be many
such quotients?  Also, for a particular kind of finite quotient, how
likely is it that the associated cover be Haken or have positive betti
number?  Although there have been many partial results on the Virtual
Haken Conjecture, these questions have been hard to address in general
by direct deductive reasoning.

Another way of thinking about these questions is from a probabilistic
point of view.  Since the set of homeomorphism types of compact
\3-manifolds is countably infinite, there is no uniform,
countably-additive, probability measure on this set.  Thus the first
issue is to define a plausible context in which we can discuss
probability.  The density of \3-manifolds with any given property will
depend on the order in which we consider them, unless the property is
either true or false for all but a finite set of \3-manifolds.  It
seems best to us to analyze orders of enumeration that are plausible
and tractable, while acknowledging that there may be other equally
plausible and tractable orders that give different answers.  In
Section~\ref{models}, we outline several reasonable models for a
random \3-manifold, and then in most of the rest of this paper
concentrate on a model of random \3-manifolds which comes from looking
at Heegaard splittings generated by random walks in mapping class
groups.

\subsection{Random Heegaard splittings}

Every closed orientable \3-manifold has a Heegaard splitting, that is,
it can be constructed by gluing two handlebodies of genus $g$ together
using a homeomorphism between their boundaries.  We will look at
\3-manifolds of a fixed Heegaard genus $g$, and consider gluings
obtained by a random word in a finite set of generators for the
mapping class group of the surface of genus $g$.  In principle, the
density of a particular property of the manifolds obtained in this way
could depend on the choice of generators for the mapping class group.
Indeed, this is plausible since in a non-amenable group such as this
one, the correlation between random words of large length in different
sets of generators usually tends to $0$.  However, the properties we
analyze have limiting densities independent of this choice of
generators.  In particular, for a finite group $Q$ the probability
that the manifold obtained from a random genus $g$ Heegaard splitting
has a cover with covering group $Q$ is well-defined
(Prop.~\ref{makes_sense}), and we denote this probability as $p(Q,
g)$.  When $Q$ is a simple group, we are able to calculate the limit
of these probabilities as the genus $g$ goes to infinity:
\begin{simpleprobtheorem}
  Let $Q$ be a non-abelian finite simple group.  Then as the genus $g$ goes to
  infinity, the probability of a $Q$-cover converges:
  \[
  p(Q, g) \to 1 - e^{-\mu} \mtext{where} \mu = \num{H_2(Q; \Z)} /\num{\Out(Q)}.
  \]
  Moreover, the limiting distribution of the number of $Q$-covers
  converges to the Poisson distribution with mean $\mu$.  
\end{simpleprobtheorem}
For example, if $Q = \PSL{2}{\F_p}$ where $p$ is an odd prime then $\mu
= 1$.  Thus for large genus the probability of a $\PSL{2}{\F_p}$ cover
is about $1 - e^{-1} \approx 0.6321$.  Hyperbolic \3-manifolds must have
infinitely many covers of this form, namely the congruence quotients.
However, at least naively, one expects far fewer congruence quotients
than given by Theorem~\ref{simple_probs}.  \comment{surely, the degree
  $d$ of the trace field goes to infinitely as the length of the walk
  does; ergo the density of primes for which there is a congruence
  quotient for a particular prime is about $1/d$.}Another interesting
example is the case where $Q = A_n$ is an alternating group; here
again $\mu = 1$ and the probability of an $A_n$ cover is greater than
$63\%$.  Moreover, we show that covers with different groups $Q_i$ do
not correlate with one another, at least for large genus.  As a
consequence, we can prove results such as the following, which is a
special case of Theorem~\ref{random-quotient-sequence}:
\begin{theoremraq}
  Let $\epsilon > 0$.  For all sufficiently large $g$, the probability that
  the \3-manifold obtained from a random genus-$g$ Heegaard splitting
  has a $A_n$-cover with $n \geq 5$ is at least $1-\epsilon$.  Moreover, the
  same is true if we require some fixed number $k$ of such covers.
\end{theoremraq}
These results should be contrasted with the analogous results for
finitely presented groups with an equal number of generators and
relators, where the relators are chosen at random.  There, the
probability of a $A_n$-quotient goes to zero like $1/n!$ as $n$ goes
to infinity, rather than remaining constant
(Theorem~\ref{thm-simple-balanced}).  Thus one way of interpreting our
results is that they affirm the belief that \3-manifolds have many
finite quotients.  In Section~\ref{sec-quotient-3-mfld}, we give some
heuristic ways to understand why this should be true, working from a
more naive point of view.  These stem from the fact that the
group relators given by a Heegaard splitting come from
disjoint embedded curves on a surface.  In particular, we
highlight special features of attaching the last two 2-handles in
forming a \3-manifold that suggest many extra finite quotients.  

\subsection{Mapping class group}
The proof of Theorem~\ref{simple_probs} boils down to understanding
the action of the mapping class group of a surface on a certain finite
set.  Let $\Sigma_g$ be a closed surface of genus $g$, and let $\M_g$ be
its mapping class group.  Consider the set $\A_g$ of epimorphisms of
$\pi_1(\Sigma_g)$ onto our fixed simple group $Q$, modulo automorphisms of
$Q$.  Then $\M_g$ acts on $\A_g$ via the induced automorphisms of
$\pi_1(\Sigma_g)$, and we show:
\begin{theorem}\label{intro_alternating}
  Let $Q$ be a non-abelian finite simple group.  Then for all
  sufficiently large $g$, the orbits of $\A_g$ under the action of
  $\M_g$ correspond bijectively to $H_2(Q, \Z) / \Out(Q)$.  Moreover,
  the action of $\M_g$ on each orbit is by the full alternating group
  of that orbit.
\end{theorem}
For a finite group $Q$ which is not necessarily simple the orbits can
be classified in the same way (Theorem~\ref{hom_classify_orbits});
that the action on an orbit is by the full alternating group is
special to simple groups $Q$ (Theorem~\ref{action_is_alternating}).
You can view Theorem~\ref{intro_alternating} as saying that when the
genus is large the action of $\M_g$ on $\A_g$ is nearly as mixing as
possible.  As such, it is directly analogous to Goldman's theorem that
the action of $\M_g$ on the $\SU{2}$-character variety is ergodic for
any genus $\geq 2$ \cite{Goldman1997}.  Perhaps surprisingly, the proof
of Theorem~\ref{intro_alternating} uses the Classification of
Simple Groups even for concrete cases such as $Q = A_n$.  However,
Theorem~\ref{simple_probs}, which is a corollary of
Theorem~\ref{intro_alternating}, also follows from a weaker version
which does not use the Classification.

What about other types of finite groups?  For abelian groups, we give
a complete picture of the distribution of $H_1(M)$ for a \3-manifold
coming from a random Heegaard splitting
(Section~\ref{homology_Heegaard}).  For a general finite group $Q$, we
do not know how to show the existence of a limiting distribution as
the genus goes to infinity, but we can show that the expected number
of $Q$-quotients does converge
(Theorem~\ref{expectations_converge}).\footnote{ Not in published
  version: Jordan Ellenberg points out that Conway and Parker
  studied the analgous question for the braid groups,
  which are mapping class groups of punctured discs.  In particular,
  they proved a version of our Theorem~\ref{hom_classify_orbits} in
  that context.  The theorem of Conway and Parker was written down and
  improved by Fried and Volklein \cite{FriedVolklein}, who used it to
  study the inverse Galois problem.  }

\subsection{Virtual positive betti number} 

These results give a good picture about the number of different types
of covers in many cases, so we now turn to the main question at hand:
\begin{virtual_haken}
  Let $M$ be a closed irreducible 3-manifold with $\pi_1(M)$ infinite.
  Then $M$ has a finite cover which is Haken, i.e.~contains an
  incompressible surface.
\end{virtual_haken}
This conjecture was first proposed by Waldhausen in the 1960s
\cite{Waldhausen68}.  It is often motivated as a way of reducing
questions to the case of Haken manifolds, where one has the most
topological tools available.  However, we prefer to view it as an
intrinsic question about the topology $M$ itself: does $M$ contain an
\emph{immersed} incompressible surface?  If so, can we lift it to an
embedded surface in some finite cover?  From a more algebraic point of
view, one of the fundamental tasks of \3-manifold topology is to
understand the special properties of their fundamental groups, as
compared to finitely presented groups in general; thus it is natural
to ask: does $\pi_1(M)$ always contain the fundamental group of a
closed surface?  (Having $\pi_1(M)$ contain a surface group is
equivalent to $M$ having an immersed incompressible surface, but is a
priori weaker than having a finite cover which is Haken. The
difference is another subtle and interesting question about $\pi_1(M)$,
namely subgroup separability.)
 
Perelman has announced a proof of Thurston's Geometrization
Conjecture using Hamilton's Ricci flow \cite{PerelmanI}
\cite{PerelmanII}; this should reduce the Virtual Haken Conjecture to
the (generic) case when $M$ is hyperbolic.  For hyperbolic $M$, the
Virtual Haken Conjecture fits nicely into a more general question of
Gromov: must a 1-ended word-hyperbolic group contain a surface group?
In any event, it seems to us that it would be very hard to prove the
Virtual Haken Conjecture without first establishing Geometrization;
for instance, the only known way to show that an ``easy to
understand'' atoroidal Haken manifold $M$ has a non-trivial finite
cover is to hyperbolize it to see that $\pi_1(M)$ is in fact a group of
matrices, hence residually finite!

Deciding whether a \3-manifold is Haken is difficult, so in the rest
of this paper we focus on the \emph{stronger} version of the
conjecture which asks that the finite cover has positive first betti
number.  (In the case of arithmetic \3-manifolds, this is also the
version that relates directly to the theory of automorphic forms.)  In
our prior work \cite{DunfieldThurston:experiments}, we found that this
conjecture holds for all $\approx 11{,}000$ of the small volume hyperbolic
\3-manifolds in the Hodgson-Weeks census.  One of our goals here is to
determine whether the patterns we observed there are in some sense
generic, or are a consequence of special properties of that sample.
For simple quotients, the results above give even larger probabilities
for such covers than those we observed in
\cite{DunfieldThurston:experiments} (see Section~\ref{compare-simple}
for a quantitative comparison).

However, for the crucial question of whether simple covers have
positive betti number, a different picture emerges for our random
manifolds here than we saw in \cite{DunfieldThurston:experiments}.  In
\cite[\S5]{DunfieldThurston:experiments}, we found that covers with a
particular fixed finite simple group $Q$ had positive betti number
with probabilities between 52--98\% depending on $Q$.  However, for
our Heegaard splitting notion of random our experimental evidence
strongly suggests that these probabilities are 0.  Moreover, we can show
\begin{theoremnohom}
  Let $Q$ be a finite abelian group.  The probability that the
  3-manifold obtained from a random Heegaard splitting of genus 2 has
  a $Q$-cover with $\beta_1 > 0$ is $0$.
\end{theoremnohom}

\subsection{Potential uses of random 3-manifolds}

In combinatorics, studying random objects is done not just for the
intrinsic interest and beauty of the subject but also for the
applications.  For instance, constructing explicit infinite families
of expander graphs is quite difficult; the first such construction was
based on the congruence quotients of $\PSL{2}{\Z}$ and uses Selberg's
$3/16$ Theorem (see e.g.~\cite{Lubotzky-expanders-book}).  On the
other hand, proofs of existence and practical construction can be done
by looking at certain classes of random graphs and showing that the
desired property occurs with non-zero probability.  Closer to the
study of \3-manifolds, Gromov initiated the study of groups coming
from certain types of random group presentations.  These have been
used to produce many examples of word-hyperbolic groups with
additional properties, such as having Property T \cite{Zuk2003}
\cite{Gromov2003}.  Very recently, Belolipetsky and Lubotzky have used
random techniques to show that given $n$ and a finite group $G$ there
exists a hyperbolic $n$-manifold whose full isometry group is exactly
$G$ \cite{BelolipetskyLubotzky2004}.

Perhaps similar techniques could
be applied to questions about \3-manifolds.  In particular, to
construct \3-manifolds with a certain list of properties, one could
try to show that these properties occur with positive probability for
a suitable model of random \3-manifold.  For such applications, the
fact that a random \3-manifold is an ill-defined concept becomes a
strength rather than a weakness, since by varying the model one can
change the characteristics of the resulting manifolds.  

Finally, another point of view on random \3-manifolds is that they
provide a \emph{quantitative} context in which to understand one of
the central questions in 3-dimensional topology: how do \3-manifold
groups differ from finite presented groups in general?  As we
mentioned, our results show that from the point of view of random
Heegaard splittings, \3-manifold groups have many more finite
quotients than finitely presented groups in general.  Recent work of
the first author and Dylan Thurston shows a similar sharp divergence
behavior with respect to fibering over the circle, where here a group
``fibers'' if is an algebraic mapping torus \cite{DunfieldDThurston}.
Surprisingly, the \3-manifolds studied there were much \emph{less}
likely to fiber than similar finitely presented groups.

\subsection{Outline of contents}

In Section~\ref{models}, we discuss several different models of random
3-manifolds and some of their basic properties.  In
Section~\ref{sec-random-balanced-pres}, we discuss groups coming from
random balanced presentations, both as a warm up for the \3-manifold
case and to provide a point of comparison.  We compute the
probabilities that such random groups have a particular abelian or
simple quotient.  Our results about random balanced presentation fit
most naturally into the context of profinite groups as we discuss in
Section~\ref{sec-profinite}.  Also in Section~\ref{sec-profinite}, we
define a profinite generalization of random Heegaard splittings.  In
Section~\ref{sec-quotient-3-mfld}, we discuss some reasons why
\3-manifolds should have many finite quotients, working from a more
naive heuristic point of view than in later sections.

The remainder of the paper,
Sections~\ref{heegaard_cover_intro}--\ref{cover-homology}, focuses on
the specific model of random Heegaard splittings and on the finite
covers of the corresponding \3-manifolds.
Section~\ref{heegaard_cover_intro} contains as much of the picture as
we could develop for an arbitrary finite covering group $Q$.  In
Section~\ref{simple_quotients}, we give much more detailed results in
the case when $Q$ is simple.  Similarly,
Section~\ref{homology_Heegaard} is devoted to the case when $Q$ is
abelian.  Finally, Section~\ref{cover-homology} discusses the homology
of a cover of a random \3-manifold.

\begin{acknowledgements}  

Part of the work for this paper was done while the authors were at
Harvard and UC Davis, respectively.  The first author was partially
supported by an NSF Postdoctoral Fellowship, NSF grant DMS-0405491,
and a Sloan Fellowship.  The second author was partially supported by
NSF grants DMS-9704135 and DMS-0072540. The first author would like to
thank Michael Aschbacher for helpful conversations.  We also thank the
referee for a very careful reading of this paper and numerous helpful
comments thereon.

\end{acknowledgements}  

\section{Models of random 3-manifolds} \label{models}

In this section we give several different models of random
\3-manifolds, and outline some elementary properties about them.  In
each case, the idea is to filter \3-manifolds in such a way so that number
of \3-manifolds with bounded complexity is finite.

\subsection{Random triangulations}  
Arguably the most natural notion of a random \3-manifold comes from
filtering by the number of tetrahedra in a minimal triangulation.
To sidestep the difficult problem of determining minimal
triangulations, we can make the triangulations themselves the basic
objects.  Let $\T_3(n)$ be the set of oriented triangulations of
closed \3-manifolds with $n$ tetrahedra.  Here, a triangulation is
just an assemblage of 3-simplices with their faces glued in pairs, and
need not be a simplicial complex in the classical sense.  In the
probabilistic setting, we are interested in the properties of the
manifolds in $\T_3(n)$ as $n$ tends to infinity.  For instance, does
the probability that $M \in \T_3(n)$ is hyperbolic go to $1$ as $n \to
\infty$?  Unfortunately, it seems very difficult to prove anything about
the manifolds in $\T_3(n)$, or even generate random elements of
$\T_3(n)$ for large $n$.  As we will explain
(Proposition~\ref{most-not-3-manifolds}), the problem is that if we start
with $n$ tetrahedra and glue their faces in pairs we almost never get
a \3-manifold.

\subsection{Random surfaces}
We Will start with the 2-dimensional case, since one gets a good
picture there and it helps explain the problem in dimension 3.  Let
$\T_2(n)$ be the set of oriented triangulations of (not necessarily
connected) surfaces; as with studying random graphs, it is convenient
to make these \emph{labeled} triangulations where each triangle is
assigned a number in $1,2,\ldots,n$ and also has an identification with
the standard 2-simplex.  Since you can not build a surface from an odd
number of triangles, we will always assume that $n$ is even.  Another
point of view on a triangulated surface is to look at the dual
1-skeleton.  This is a trivalent graph with labeled vertices where the
incoming edges at each vertex are labeled by $1,2,3$ according to
the label of the edge that they are dual to.  Conversely, such a labeled
trivalent graph gives a triangulation in $\T_2(n)$.  This
triangulation is unique since there is only one way to glue a pair of
sides on two oriented triangles compatible with the orientations.

We can generate elements of $\T_2(n)$ with the uniform distribution
easily, indeed in time linear in $n$: simply start with $n$ triangles
and pick pairs of sides at random and glue.  Note that the dual
1-skeleton of a random element of $\T_2(n)$ is a \emph{uniformly
  distributed} random trivalent graph with $n$ vertices.  Thus we can
directly apply results about random regular graphs to study properties
of $\T_2(n)$ (see \cite{Wormald1999} for a survey of regular random
graphs).  For instance, it follows that the probability that $\Sigma \in
T_2(n)$ is connected goes to $1$ as $n \to \infty$.

\subsection{Euler characteristic} 

We will discuss further consequences of the structure of the dual
graph later, but first we will explain why the expected genus of a
random surface is close to the maximum possible.  With a slightly
different model, that of gluing sides of an $n$-gon, the genus
distribution is needed to compute the Euler characteristic of the
moduli space of Riemann surfaces of a fixed genus.  For this reason,
it was studied in detail by Harer, Zagier, and Penner
\cite{HarerZagier86} \cite{Penner88} \cite{Zagier95}.

To compute the Euler characteristic of $\Sigma \in \T_2(n)$, we just need
to know the number of vertices $v$ as $\chi(\Sigma) = -n/2 + v$.  So what is
the expected number of vertices?  Take the point of view of randomly
gluing triangles together, and think of how the links of the final
vertices are built up by the gluing process.  We start with $3 n$ link
segments in the corners of the $n$ triangles.  These link segments
have an orientation induced from the orientation on the triangles.  At
each stage, we have some number of arcs and circles built up out of
these segments.  At each gluing of triangles, two pairs of endpoints
of arcs are glued together, respecting the orientations.

If we were not gluing the link arcs two at a time, we would have
exactly the same situation as counting the number of cycles of a
random permutation (see e.g.~\cite[\S X.6(b)]{Feller1}).  We Will first
describe what would happen with this simplification and discuss the
full picture below.  With this simplification, at the $k^{\text{th}}$
arc gluing we have $(3 n - k + 1)$ choices of where to glue the
positive end of the given link arc, and exactly one of these choices
creates a closed link circle.  Thus we expect $1/(3 n-k+1)$ final
vertices to be created per gluing, and the expected total number of
vertices is $\sum_{k = 1}^{3 n} 1/k \approx \log( 3 n) \approx \log(n)$.  Note that this says
that the expected genus is about $n/4 - \log(n)$ where the maximum
genus possible with $n$ triangles is $\lfloor n/4 \rfloor$.  Thus the expected
genus is quite close to the maximal one.  In particular, the
probability that $\Sigma \in \T_2(n)$ has any fixed genus $g$ goes to zero
as $n \to \infty$.

Unfortunately, we do not know how to prove that the expected number of
vertices is $\log(n)$; as this is somewhat tangential, we content
ourselves with the following upper bound, which gives many of the same
qualitative results:
\begin{theorem}
  The expected number of vertices for a random surface $\Sigma \in\T_2(n)$ is at most $(3/2) \log(n) + 6$.
\end{theorem}
\begin{proof}
  First, we will explain where the problem is with the argument we gave
  above.  Call an unglued edge of a triangle \emph{bad} if the two
  link arcs which intersect it are actually the \emph{same} arc.
  Gluing a bad edge creates a link circle if and only if it is glued
  to another bad edge.  Thus, the expected number of circles created
  by such a gluing depends on the number of preexisting bad edges,
  which could conceivably be large.
  
  We can deal with this problem as follows.  At each stage, we always
  pick a good edge as the first edge in the pair to glue, if there is
  one.  Because we allow the first chosen edge to be glued to
  \emph{any} edge, good or bad, every $\Sigma$ is generated by this
  process and we have not changed the distribution on $\T_2(n)$.
  However, we have made counting easier.  Let $G_k$ be the random
  variable which is the number of link circles created at the
  $k^{\text{th}}$ stage by a gluing which contains at least one good
  edge.  Let $B_k$ be the number of bad edges created at the
  $k^{\text{th}}$ gluing.  (Both of these variable are set to $0$ once
  we have exhausted the good edges.)  The number of vertices in the
  final surface is equal to $\sum G_k$ plus half the number of bad edges
  left at the end of the good gluings.  The number of bad edges at the
  end is at most the number created during the whole process.  Thus as
  expectations always add, the expected number of vertices is bounded
  by
  \[
  \sum_{k= 1}^{3 n/2} E(G_k) + \frac{1}{2} \sum_{k=1}^{3 n/2} E(B_k).
 \]
 We claim that if there is a good edge left, then $E(G_k) = 2/(3n - 2k +
 1)$, where the denominator is the number of choices for an edge to glue
 to.  There are two cases to consider, depending on whether the link
 arcs of our chosen good edge have other endpoints on the same edge,
 but the probability is the same in both cases.  Similarly, you can
 see $E(B_k) \leq 2/(3n - 2k + 1)$.  Combining, we get that the expected
 number of vertices is less than $3 \sum_{k=1}^{3 n/2} 1/(2k - 1) \leq  3/2 \log(n) + 6$.
\end{proof}

We conclude with an outline of how to turn the problem of precisely
computing the expected Euler characteristic into a problem about the
character theory of the symmetric group.  Suppose we are to create a
surface from $n$ triangles.  Label the oriented sides of the triangles
by $1,2,\ldots,3 n$.  A pairing of the sides can be thought of as a
fixed-point free involution $\pi$ in the symmetric group $S_{3 n}$.
Label the oriented vertex links of the triangles by saying that such a
vertex link has the same number as the edge which contains its
positive endpoint.  We label each triangle so that the 3-cycle $( 3k -
2, 3k - 1, 3k)$ rotates the edges of the $k^\mathrm{th}$ triangle by
$1/3$ of a turn in the direction of the vertex links.  Let $\sigma = (1 2
3) (4 5 6) (7 8 9) \ldots \in S_{3n}$ be the element which rotates the
sides of each triangle in this way.  Then the vertex link $k$ is glued to the
vertex link $(\sigma \pi) (k)$.  Thus the number of vertices of the surface
corresponding to the gluing permutation $\pi$ is just the number of
cycles of $\sigma \pi$.  Hence if $C$ is the conjugacy class of fixed-point free
involutions in $S_{3 n}$, then the average number of vertices is:
\[
\frac{1}{\num{C}} \sum_{\pi \in C}(\mbox{num cycles in $\sigma \pi$}).
\]
There are different ways of attacking problems of this kind, and an
elementary approach is to use the character theory of the symmetric
group, see \cite{Jackson87} \cite{Zagier95}.  For other approaches, based on
random matrices, see \cite{HarerZagier86} \cite{Penner88}
\cite{ItzyksonZuber90}.
 
\subsection{Local structure} 

An interesting property of random regular graphs is that most
vertices have neighborhoods which are embedded trees.  More precisely,
fix a radius $r$ and let the number of vertices $n$ get large.  Then
with probability approaching 1, the proportion of vertices which have
neighborhoods which are embedded trees of radius $r$ is very near 1.
The distribution of short cycles in a random regular graph is also
understood, and as the following theorem shows, the distribution is
essentially independent of the size of the graphs \cite[Thm
2.5]{Wormald1999} :
\begin{theorem}[Bollob\'as]\label{bollo-thm}
  Consider regular graphs where the vertices have valence $d$.  Let
  $X_{i,n}$ be the random variable which is the number of cycles of
  length $i$ in a random such graph with $n$ vertices.  Then for $i$
  less that some fixed $k$, the $X_{i,n}$ limit as $n \to \infty$ to
  independent Poisson variables with means $\lambda_i = \frac{ (d -1)^i}{2
    i}$.
\end{theorem}
One consequence is that if we fix $r$ and pick $\Sigma \in T_2(n)$ with $n$
large, there is a non-zero (albeit small) probability that the
shortest cycle in the dual 1-skeleton has length $\geq r$. That is, the
``combinatorial injectivity radius'' of a random triangulated surface
is large a positive proportion of the time.

\subsection{Triangulations of 3-manifolds}

Now we return to trying to understand random triangulations in
$\T_3(n)$.  Unlike the surface case, we can not study this question by
studying random gluings of tetrahedra:  

\begin{proposition}\label{most-not-3-manifolds}
  Let $X$ be the cell complex resulting from gluing pairs of faces of
  $n$ tetrahedra at random.  Then the probability that $X$ is a
  \3-manifold goes to $0$ as $n \to \infty$.
\end{proposition} 

\begin{proof}
  The link of a vertex in $X$ is always a surface.  Moreover, $X$ is a
  \3-manifold if and only if the combinatorial link of every vertex is
  a sphere (the ``only if'' direction follows for Euler characteristic
  reasons).  Intuitively, since we are gluing the tetrahedra at random, a
  link surface should be a random surface in the above sense.  If this
  were the case, then the probability that the link is a sphere goes
  to $0$ as $n \to \infty$, and so $X$ would almost never be a manifold.
  However, every time we glue a pair of tetrahedra, we are gluing 3
  pairs of link surface pieces at once in a correlated way.
  
  We will finesse this issue by using the fact that if $X$ is a manifold then
  the average valence of an edge is uniformly bounded, and contrast
  this with the fact that the dual 1-skeleton of $X$ is a random
  $4$-valent graph.  For the first point, Euler's formula implies that
  the average valence of a vertex in a triangulation of $S^2$ is less
  than $6$.  The average valence of an edge in $X$ is equal to the
  average valence of a vertex in the vertex links; thus when $X$ is a
  manifold the average edge valence is less than $6$.  In particular,
  since every edge has positive valence, this implies that at least
  $1/6$ of the edges have valence $\leq6$.  Let $\Gamma$ be the dual
  1-skeleton of $X$.  An edge of valence $k$ in $X$ gives a cycle in
  $X$ of length $k$.  Thus if $X$ is a \3-manifold, the number of
  distinct cycles in $\Gamma$ of length $\leq 6$ is a definite multiple
  of the number of vertices.  But $\Gamma$ is a random 4-valent graph,
  and by Theorem~\ref{bollo-thm}, the distribution of the number of
  cycles of length $\leq 6$ is essentially independent of $n$.  Thus as
  $n \to \infty$, the probability that $X$ is a manifold goes to $0$.
 \end{proof}

\begin{remark}
  The proof just given also shows that the probability that a 4-valent
  graph is the 1-skeleton of \emph{some} triangulation of a
  \3-manifold goes to $0$ as the number of vertices goes to infinity.
\end{remark}

All the properties of random surfaces that we described were
consequences of the fact that the uniform distribution on $\T_2(n)$ was
generated by randomly gluing triangles.  In the \3-manifold case, we
are deprived of this tool, and it seems difficult to say anything at
all about a random element of $\T_3(n)$.  For instance, we do not even
know the expected number of vertices for an $M \in \T_3(n)$, much less
whether we should expect $M$ to be irreducible or hyperbolic.  

Even if one could not say much theoretically, it would be very useful
to be able to generate elements of, say, $\T_3(100)$ with the uniform
distribution, even approximately or heuristically.  It would also be
interesting to understand the complexity of uniformly generating
elements of $\T_3(n)$; perhaps there is simply no polynomial-time
algorithm to do so.  It is interesting to note that while spheres make
up a vanishingly small proportion of $\T_2(n)$ as $n \to \infty$, it is
actually possible to generate triangulations of $S^2$ in time linear
in the number of triangles \cite{PoulalhonSchaeffer2003}.

\subsection{Random Heegaard splittings}

Every closed orientable \3-manifold has a Heegaard splitting, that is,
can be obtained by gluing together the boundaries of two handlebodies.
Considering such descriptions gives us another notion of a random
\3-manifold, and this is the one we focus on in this paper.  The
version that we will use here takes the following point of view, based
on the mapping class group.  Fix a genus-$g$ handlebody $H_g$, and
denote $\partial H_g$ by $\Sigma$.  Let $\M_g$ be the mapping class group of
$\Sigma$.  Given $\phi \in \M_g$, let $N_\phi$ be the closed \3-manifold
obtained by gluing together two copies of $H_g$ via $\phi$.  Fix
generators $T$ for $\M_g$.  A random element $\phi$ of $\M_g$ of
complexity $L$ is defined to be the result of a random walk in the
generators $T$ of length $L$.  Then we define the \emph{manifold of a
  random Heegaard splitting} of genus $g$ and complexity $L$ to be
$N_\phi$, where $\phi$ is a random element of $\M_g$ of complexity $L$.
We are then interested in the properties of such random $N_\phi$ as $L
\to \infty$.  A priori, this might depend on the choice of generators
for $\M_g$.  We will show, however, that certain properties do have
well-defined limits independent of this choice
(Sections~\ref{heegaard_cover_intro}--\ref{cover-homology}).

Random Heegaard splittings are much more tractable than random
triangulations, in part because every random walk
in $\M_g$ actually gives a \3-manifold.  Also, for many problems we
can reduce the question to a 2-dimensional one, that is, a question
about the mapping class group.  The disadvantage is that the need to
fix the Heegaard genus feels artificial from some points of view.  For
instance, it means that the injectivity radius of a hyperbolic
structure on $N_\phi$ is uniformly bounded above \cite{White2002}.

A very natural question is how often does the same 3-manifold appear
as we increase $L$?  For instance, there are arbitrarily long walks
$\phi$ in $\M_g$ for which $N_\phi = S^3$.  Thus you might worry that
some small number of manifolds dominate the distribution, and so our
notion of random is not very meaningful.  However, we will show later
that if $\mathcal{F}$ is any \emph{finite} set of 3-manifolds, then
the probability that $N_\phi \in \mathcal{F}$ goes to $0$ as $L \to
\infty$.  This follows from Corollary~\ref{cor-things-dont-repeat},
which shows that $H_1(N_\phi, \Z)$ is almost always finite, but the
expected size grows with $L$.

The next obvious question is: what is the probability that $N_\phi$ is
hyperbolic?  We believe
\begin{conjecture}\label{conj-hyp}
  As $L \to \infty$ the probability that $N_\phi$ is hyperbolic goes to
  $1$.  Moreover, the expected volume of of $N_\phi$ grows linearly in
  $L$.
\end{conjecture}
One expects that the hyperbolic geometry of $N_\phi$ away from the cores
of the handlebodies should be close to that of a ``model manifold'' of
the type used in the proof of the Ending Lamination Conjecture.
Despite this heuristic picture, a proof of Conjecture~\ref{conj-hyp}
is likely to be quite difficult.  One approach would be to try to show
that the expected distance of the Heegaard splitting defining $N_\phi$
(in the sense of Hempel \cite{Hempel2002}) is greater than
$2$.\footnote{Joseph Maher has recently announced \cite{MaherWalk} a proof that
  the probability of a Heegaard splitting having distance less than a
  fixed $C$ goes to $0$ as $L \to \infty$; this would establish the first
  part of Conjecture~\ref{conj-hyp}, assuming geometrization.}
Namazi's results connecting Heegaard splittings to hyperbolic geometry
are also relevant here \cite{Namazi2005}.
 
Finally, there are other notions of random Heegaard splittings you
could consider.  For instance you could think of specifying a random
Heegaard diagram of complexity $L$ by choosing one uniformly from
among the finite number of such where the number of intersections of
pairs of defining curves is $\leq L$.  

\subsection{Universal link related notions}

There are links $L$ in $S^3$ such that every closed orientable
\3-manifold is a cover of $S^3$ branched over $L$.  One such link is
the figure-8 knot \cite{HildenLozanoMontesinos85}.  Let $K$ be this
knot and $M$ be its exterior.  There are only finitely many branched
covers of $(S^3, K)$ of degree $\leq L$, since such a cover corresponds
to a finite-index subgroup of $\pi_1(M)$.  Thus another notion of
random \3-manifold is to choose uniformly among all conjugacy classes
of subgroups of $\pi_1(M)$ of index $\leq L$ and build the corresponding
manifold.  As with random triangulations, it is unclear if there are even
efficient ways to generate such covers experimentally.  While
enumerating all subgroups of index $\leq L$ is certainly algorithmic
\cite{Sims94}, the number of such subgroups grows super-exponentially
in our case.  This is because $\pi_1(M)$ virtually surjects onto a free
group on two-generators.  So if we wanted to sample branched covers of
index $\leq L$ for large $L$, we would need some way of picking out the
subgroups without enumerating all of them.  With current technology,
it would be difficult to enumerate all subgroups of $\pi_1(M)$ for
indices beyond the low 20s.  For more about \3-manifolds from the
point of view of branched covers of the figure-8 knot, see
\cite{Hempel1990}.

\subsection{ Random knots based notions}

Another notion of random \3-manifold would be to take a Dehn surgery
point of view.  That is, one could take some notion of a random knot
or link in $S^3$ and do Dehn surgery on it, where the Dehn surgery
parameters are confined to some finite range at each stage.  For this,
one would need a good notion of a random knot or link.  One could use
models based on choosing a random braid and either taking the closed
braid or making a bridge diagram.  Or you could look at all planar
diagrams with a fixed number of crossings.  These can be efficiently
generated \cite{PoulalhonSchaeffer2003}. Another reasonable notion is
to build a knot out of a uniformly distributed collection of unit
length sticks stuck end to end (see,
e.g.~\cite{DiaoPippengerSumners94}).  These can be generated
efficiently, but have the disadvantage that they are typically
satellite knots \cite{Jungreis94}.

\section{Random balanced presentations}\label{sec-random-balanced-pres}

For a finite presentation of a group, the \emph{deficiency} is the
difference $g - r$ between the number of generators and the number of
relators.  In the case of a closed \3-manifold group, the natural
presentations coming from cell divisions or Heegaard splitting have
deficiency $0$.  Deficiency $0$ presentations are also called
\emph{balanced}.  If a group has a presentation with positive
deficiency, then it already has positive first betti number, so
deficiency 0 is the borderline case for the virtual positive betti
number property of a finitely presented group.  In this section, we
study groups defined by random presentations of deficiency $0$, and
otherwise ignore the constraints coming from the topology of
\3-manifolds. In particular, we compute the probabilities that they
admit epimorphisms to certain finite groups.  In later sections, we will
contrast these results with those specific to \3-manifold groups.

First let us choose a suitable meaning for a ``random presentation'' by
giving a definition of a random relator.  Consider the free group
$F_g$ on $g$ generators $a_1, \dots, a_g$.  Given an integer $n > 0$,
consider all unreduced words of length $n$ where each letter is either
a generator $a_i^{\pm 1}$ or the identity; there are $(2 g + 1)^n$ such
words.  A \emph{random relator} of length $n$ is such a word selected
at random, with each word equally likely.  If we fix a number $g$ of
generators and number $r$ of relations, a \emph{random presentation}
of complexity $n$ is the group $G= \spandef{F_g}{R_1, \ldots, R_r}$, where
each $R_i$ is a random relator of complexity $n$.  Such random
presentations have been studied extensively by Gromov and others.  In
particular, Gromov showed that the probability  that $G$ is
word-hyperbolic goes to $1$ as $n \to \infty$ \cite{Gromov93,Olshanskii92}.

In the rest of this section, we consider the probabilities that
such random groups have different kinds of finite quotients, focusing
on the case of deficiency $0$.  What we do here fits well into the
context of profinite groups, as we describe later in
Section~\ref{sec-profinite}, and that point of view provides additional
motivation for this section.  

\subsection{Quotients of a fixed type}\label{sec-quo-fixed-type}

Let $Q$ be a finite group.   We want to consider the probability that a
random $g$-generator $r$-relator group $G$ has a epimorphism onto $Q$.
We begin by showing that this probability makes sense, and, in later
subsections, calculate it for certain classes of $Q$.  First, consider a
fixed epimorphism $f \maps F_g \to Q$; what is the probability that $f$
extends to $G$, equivalently that $f(R_i) = 1$ for all $i$?  One way
to think of $R_i$ is as the result of a random walk of length $n$ in
the Cayley graph of $F_g$.  In this random walk, each edge is equally
likely as the next step, and there is a $1/(2 g + 1)$ probability of not moving at
each stage.  The image $f(R_i)$ is thus the result of the analogous
random walk in the Caley graph of $Q$ with respect to the generators
$\{ f(a_i) \}$.  The next lemma says that as $n \to \infty$, the result of
such a random walk on a finite graph is nearly uniformly distributed;
thus the probability that $f(R_i) = 1$ converges to $1/\num{Q}$ as $n
\to \infty$.

\begin{lemma}\label{uniform_measure}
  Let $\Gamma$ be a connected finite graph.  Consider random walks on
  $\Gamma$ with fixed transition probabilities.  Suppose that at each
  vertex the probability of taking any given edge is positive, as is
  the probability that the walk stays at the vertex.  Then the
  distribution of the position of the walk after $n$ steps converges
  to the uniform distribution as $n \to \infty$.
\end{lemma}

The reason for requiring a positive probability for pausing at
each stage is to avoid parity issues, as happens when $\Gamma$ is a
cycle of even length; therefore, we will always include the identity
in the set of generators when constructing a Caley graph.  The lemma
is completely standard, but as we use it repeatedly we include a proof.

\begin{proof}
  Consider the vector space $F$ of functions from the vertices of
  $\Gamma$ to $\C$.  Let $L$ be the linear transformation of $F$ which
  averages a function over the radius one neighborhood of a vertex,
  weighted according to the transition probabilities.  That is, if $f
  \maps \Gamma \to \C$ then
  \[
  L(f)(v) = \sum_{w \in B_1(v)} (\mbox{Probability of transition from $v$ to $w$}) f(w).
  \]
  Let $\delta \in F$ be the characteristic function for the initial point
  of the walk.  Then the probability distribution for the position of the walk
  after $n$ steps is $L^n \delta$.  Note that $L$ is a non-negative linear
  matrix, and that if $n$ is greater than the diameter of $\Gamma$ then
  $L^n$ has strictly positive entries.  Constant functions are
  eigenvectors of $L$, and by the Perron-Frobenius theorem, the
  successive images by $L$ of any non-negative and non-zero function
  converge to this eigenspace.  In particular, $L^n \delta$ converges to
  the uniform measure.
\end{proof}

Now let us use the same idea to show that the probability that a
random $G$ has a $Q$-quotient is well-defined.  More precisely, let
$p(Q, g, r, n)$ be the probability that the group of a $g$-generator,
$r$-relator presentation of complexity $n$ has a $Q$-quotient (note there
are only finitely many such presentations).  Then
\begin{proposition}\label{prop-probs-make-sense}
  The probabilities $p(Q, g, r, n)$ converge as the complexity $n$ of
  the presentation goes to infinity.  Moreover, the distribution of
  the number of quotients also converges.
\end{proposition}
In discussing the distribution of the quotients, it is natural to
consider two epimorphisms to $Q$ as the same if they differ by an
automorphism of $Q$, and so we will adopt this convention in our
counts.  This is equivalent to counting normal subgroups with quotient
$Q$.
\begin{proof}  
  Let $\E$ be the set of epimorphisms of the free group $F_g$ onto
  $Q$, modulo automorphisms of $Q$.  We will fix one representative in
  each equivalence class in $\E$, and so regard the elements as actual
  epimorphisms from $F_g$ to $Q$.  Consider the group $Q^\E$, and let
  $P \maps F_g \to Q^\E$ be the induced homomorphism where the
  $f$-coordinate of $P(w)$ is $f(w)$.  Let $S \leq Q^\E$ be the image of
  $F_g$ under $P$, and let $\Gamma$ be the Caley graph of $S$ with respect
  to the generators $\{P(a_i)\}$.  Suppose $R$ is a word in $F_g$.
  Then the $f \in \E$ which kill $R$ are exactly those where the
  $f$-coordinate of $P(R)$ is 1.  If $R$ is a random relator with high
  complexity, then by Lemma~\ref{uniform_measure} the element $P(R)$
  in $S$ is nearly uniformly distributed in $S$.  Thus the probability
  that some $f \in \E$ kills $R$ is approximately the ratio
  \[
   \alpha = \frac{\num{ \setdef{ s \in S}{\mbox{ $s_f = 1$ for some $f$}}}}{\num{S}}.
    \]
    As the relators are chosen independently, the probabilities $p(Q,
    g, r, n)$ converge\footnote{ Not in published version: Michael
      Bush kindly points out that the limiting probability $p(Q, g, r,
      \infty)$ is not usually $\alpha^r$ as we incorrectly claim here.  The
      problem is that for there to be a homomorphism of the resulting
      group, the $r$ relators must all map under $P$ to elements which
      are the identity at the \emph{same} component $f$.  However, the
      $p(Q,g,r,n)$ still have a limit, and the correct expression for
      $p(Q, g, r, \infty)$ can be computed by counting the corresponding
      subset of $S^r$.  The lemma itself is correct, and the erroneous
      formula for $p(Q,q,r,\infty)$ is not used elsewhere in this
      paper.}  to $\alpha^r$ as $n \to \infty$.  Similarly, the probability
    that we have a fixed number $k$ of $Q$-quotients converges for
    each $k$.
\end{proof}  

We will call the limiting probability $p(Q, g, r)$.  As we saw, it
only depends on the finite sets $\E$ and $S$, so we turn now to
understanding them.  First, homomorphisms from $F_g$ to $Q$ which are
not necessarily onto are parameterized by $Q^g$.  The set $\E$ is the
quotient of the subset of $Q^g$ consisting of $g$-tuples which
generate $Q$, under the diagonal action of $\Aut(Q)$.  As $\Aut(Q)$
acts freely on this proper subset, we get that $\num{\E} <
\num{Q}^g/\num{\Aut(Q)}$.  This over-estimate will actually be close to
$\num{\E}$ if $g$ is large; as $g \to \infty$ the proportion of $g$-tuples
in $Q^g$ which do not generate goes to $0$.  Understanding $S$ in
general is complicated as it is typically not all of $Q^\E$.  However,
it is easy to compute the expected (average) number of $Q$-quotients of
such random $G$.  Note that for fixed $f \in \E$ the probability that
the $i^\mathrm{th}$ relator is in the kernel of $f$ is $1/\num{Q}$.
As the relators are chosen independently, the probability that $f$
extends to our random group with $r$ relators is $1/\num{Q}^r$.  Thus
the expected number of such quotients coming from $f$ is
$1/\num{Q}^r$; as expectations add, the expected number of
$Q$-quotients for $G$ is $\num{\E}/\num{Q}^r$.  For any non-negative
integer-valued random variable, the chance it is positive is less than
or equal to its expectation.  Thus $p(Q, g, r) \leq \num{\E}/\num{Q}^r <
\num{Q}^{g - r}/\num{\Aut(Q)}$.  Now, we are most interested in balanced
presentations, and this gives:
\begin{theorem}\label{thm-gen-bal-bound}
  Let $Q$ be a finite group.  The probability $p(Q, g, g)$ that a
  random $g$\nobreakdash-\hspace{0pt}generator balanced group has a
  epimorphism to $Q$ is $< 1/\num{\Aut(Q)}$.
\end{theorem}
Now the number of finite groups with $\num{\Aut(Q)}$ bounded is finite
\cite{LedermannNeumann1956a}, and so the theorem implies that $p(Q, g,
g) \to 0$ as $\num{Q} \to \infty$.  Thus the larger $Q$ is, the less
likely a random balanced $G$ is to have $Q$ as a quotient.  In the
rest of this section we refine our picture for the classes of abelian and
simple groups.

\subsection{Non-abelian quotients}

We start with the case of a non-abelian simple group $Q$, where we
develop a complete picture.  As in Section~\ref{sec-quo-fixed-type},
consider the set $\E$ of epimorphisms from $F_g \to Q$, modulo
$\Aut(Q)$. Most collections of $g > 1$ elements of a finite simple
group $Q$ generate it, especially if $g > 2$ or if $Q$ is not too
small.  To get a rough idea of the probability that a random
collection of $g$ elements of $Q$ generates $Q$, consider the contrary
hypothesis.  If the elements fail to generate, then there is some
maximal subgroup $H$ of $Q$ that contains them all.  For a particular
$H$, the chance that $g$ elements lie in $H$ is
$1/\left[Q:H\right]^g$.  The sum over all maximal subgroups gives an
upper bound for the proportion that do not generate, substantially
less than 1 in all but a few small cases.  These upper bounds with
$g=2$ in a few of the small cases are $(A_5, .53)$, $(A_6,.57)$,
$(A_7,.35)$, $(A_8,.34)$, $ ( A_9,.18)$, $(\PSL{2}{\F_7},.41)$,
$(\PSL{2}{\F_8}, .17)$, $(\PSL{2}{\F_9}, .57)$, $(\PSL{2}{\F_{11}},
.28)$ and $(\PSL{2}{\F_{13}}, .11)$.  As the size of the simple group
gets larger, the probability of $2$ elements generating goes to 1; see
the references in \cite[\S1.1]{Pak2001}.

The automorphism group of a non-abelian finite simple group contains
$Q$ itself as the group of inner automorphisms; the quotient group
is the outer automorphism group, which is generally rather small.  The
upper bound we gave earlier is thus $\num{\E} \leq
\num{Q}^{g}/\num{\Aut(Q)} = \num{Q^{g-1}}/\num{\Out(Q)}$.  The
preceding paragraph indicates that this bound is actually quite
accurate except for small $Q$ and $g$.  As in
Section~\ref{sec-quo-fixed-type}, we now have that the expected number
of $Q$ quotients of a $g$-generator balanced group is
$\num{\E}/\num{Q}^g$ and that this is a bound on the probability
$p(Q,g,g)$ for having a $Q$-quotient.  Thus we have
\begin{equation}\label{simp_quo}
p(Q,g,g) \leq \frac{\num{\E}}{\num{Q}^g} < \frac{1}{\num{Q} \num{\Out(Q)}} \leq \frac{1}{\num{Q}}.
\end{equation}

In order to compute $p(Q,g,g)$ exactly, we need to understand the
image $S$ of the induced product map $F_g \to Q^\E$ used in the
proof of Proposition~\ref{prop-probs-make-sense}.  In this case $S$ is
actually all of $Q^\E$:
\begin{lemma}[\cite{Hall1936}]
\label{lem-simple-surjects-to-product}
  Consider epimorphisms $f_i \maps F_g \to Q_i$, where each $Q_i$ is a
  non-abelian finite simple group.  Suppose no pair $(f_i, f_j)$ are
  equivalent under an isomorphism of $Q_i$ to $Q_j$.  Then the
  product map $F_g \to \prod Q_i$ is surjective.
\end{lemma}
It is important in this lemma that $Q_i$ be non-abelian.  For
instance, if we take $Q = \Z/2$, then $\num{\E} = 2^g - 1$ and so
$Q^\E$ has $2^{2^g - 1}$ elements.  In contrast, the image of $F_g \to
(\Z/2)^\E$ is generated by $g$ elements and thus has size at most
$2^g$.
\begin{proof}
  This lemma was first proved by P.~Hall \cite{Hall1936}.  As it is
  crucial for us, and unfamiliar to most topologists, we include a
  proof.
  
  We begin with case $n=2$.  Suppose $F_g \to Q_1 \times Q_2$ is not
  surjective.  Let $S$ be the image; we will show that $S$ is the
  graph of an isomorphism between $Q_1$ and $Q_2$ compatible with the
  $f_i$.  Consider the projection $\pi \maps S \to Q_1$, and let
  $\overline{Q}_2$ denote the subgroup $\{1 \} \times Q_2$ in $Q_1 \times Q_2$.
  Let $K$ be the kernel of $\pi$, that is $K = S \cap \overline{Q}_2$.
  Note that the conjugation action of a $q \in \overline{Q}_2$ on
  $\overline{Q}_2$ can also be induced by conjugating by some $s \in
  S$, as the projection $S \to Q_2$ is onto; therefore, as $K$ is
  normal in $S$ it must be normal in $\overline{Q}_2$.  As $Q_2$ is
  simple, $K$ must either be $1$ or $\overline{Q}_2$.  In the latter
  case, $S$ contains $\overline{Q}_2$ which implies $S = Q_1 \times Q_2$ as
  $\pi$ is onto.  Thus $K = 1$ and $\pi$ is an isomorphism.  Similarly,
  the projection $\pi' \maps S \to Q_2$ is an isomorphism.  Thus $f_1$
  and $f_2$ are equivalent under the isomorphism $\pi' \circ \pi^{-1}$.
  
  The $n = 2$ case did not use that the $Q_i$ are non-abelian.  That
  hypothesis is used in the form of
  \begin{claim}
    Let $N \unlhd Q_1 \times \cdots \times Q_k$ be a normal subgroup, where the $Q_i$
    are non-abelian simple groups.  Then $N$ is a direct product of a
    subset of the factors.
  \end{claim}
  As before, let $\overline{Q}_i$ denote the copy of $Q_i$ in the
  product.  To see the claim, first observe that $N \cap \overline{Q}_i$
  is either $1$ or all of $\overline{Q}_i$.  If the latter case, mod
  out by $\overline{Q}_i$ to get a case with smaller $k$.  So we can
  assume $N \cap \overline{Q}_i = 1$ for each $i$.  But then $[ N,
  \overline{Q}_i] \leq N \cap \overline{Q}_i$ as both subgroups are
  normal, and so $[N, \overline{Q}_i] = 1$.  But then $N$ is central,
  and thus trivial, proving the claim.  
  
  To conclude the proof of the lemma, choose the smallest $n$ such
  that $F_g \to Q_1 \times \cdots \times Q_n$ is not surjective, and let $S$ be the
  image.  Then as in the $n=2$ case, the projection $\pi \maps S \to Q_1
  \times \cdots \times Q_{n-1}$ is an isomorphism.  Let $\alpha \maps Q_1 \times \cdots \times Q_{n-1}
  \to Q_n$ be the composition of $\pi^{-1}$ with projection onto $Q_n$.
  By the claim, the kernel $N$ of $\alpha$ is a direct product of some of the
  factors.  After reordering, we can assume $N = 1 \times Q_2 \times \cdots \times Q_{n-1}$.
  But then the map $F_g \to Q_1 \times Q_n$ is not surjective, and we are
  back in the $n=2$ case.
\end{proof}

We saw above that the limiting probability of getting exactly $k$
quotients with group $Q$ is simply the density of $s \in S$ with
exactly $k$ of the coordinates equal to $1$.  Thus as $S = Q^\E$ the
limiting distribution is the binomial distribution:
\begin{equation}\label{eq-prob-dist-balanced}
\left\{\num{\mbox{$Q$-quotients}} = k \right\}= \binom{n}{k} p^k \left( 1 - p \right)^{n - k}
\end{equation}
where $p = 1/\num{Q}^g$ and $n = \num{\E}$.
This binomial distribution is well-approximated by the Poisson
distribution.  Recall that the Poisson distribution with mean $\mu > 0$
is a probability distribution on $\Z_{\geq 0}$ where $k$ has probability
$\frac{\mu^k}{ k!} e^{-\mu}$.  Roughly, the Poisson distribution
describes the number $k$ of occurrences of a preferred outcome in a
large ensemble of events where, individually, the outcome is rare and
independent, but in aggregate the expected number of occurrences is
$\mu > 0$.  For instance, it is the limit of the binomial distribution we
have here, if $\mu = \num{\E}/\num{Q}^g$ is kept constant and $n =
\num{\E}\to \infty$.  The difference between (\ref{eq-prob-dist-balanced})
and the Poisson distribution is usually negligible even in small
cases.  For instance, $p(A_5, 2, 2)$ is $0.0052646\ldots$ whereas the
Poisson approximation is $1 - e^{-\mu} \approx 0.0052638$.  Summarizing, we have:
\begin{theorem}\label{thm-simple-balanced}
  Let $Q$ be a non-abelian finite simple group.  Let $n$ be the number
  of epimorphisms from the free group $F_g$ to $Q$ (modulo $\Aut(Q)$),
  and let $\mu = n/\num{Q}^g$.  The probability that a random
  $g$-generator balanced group has a $Q$-quotient is
  \[
  p(Q, g, g) = 1 - \left(1 -   \num{Q}^{-g}\right)^{\num{Q}^g \mu} \approx 1 - e^{-\mu}.
  \]
  and the distribution of the number of quotients is nearly
  Poisson with mean $\mu$.  
  
  Moreover, as $g$ goes to infinity $\mu \to 1/\num{\Aut{Q}}$, and the
  distributions limit to the Poisson distribution with mean
  $1/\num{\Aut{Q}}$.
\end{theorem}

We end this subsection with Table~\ref{expected_table} which
summarizes the situation for the first few finite simple groups.  As
you can see, all the probabilities are very low; we will see that
this is not the case for \3-manifold groups.  

\begin{table}
\begin{center}
\begin{tabular}{rrrrrr}
\textbf{Quotient}  & \textbf{Order} &  \textbf{Gen pairs} & \textbf{Out} &  \textbf{Exp 2-gen}& \textbf{Exp n-gen} \\
$A_5$     &    60 & 19 & 2&  .005278 & .008333\\
$\PSL{2}{\F_7}$  &   168 & 57 & 2& .002020 & .002976\\
$A_6$     &   360 & 53 & 4& .000409 & .000694\\
$\PSL{2}{\F_8}$  &   504 & 142 & 3& .000559 & .000661\\
$\PSL{2}{\F_{11}}$ &   660 & 254 & 2& .000583 & .000758\\
$\PSL{2}{\F_{13}}$ &  1092 & 495 & 2& .000415 & .000458 \\
$\PSL{2}{\F_{17}}$ &  2448 & 1132 & 2&.000189 & .000204 \\
$A_7$     &  2520 & 916  &2& .000144  & .000198\\
$\PSL{2}{\F_{19}}$ &  3420 & 1570 &2& .000134 & .000146\\
$\PSL{2}{\F_{16}}$ &  4080 & 939 &4& .000056  & .000061\\
$\PSL{3}{\F_3}$  &  5616 & 2424  & 2& .000077  & .000089\\
$U_3(\F_3)$  &  6048 &  2784& 2&  .000076& .000083\\
$\PSL{2}{\F_{23}}$ &  6072 & 2881 & 2 &.000078 & .000082\\
$\PSL{2}{\F_{25}}$ &  7800 & 1822 & 4 & .000030& .000032\\
$M_{11}$    &  7920 &  6478& 1& .000103& .000126\\
$\PSL{2}{\F_{27}}$ &  9828 & 1572 & 6&.000016  & .000017\\
$\PSL{2}{\F_{29}}$ & 12180 & 5825 & 2& .000039 & .000041\\
$\PSL{2}{\F_{31}}$ & 14880 &7135  & 2& .000032 & .000034\\
$A_8$     & 20160 & 7448 & 2&  .000018& .000024\\
$\PSL{3}{\F_4}$  & 20160 &1452  & 12 &000004  & .000004\\
$\PSL{2}{\F_{37}}$ & 25308 & 12291 & 2& 000019 & .000020\\
$U_4(\F_2)$  & 25820 &11505 & 2&  .000017& .000019\\
$Sz(\F_8)$   & 29120 &  9534& 3& .000011 & .000011\\
$\PSL{2}{\F_{32}}$ & 32736 &6330  & 5&  .000006& .000006\\
 & & &  & & \\
\end{tabular}
\end{center}

\caption{
This table gives values and bounds for the expected number of
epimorphisms from a random deficiency 0 group to a finite simple group
$Q$.  The simple group $Q$ is listed in the first column. The second
column is the order of $Q$, and the third column \textbf{Gen pairs} is
the number of pairs of elements that generate $Q$, up to automorphisms
of $Q$. Column \textbf{Out} gives the order of the outer automorphism
group of $Q$.  Column \textbf{Exp 2-gen} is the expected number of
epimorphisms to $Q$ among groups with 2 generators and 2 relators.
Column \textbf{Exp n-gen} is $1/|\Aut(Q)|$, which is an upper bound
for the expectation independent of the number of  generators, the limit 
of the expectation as the number $n$ of generators goes to infinity, and a
good approximation to the expectation when $n > 2$.
 }
\label{expected_table}
\end{table}

\subsection{Probability of some simple quotient}\label{subsec-any-simp-quo}

Now, let us consider the more global question: What are the chances
that a finitely presented group of deficiency 0 admits an epimorphism
to \emph{some} non-abelian finite simple group?  First consider a
finite collection $\mathcal{C}$ of simple groups. Let $G$ be a random
$g$-generator balanced group with complexity $n$.  For a fixed $Q \in
\mathcal{C}$, Theorem~\ref{thm-gen-bal-bound} implies that the
probability of $G$ having a $Q$-quotient is $< 1/\num{\Aut(Q)} \leq
1/\num{Q}$, as long as $n$ is large enough.  As there are finitely
many $Q$, we get that for large $n$ the probability that $G$ has a
$Q$-quotient for some $Q$ in $\mathcal{C}$ is less than
\begin{equation}\label{eq-joint-simple-est}
 \sum_{Q \in \mathcal{C}}  1/\num{\Aut(Q)} \leq  \sum_{Q \in \mathcal{C}}  1/\num{Q}.
\end{equation}
In fact, quotients for different $Q$ are independent (this follows
from Lemma~\ref{lem-simple-surjects-to-product}, just as in the proof
that $S = Q^\E$ in the context of Theorem~\ref{thm-simple-balanced}).
Therefore, we could replace (\ref{eq-joint-simple-est}) by $1 -
\prod\left(1 - 1/\num{\Aut(Q)}\right)$, but the former will do for us
here. 

If we were to formally carry out this calculation for the collection
of all finite simple groups, we would get that the probability of
having a non-abelian simple quotient is less than $\sum 1/\num{\Aut(Q)}
\leq  \sum 1/|Q|$, where the sum is over all such groups.  By the
Classification of Finite Simple Groups, eventually nearly all
non-abelian simple groups up to a given size are of the form
$\PSL{2}{\F_q}$.  As $\PSL{2}{\F_q}$ has size $(q-1)q(q+1)/2$ for odd
$q$ and twice that for even $q$, the sum of $1/|\Aut(Q)|$ over all
non-abelian finite simple groups $Q$ is finite.  It is not even very
large: approximately $0.015$.  However, this does not give a proof
that many random balanced groups have no non-abelian
simple quotients; for a fixed group $G$, the relators certainly do not
map to uniformly distributed random elements of $Q$ as $|Q| \to \infty$.
For one thing, a relator of length $R$ is confined to the ball of
radius $R$ in the Cayley graph of $Q$, and this ball has fewer than
$((2g)^R -1)/(2g-1)$ elements.  If the relators were uniformly
distributed in these balls, then the probability of an epimorphism to
$Q$ would be bounded below (although very small), so one would expect
there to eventually be an epimorphism of $G$ to some phenomenally
large simple group $Q$.  But this argument is also invalid, since
among all finite simple groups with a choice of a sequence of $g$
generators, there are only finitely many isomorphism classes of balls
of radius $R$, so we have only finitely many chances to find an
epimorphism.

The real question, which appears to be a difficult issue, is how many
different isomorphism classes of balls of radius $R$ exist among all
non-abelian finite simple groups.  It seems reasonable to us that most
of these balls fall into patterns with relatively few new variations
when $\log \num{Q}$ is large compared to $R$. A good estimate of this
sort could imply that most groups of deficiency $0$ have only finitely
many epimorphisms to finite simple groups, often no such
homomorphisms.  Since a random group in our sense is word hyperbolic,
this would imply that there are word hyperbolic groups that are not
residually finite.

In any case, for a typical deficiency $0$ group that has no quotients
among the first few non-abelian simple groups, it is clear that if it
has any such quotient, the index must be so astronomically large as to
be far beyond brute force computation.  From the calculation above, a
random balanced presentation with 3 or more generators has about 1.5\%
probability to admit an epimorphism to a non-abelian simple group of
manageable size, and a 2-generator group has about 1.3\% probability.
To test our thinking, we made 1000 random 2-generator presentations
with statistics similar to the census manifolds used in our paper
\cite{DunfieldThurston:experiments}, and computed all epimorphisms to
the first few non-abelian simple groups.  Only 15 of these groups had
any such quotients, and only 4 had more than one such quotient.  This
fits reasonably well with the estimate of 1.3\% above.

\subsection{Abelianization}\label{sec-balanced-abelian}

Let us begin by looking at the abelianization of a random balanced
group from a different, more global, perspective.  The abelianization
is the quotient of $\Z^g$ by the subgroup generated by the
abelianization of each relator. In other words, we can make a matrix
$M_R$ whose columns correspond to the relators so that the $(i,j)$
entry is the exponent sum of the occurrences of generator $g_i$ in
relator $R_j$.  If the abelianization is infinite, the determinant of
$M_R$ is 0, otherwise the determinant of $M_R$ is the order of the
abelianization.  For a random relator $R_j$ of length $n$, the
corresponding column is just the result of a suitable random walk in
the integer lattice $\Z^g$.  Individual entries can also be thought of
as generated by random walks, and like all 1-dimensional walks their
absolute value is proportional to $\sqrt{n}$. Thus the typical
determinant grows large as $r$ grows large---indeed, it grows as
$n^{g/2}$ (to see this rigorously, note that the distribution of
$(1/\sqrt{n}) M_R$ converges to that of matrices with independent
Gaussian entries).

However, for our purposes it is more important to determine the
probability that a random presentation admits a finite-sheeted
covering of a given type, as we did in the previous subsection.  Any
finite abelian group is the product of its $p$-Sylow subgroups, so we
focus in on just one prime.  We Will think about this from the point of
view of the $p$-adic integers $\Z_p$.  As rational integers which are
coprime to $p$ have inverses in $\Z_p$, the $p$-Sylow subgroup of the
cokernel of our matrix $M_R$ is the same as the quotient of $\Z_p^g$
by the $\Z_p$-submodule generated by the columns of $M_R$.  We are
interested in asymptotics as the complexity $n$ of our presentation
goes to infinity, and so we want to understand the limiting
distribution of the $M_R$.  More precisely, let $m_n$ be the
probability measure on elements of $F_g$ coming from random walks of
length $n$.  This gives us a measure on the set of balanced
presentations with $g$-generators, with finite support.  Then one has
\begin{lemma}
  The push-forward of the measure $m_n$ to the space of $p$-adic $g \times
  g$ matrices converges weakly to the uniform distribution, i.e.~in the
  limit the entries are elements of $\Z_p$ chosen uniformly and
  independently with respect to Haar measure.
\end{lemma}
\begin{proof}
  The Haar measure on $\Z_p$ can be understood by thinking of $\Z_p$
  as the inverse limit of $\Z/p^k\Z$.  Showing that the weak limit is
  Haar measure is tantamount to checking that, for each $k$, the
  distribution of the entries modulo $p^k$ converges to the
  uniform distribution as $n \to \infty$.  The mod $p^k$ abelianization of
  a random relator is the same as going for a random walk in the Caley
  graph of $(\Z/p^k \Z)^g$.  As always, that distribution becomes the
  uniform one as $n \to \infty$, proving the lemma.
\end{proof}
This $p$-adic point of view is equivalent to considering $(\Z/p^k
\Z)^g$ modulo the subgroup generated by a random sequence of $g$
elements, and looking at the limiting distribution of quotient groups
as $k$ goes to infinity; however, it has the advantage of giving us a
concrete limiting object to calculate with.

First, let us compute the distribution for the orders of the $p$-Sylow
subgroups.  This is just the largest power of $p$ which divides
$\det(M_R)$; if $| \cdot |_p$ denotes the $p$-adic norm, this is the same
as $1/| \det(M_R)|_p$.  Thus, we need to understand the distribution
of $\det(M_R)$, where $M_R$ is a $g \times g$ matrix with entries in $\Z_p$,
chosen uniformly.  The easy case is when $g = 1$, for then $\det(M_R)$
is just uniformly distributed.  As an element in $\Z/p^k \Z$ has a
$p^{-k}$ chance of being $0$, an element in $\Z_p$ has a $p^{-k}$
chance of being in $p^k \Z_p$.  Thus the chance that $z \in \Z_p$ has
$|z|_p = p^{-k}$ is $c_k = p^{-k} ( 1 - 1/p)$.  A useful way to encode
the sequence $\{ c_k \}$ is to use a generating function:
\[
\sum_{k=0}^\infty c_k t^k = \frac{p - 1}{p - t}.  
\]
In the general case we will show:
\begin{proposition}\label{prop-homology-random-group}
  Let $d_k$ be the asymptotic probability that the order of the
  $p$-Sylow subgroup of the abelian group defined by a random
  $g$-generator balanced presentation is $p^k$.  The generating
  function for the sequence $\{d_k\}$ is
\[ \frac{ (p-1)(p^2-1) \dots (p^g-1)}{(p-t)(p^2-t) \dots (p^g-t)}. \]
\end{proposition}
Except for small primes $p$, this is close to the distribution for the
case $g=1$.  Thus, except for small $p$, the probability that the
$p$-Sylow subgroup is non-trivial is close to $1/p$.  This is the same
as asking that the group surject onto $\Z/p$, and by
Theorem~\ref{thm-gen-bal-bound} we already knew this probability was
$< 1/\num{\Aut{Q}} = 1/(p-1)$.  Thus in this case the general
estimate is close to correct.  We now prove the proposition.
\begin{proof}
  A vector in $\Z_p^g$ has probability $1/p^{kg}$ to be in
  $p^k\Z_p^g$.  This tells us the distribution of the maximal $p$-adic
  norms of an element in the first column of $M_R$: the probability
  that the maximal $p$-adic norm equals $1/p^k$ is a geometric
  progression, with generating function $(p^g-1)/(p^g-t)$.  If the
  first column equals $p^k W$ where $k$ is as large as possible, then
  $\Z_p^g /\left\langle W \right\rangle$ is isomorphic to $\Z_p^{g-1}$.
  Moreover, the remaining columns map to independent random elements
  of this module.  Thus $| \det(M_R) |_p$ is the product of $p^{-k}$
  with $| \det( N ) |_p$ where $N$ is a random $(g - 1) \times (g-1)$
  matrix.  Therefore we can get the generating function for $|
  \det(M_R) |_p$ by multiplying the generating functions for these two
  things.  Inducting on $g$ completes the proof.
\end{proof}
\begin{remark}\label{remark-p-independence}
It is worth noting that the proof shows that the Sylow subgroups for
distinct primes $p$ and $q$ are independent, essentially because the
quotient maps from $\Z$ to $\Z/p$ and $\Z/q$ induces a surjection $\Z
\to \Z/p \times \Z/q$; thus a random walk in $\Z$ pushes forward to the
(nearly) uniform distribution on $ \Z/p \times \Z/q$.
\end{remark}

Now, we will delve further and determine the typical isomorphism type
for the $p$-part of the homology.  One way to describe the isomorphism
class of an abelian $p$-group $A$ is to specify the sequence of ranks
of $\rho_i(A) $ of $p^i A / (p^{i+1} A)$.  For example, the group
$(\Z/p)^2 \oplus \Z/p^2 \oplus \Z/p^5$ corresponds to $4, 2, 1, 1, 1, 0$, with
all subsequent terms also $0$.  Introducing a variable $t_k$ to
denote an instance of $(\Z/p)^k$, then an isomorphism class
corresponds to a monomial in the $t_k$; the example corresponds to
$t_4 t_2 t_1^3$.  With this notation, there is a fairly nice and
straightforward computation for the power series in $t_1, \dots, t_g$
whose coefficients give the asymptotic probability that a
$g$-generator, $g$-relator group has the particular isomorphism type
of $p$-Sylow subgroup of its abelianization; this series is a rational
function $\mathit{AFP}_g$. This is a bit of a digression for studying
\3-manifolds, so we will content ourselves with stating the formulae
for 1, 2, and 3-generator groups:
\[ \mathit{AFP}_1 = \frac{p-1}{p-t_1} \]
\[ \mathit{AFP}_2 = \frac{(-1 + p)^2(1 + p)(p^2 + t_1)}{(p - t_1)(p^4 - t_2)} \]
\[ \mathit{AFP}_3 = \frac{(-1 + p)^3(1 + 2p + 2p^2 + p^3)(p^8 + p^5t_1 + p^6t_1 + 
  p^2t_2 + p^3t_2 + t_1t_2)}{(p - t_1)(p^4 - t_2)(p^9 - t_3)} . \] For
instance, the 2-Sylow subgroup of the abelianization of a random
2-generator 2-relator group has probability 3/8 to be trivial, 9/32 to
be $\Z/2$, 9/64 to be $\Z/4$, 3/128 to be $(\Z/2)^2$, etc.
Independently, the 3-Sylow subgroup has probability 16/27 to be
trivial, 64/243 to be $\Z/3$, 64/729 to be $\Z/9$, 16/2187 to be
$(\Z/3)^2$, and so on.  To see how to compute $\mathit{AFP}_g$, note
that in our case, where $A$ is the quotient of $\Z_p^g$ by the
subgroup generated by a random sequence of $g$ elements, the
probability that $\rho_0(A) = k$ is the probability that $g$ elements of
$(\Z/p)^g$ generate a subgroup of rank $g-k$.  Similarly, when
$\rho_i(A) = h$, the conditional probability that $\rho_{i+1} = k$ is the
probability that a random sequence of $h$ elements of $(\Z/p)^h$
generate a subgroup of rank $h-k$.

\section{The profinite point of view}\label{sec-profinite}

In the last section, when we studied the finite quotients of a
``typical'' balanced group, we worked with asymptotic probabilities
$p(Q,g,g)$, which were limits of finite probabilities as the size of
the presentation increases.  In the case of abelian groups, we saw
that these probabilities could be thought of as probabilities on a
certain $p$-adic object, where the notion of probability came from the
natural Haar measure (Section~\ref{sec-balanced-abelian}).  In this
section, we explain how this picture holds true in general by
considering random quotients of profinite free-groups; this helps
clarify why we got well-defined probabilities such as $p(Q,g,g)$.  At
the end, we discuss a natural analog of a Heegaard splitting in the
profinite context.

\subsection{Profinite completions}
We begin with a brief sketch of the theory of profinite groups and
completions (for more, see e.g.~\cite{Wilson1998}, \cite{RibesZalesskii},
 and \cite{DixonEtAl1999}).  Let $G$ be a finitely generated group.  The
\emph{profinite completion} $\overline G$ of $G$ is a compact
topological group defined as the inverse limit of the system of all
finite quotients of $G$.  (Note that whenever $Q_1$ and $Q_2$ are any
two finite quotients, both quotients factor through the image of $G$
in the product map to $Q_1 \times Q_2$, so the set of finite quotients does
form an inverse system.)  If $G$ has only finitely many finite
quotients, then $\overline G$ is a finite group (possibly
trivial).  Otherwise, $\overline G$ has the topology of a Cantor set,
whose stages of refinement give particular finite quotients.  The
natural map $G \to \overline G$ is injective if and only if $G$ is
residually finite.  To reconstruct the finite quotients of $G$, take
small open and closed neighborhoods $V$ of the identity in $\overline
G$, form the subgroup $W$ generated by $V$, and then pass to the
intersection of the finitely many conjugates of $W$ to obtain an open
and closed neighborhood $X$ which is a normal subgroup.  The quotient
$\overline G / X$ is a finite group, and the finite quotients obtained
in this way from any neighborhood basis of $1$ are cofinal among all
finite quotients of $G$.  In general, a \emph{profinite group} is any
compact topological group that has a neighborhood basis of the
identity consisting of open and closed subgroups.  Equivalently, a
profinite group is a group that is the inverse limit of finite groups.
Since a profinite group $\overline G$ is a compact topological group,
it has a unique bi-invariant probability measure, its Haar measure.
This measure is the inverse limit of the counting measures on its
finite quotients.  Thus, any property of elements of $\overline G$ has
a well-defined probability (provided the set of such elements is
measurable).

\subsection{Profinite presentations}
In the profinite context, a finitely presented group is the following.
Consider the profinite completion $\overline{F}_g$ of the free group
on $g$ generators.  Given a finite set of elements $\left\{ R_1, \ldots ,
  R_r \right\}$ of $\overline{F}_g$, let $K$ be the topological
closure of the normal subgroup they generate.  The quotient
topological group $\overline{G} = \overline{F}_g/K$ is the group of
the \emph{profinite presentation} with $g$ generators and relations
$\{ R_i \}$.  Now focus on the set $\mathcal{B}_g$ of all
$g$-generator balanced profinite presentations, which is just the
product of $g$ copies of $\overline{F}_g$, one for each relator.  As
such, it has a natural probability measure coming from the product of
Haar measures on each factor; equivalently, we are thinking of each
relator as being chosen independently at random.  Thus we can
talk about the probability that $\overline{G} \in \mathcal{B}_g$ has
some particular property.  In the case of the property of having a
epimorphism to a finite group $Q$, this is really the same question we
encountered before:
\begin{theorem}\label{thm-profinite-prop-same}
  Let $Q$ be a finite group.  Let $\overline{G}$ be the group defined
  by a randomly chosen $g$-generator balanced profinite presentation.
  Then the probability that $\overline{G}$ has a epimorphism to $Q$ is
  $p(Q,g,g)$.
\end{theorem}
The quickest way to see this would be to repeat the proof of
Proposition~\ref{prop-probs-make-sense} in this context, and see that
one gets the same answer.  We will phrase it a little differently to
make clear why we get the same answer---after all, the set of regular
(non-profinite) presentations has measure 0 in $\mathcal{B}_g$, and so
it is hardly given that asymptotic probabilities of regular presentations
are the same as the corresponding probabilities for profinite
presentations.

Consider random walks on $F_g$ of length $n$, and let $m_n$ be the
probability measure on $F_g$ given by the endpoints of such walks.  We
can also think of $m_n$ as a measure on $\overline{F}_g$.  Then we have:
\begin{lemma}\label{lem-relator-measures-converge-weakly}
  The measures $m_n$ converge weakly to Haar measure on
  $\overline{F}_g$.
\end{lemma}
\begin{proof}
  On a totally disconnected space such as $\overline{F}_g$, locally
  constant functions are uniformly dense among continuous functions.
  So it suffices to check that for a locally constant function $f
  \maps \overline{F}_g \to \R$, the integrals of $f$ with respect to
  $m_n$ converge to the integral of $f$ with respect to Haar measure.
  Since every locally constant function on $\overline{F}_g$ is the
  pullback from a function on some finite quotient, this lemma follows
  from Lemma~\ref{uniform_measure}.
\end{proof}
If $S \subset \mathcal{B}_g$ is both open and closed, then its
characteristic function is continuous.  Hence, if we look at regular
(non-profinite) presentations
Lemma~\ref{lem-relator-measures-converge-weakly} implies
\[
\lim_{n \to \infty} P\setdef{G \in S}{\mbox{$G$ a random balanced group of complexity $n$}} = \mu( S ),
\]
where $\mu$ is the natural measure on $\mathcal{B}_g$.  For instance,
the property of having an epimorphism to a fixed finite group $Q$ is
both open and closed; thus Theorem~\ref{thm-profinite-prop-same}
follows from Lemma~\ref{lem-relator-measures-converge-weakly}.
\comment{Should the first half of the preceding sentence have more
  justification? -N}

Passing to random profinite presentations makes it possible to
estimate the probability that a group has no non-abelian simple
quotients at all.  For regular presentations, we weren't able to show
this probability was positive, but the formal calculation in
Section~\ref{subsec-any-simp-quo} actually applies in the profinite
context.  In particular, the subset of $\mathcal{B}_g$ consisting of
groups which surject onto $Q$ has measure less than $1/\num{\Aut{Q}}$.
As $\sum 1/\num{\Aut{Q}}$ is finite and indeed about $0.015$ we have:
\begin{theorem}
  Let $\overline{G}$ be the group defined by a random $g$-generator
  profinite balanced presentation.  Then with probability $1$, the
  group $\overline{G}$ has only finitely many non-abelian finite
  simple quotients.  If $g \geq 3$, the probability that $\overline{G}$
  has no such quotients is about 98.5\%; if $g=2$, about 98.7\%.  
\end{theorem}

The abelian quotients of a random balanced $\overline{G}$ can be
understood directly from Section~\ref{sec-balanced-abelian}.  Usually,
but not almost always, the abelianization $A$ of a random balanced
profinitely presented group is the inverse limit of cyclic groups.
This is equivalent to the condition that there is no prime $p$ such
that $A$ admits a continuous epimorphism to $\Z/p \times \Z/p$.  Most of
the exceptions are for $p=2$, with most of the remaining exceptions
for $p=3$; the probability for the existence of such a homomorphism is
only about $1/p^4$ for larger $p$.  Among balanced profinitely
presented groups with 1 through 5 generators, the probabilities that
all finite abelian quotients are cyclic are about $1.0, 0.924, 0.885,
0.865, 0.856$. The limiting value for a large number of generators is
about $0.847$.

\comment{ The profinite point of view is nice because it gives us a
  concrete object to think about instead of thinking about asymptotics.
  There's another limit that we take both in
  Section~\ref{sec-random-balanced-pres} and the case of 3-manifold
  groups, namely letting the number of generators go to infinity; at
  least in the cases of certain kinds of quotients one also gets a
  limiting behavior.  I wonder if there is a limiting object there as
  well --- it might give a proof that there is a limiting behavior for
  groups which are not simple (at the moment, we only know how to show
  that the expectations converge, but not the distributions).  -N}

\subsection{Profinite generalizations of 3-manifold groups}

In this subsection, we define a class of profinite groups that
includes the profinite completions of all \3-manifold groups; this
class comes with a natural probability measure.  While we will not
make direct reference to these ideas elsewhere in this paper, they
provide a natural context for the results of
Sections~\ref{heegaard_cover_intro}-\ref{homology_Heegaard}, just as
 groups with balanced profinite presentations do for the
results of Section~\ref{sec-random-balanced-pres}. 

Consider a Heegaard diagram of a \3-manifold, and let $S_g$ be the
fundamental group of the Heegaard surface.  Looking at the fundamental
groups of the two handlebodies gives us a diagram of groups
\[ 
F_g \gets S_g \to F_g.
\]
There is a corresponding diagram of profinite completions:
\[
\overline F_g \gets \overline S_g \to \overline F_g.
\]
The profinite completion of the fundamental group of the \3-manifold
is the quotient of $\overline S_g$ by the topological closure $K$ of
the normal closure of the kernel of the two homomorphisms. 

Since $S_g$ is finitely generated and residually finite, there is a
neighborhood basis for $1$ in $\overline  S_g$ that consists of
\emph{invariant} subgroups of finite index, that is, subgroups
invariant under all automorphisms. For any invariant subgroup, the
mapping class group $\M_g$ acts as an automorphism of the quotient
group. Therefore $\M_g$ is also residually finite, and furthermore, the
action of $\M_g$ on $\overline S_g$ extends to a continuous
action of $\overline \M_g$.

\begin{rmk} 
  These actions are not necessarily faithful, in particular it is not
  for $g=1$.  The torus case is the first case of the congruence
  subgroup problem: does every finite index subgroup of $\SL{k}{\Z}$
  contain a principal congruence subgroup?  (A principal congruence
  subgroup is the kernel of a reduction mod $n$ to $\SL{k}{\Z/n\Z}$.)
  For $\SL{2}{\Z}$ the answer is no, basically since $\SL{2}{\Z}$ is
  virtually a free group and thus it is easy to find quotients which
  are simple groups not isomorphic to $\PSL{2}{\F_p}$.  (Tangentially,
  the answer to the congruence subgroup problem is yes for $k \geq 3$.)
  It is unknown if the action of $\overline \M_g$ is faithful in genus
  greater than 1.
\end{rmk}

This picture gives us some justification in considering
\emph{profinite Heegaard diagrams } $\overline F_g \gets \overline S_g \to
\overline F_g$ which are limits of diagrams of actual \3-manifolds; in
other words, they are obtained by gluing two copies of the standard
map $\overline S_g \to \overline F_g$ by an element of $\overline
\M_g$.  Associated with such a diagram is a locally compact, totally
disconnected topological group, which we will refer to as a
\emph{profinitefold group}: the quotient of $\overline S_g$ by the
smallest normal, closed subgroup $K$ containing the kernels of the two
homomorphisms to $\overline F_g$.  (This construction is special to
dimension 3, so we will not bother with a dimension indicator such as
``3-profinitefold group''.)

Let $\T_g$ be the subgroup of $\M_g$ consisting of homeomorphisms of
the surface that extend to homeomorphisms of the handlebody.  Two
elements $f_1, f_2 \in \M_g$ define equivalent Heegaard diagrams if
$\T_g \backslash f_1 / \T_g = \T_g \backslash f_2 / \T_g $.  Similarly, it makes sense
to define two Heegaard profinite diagrams to be equivalent if the
gluing automorphisms are in the same double coset in $\overline \T_g
\backslash \overline \M_g / \overline \T_g$.  Haar measure on $\overline \M_g$
pushes forward to a measure on this double coset space. \comment{I
  want to say that each double coset has measure 0, i.e. the quotient
  measure has no atoms. In particular, there are uncountably many
  equivalence classes. One proof would use the independent
  probabilities that any finite quotient occurs.}  This gives a
probability measure on the set of profinitefold groups, which we will
use to make sense of statements about random profinitefold groups.

The first homology of any finite sheeted cover of any irreducible
\3-manifold $M$ can be reconstructed from $\overline {\pi_1(M)}$: if
$\Gamma$ is the fundamental group of this finite sheeted cover, the
profinite completion of the abelianization of $G$ is the same as the
abelianization of $\overline \Gamma$, which is the corresponding subgroup
of finite index in $\overline{\pi_1(M)}$.  In particular, the
abelianization of $\Gamma$ is infinite if and only if the abelianization
of $\overline \Gamma$ admits a continuous epimorphism to $\overline \Z$.

In the profinite context, first consider a group $\overline{G}$
defined by a random balanced profinite presentation.  We claim that
with probability $1$, $\overline{G}$ has no continuous homomorphism to
$\overline \Z$.  In this context,
Proposition~\ref{prop-homology-random-group} says that the probability
that $\overline{G}$ has a continuous homomorphism to $\Z/p$ is about
$1/p$.  As these probabilities are independent as we vary $p$ (see
Remark~\ref{remark-p-independence}), the probability that we have one
to all $\Z/p$ is zero. 

Turning now to the case of profinitefold groups, the natural analog of
the virtually positive betti conjecture is
\begin{question}
\label{profinite virtually betti}
Does the profinitefold group $\overline G$ defined by almost every
profinite Heegaard diagram have a subgroup of finite index with a
continuous epimorphism to $\overline \Z$?  \end{question}

It is too weak a condition merely to require that $\overline G$ have
an infinite abelianization.  In fact, the abelianization of almost
every profinitefold group is indeed infinite, because the first
homology group of a \3-manifold is typically a finite group that is
large if the manifold is complicated. 

There are uncountably many isomorphism classes of Heegaard profinite
diagrams up to isomorphism, so the countable set coming from profinite
completions of actual Heegaard diagrams forms a set of measure $0$.
Thus, Question \ref{profinite virtually betti} and the question of
whether all \3-manifolds with infinite fundamental group have
virtually positive betti number does not appear to have any easy logical
implication one way or the other --- the divergence between them
involves different orders of taking limits.  Nevertheless, they are
intuitively and heuristically connected, and so it would be quite
interesting to settle Question~\ref{profinite virtually
  betti}.

\section{Quotients of 3-manifold groups}\label{sec-quotient-3-mfld}

Group presentations coming from Heegaard splittings of \3-manifolds
differ substantially from random deficiency-0 presentations because
the relators, rather than being generic elements in the free group,
are given by a $g$-tuple of disjoint simple closed curves on a genus
$g$ handlebody.  Indeed, \3-manifold presentations are a vanishingly
small proportion of all deficiency-0 presentations since the number of
simple closed curves with word length $R$ grows polynomially in $R$
rather than exponentially.  Geometrically, the curves' embeddedness
forces the words to be far from independent, and typically there are
many repeating syllables at varying scales (for a graphical
illustration of this, see \cite[Fig.~1.5]{DunfieldDThurston}).

In this section, we try to explain why these geometric properties
force there to be more finite quotients than for a general
deficiency-0 group.  Later we will examine this question from the
point of view of random Heegaard splittings
(Sections~\ref{heegaard_cover_intro}--\ref{cover-homology}), but in
this section we take a more naive heuristic point of view.  In
particular, we try to explain why a given quotient $f \maps F_g \to Q$ is
much more likely to extend over the last 2 relators. 

\subsection{Last relator}
One reason to expect \3-manifold groups to have more finite quotients
than random deficiency-0 presentations has to do with the last
relator. To describe this topologically, if we attach 2-handles to the
handlebody along $g-1$ of the curves, we obtain a \3-manifold $M$
whose boundary is a torus; the remaining curve is a simple closed
curve on the torus.  Thus, if an epimorphism $F_g \to Q$ satisfies the
first $g-1$ relators, the remaining relator is restricted to an
abelian subgroup $A$ of $S$ that can be generated by at most 2
elements.  Assuming that the distribution in $A$ is nearly uniform,
this suggests that there is approximately a $1/\num{A}$ chance that
the last relator is satisfied, as compared to a $1/\num{Q}$ chance for
a general relator.  Actually, the situation is more complicated
because the last relator is a \emph{simple} closed curve on the torus
$\partial M$.

Consider a torus $T$ and a finite quotient $f \maps \pi_1(T) \to A$.
Simple closed curves on $T$ correspond to primitive elements of
$\pi_1(T) = \Z^2$, and so we are interested in the probability that a
primitive element lies in the kernel of $f$.  If $A$ is cyclic of
order $a$, then one can change basis so that $f$ is the
factor-preserving map $\Z \oplus \Z \to \Z/a \oplus 0$.  Now look at all primitive
lattice points in $\Z^2$ in some large ball; we want to know the
proportion of them which lie in the kernel of $f$.  It turns out that
this is not quite $1/a$, but rather $1/\beta(a)$ where $\beta$ is a
function of the prime decomposition of $a$ given by
\[
\beta( p_1^{k_1} p_2^{k_2} \cdots p_n^{k_3} ) = p_1^{k_1-1}(p_1 + 1)
p_2^{k_2 - 1} (p_2 + 1) \cdots p_n^{k_n - 1} (p_n + 1).
\]
On the other hand, if $A$ is non-cyclic then you can change basis for
$\pi_1(T) = \Z^2$ so that $f$ is the factor-preserving map $\Z \oplus \Z \to
\Z/a \oplus\Z/ ab$; thus every element of the kernel is divisible by $a$,
and so there are \emph{no} primitive elements in the kernel.  

Returning to our original situation, suppose we are attaching the last
of $g$ relators and want to know if a given epimorphism $F_g \to Q$
extends over this final handle.  This leads us to ask: what is the
distribution of possible subgroups $A$ of $Q$ which are the image of
the fundamental group of the remaining torus $T$?  Not all 2-generator
abelian subgroups can occur.  For instance, the image $H_2(A) \to
H_2(Q)$ must be trivial, since the torus $T$ is the boundary of a
\3-manifold and $H_2(T) \to H_2(A)$ is surjective.  This condition
reduces the number of non-cyclic $A$ we need to consider (though it
need not eliminate them completely), which is good since those never
extend to the resulting manifold.

For example, in $A_5$ the subgroups isomorphic to $\Z/2 \oplus \Z/2$ are
eliminated by this criterion, and so the relevant abelian subgroups
are just the cyclic subgroups, which have orders $1$, $2$, $3$ or $5$.
If each type of cyclic subgroup occurs equally often, this would lead
to the guess that the last relator would be satisfied $35\%$ of the
time.  But even if the cyclic group of order 5 occurs much more
frequently than the others, this would still give an estimate that the
last relator would be satisfied 16.7\% of the time, far more than the
1.7\% predicted for a random relator.  If one looks at a random map
from $\Z^2 \to A_5$, then the cyclic group of order 5 is indeed the
most common image, occurring about half the time.

\subsection{Genus 2} \label{genus-2-relator}

There is also a special argument that sometimes applies to the
next-to-last relator, which works in particular for epimorphisms to
$\PSL{2}{\F_q}$.  The surface $\Sigma$ of genus 2 has a special
homeomorphism of order 2, the hyperelliptic symmetry $\tau$, that is
centralized by the entire mapping class group of $\Sigma$.  The
quotient $\Sigma/\tau$ orbifold is a sphere with 6 elliptic points of order
2.  Any simple closed curve on $\Sigma$ can be isotoped to be set-wise
invariant under $\tau$.  If the curve is non-separating, then it is
mapped to itself with reversed orientation.  If we fix a hyperbolic
metric on $\Sigma$ which is invariant under $\tau$, then the
geodesic representative of a non-separating curve passes through
exactly 2 of the 6 fixed points of $\tau$.  

The consequences of this are easiest to describe for the boundary of a
genus 2 handlebody $H$: Any non-separating simple closed curve on
$\partial H$ describes a circular word in the free group $F_2 = \pi_1(H)$
that is the same read backward or forward.  This is because the
hyperelliptic symmetry $\tau$ of $\partial H$ extends over $H$; the induced
action $\tau_* \maps \pi_1(H) \to \pi_1(H)$ sends standard generators $\{a,
b\}$ of $\pi_1(H)$ to their inverses $\{a^{-1}, b^{-1}\}$.  As
mentioned, $\tau$ sends a non-separating curve on $\partial H$ to itself
with reversed orientation.  Thus if $w$ is a word in $\pi_1(H)$ which
is the image of a non-separating curve, we have that $\tau_*(w)$ is
conjugate to $w^{-1}$; if $w$ is regarded as a circular word this is
the same as saying that it is the same read backward or forward.  It
may or may not be possible to conjugate the linear word $w$ to read
the same backward and forward.  This depends on which pair of the 6
fixed points of $\tau$ the curve passes through, as we now explain. Pick
dual discs for our chosen basis of $\pi_1(H)$ which are invariant under
$\tau$; these contain 4 of the fixed points of $\tau$.  Running around the
geodesic representing $w$ looking at intersections with the discs
reads off the word $w$.  If the geodesic goes through one of the
middle 2 fixed points of $\tau$ which are not near the dual disks, then
reading off starting at one of those points results in a linearly
palindromic $w$.  If instead we start reading from a fixed point of
$\tau$ on one of the discs, then we get a $w$ so that $\tau(w) = s w^{-1}
s^{-1}$ where $s$ is one of the generators.  Thus, it is always
possible to conjugate $w$ so that $\tau_*(w) = s w^{-1} s^{-1}$ where $s
\in \{ 1, a^{\pm 1}, b^{\pm 1}\}$.  In this case, we will say that $w$ is
in \emph{standard form}.

In the case of an epimorphism $\pi_1(H) \to \PSL{2}{\F_q}$, we will show
that the involution $\tau_*$ on $\pi_1(H)$ pushes forward to one on
the image group.  We will use this fact to greatly restrict the
possibilities for the image of a non-separating curve in
$\PSL{2}{\F_q}$.  In particular, we will show:
\begin{theorem}~\label{genus_2_image}
  Let $H$ be a handlebody of genus 2, and $\rho \maps \pi_1(H) \to
  \PSL{2}{\F_q}$ be an epimorphism where $q$ is an odd prime power.
  If $w$ is a word in standard form coming from a non-separating
  embedded curve in $\partial H$, then the image of $w$ under $\rho$ lies
  in a subset of $\PSL{2}{\F_q}$ of size at most $(1/2)(q^2 + q + 2)$.
\end{theorem}
For comparison, the order of $\PSL{2}{\F_q}$ is $(1/2)q(q+1)(q-1)$.
Note also there is always a non-separating curve on $\Sigma$ in the kernel
of $\rho$, as there is one in the kernel of $\pi_1(\Sigma) \to \pi_1(H)$.  Thus
from this theorem one would naively expect that a non-separating curve
is about $q$ times more likely to be in the kernel of $\rho$ than a
random word in $\pi_1(H)$.  In the theorem, the subset of
$\PSL{2}{\F_q}$ mentioned depends on the standard form of $w$, more
precisely on the $s$ such that $\tau_*(w) = s w^{-1} s^{-1}$.  If you
prefer a statement which is independent of $s$, just multiply the size
of the subset by $5$ (the number of possibilities for $s$).

Before proving Theorem~\ref{genus_2_image}, let us further contrast the
picture it gives with that of random words in $\pi_1(H)$.  Consider
\3-manifolds $M$ obtained by attaching a single 2-handle to $H$ along
a non-separating curve (these are examples of tunnel-number one
\3-manifolds).  For comparison, look at two-generator, one-relator
groups where the generator is chosen at random.  For such a random
group, we can work out the probability of a $Q = \PSL{2}{\F_q}$
quotient just as we did before; there are $\approx \num{Q}/\num{\Out{Q}}$
epimorphisms from $F_2$ onto $Q$, and each factors over the relator
with probability $1/\num{Q}$.  Thus the number of $Q$-quotients should
be roughly Poisson distributed with mean $1/\num{\Out{Q}}$.  If we
specialize to the case that $q$ is prime, then $\num{\Out{Q}} = 2$ and
so the probability of a $Q$-quotient for the random group is $\approx 1 -
e^{-1/2} \approx 39\%$.  In particular, this probability is essentially
independent of $Q$.  In contrast, Theorem~\ref{genus_2_image} suggests
that the number of quotients of a \3-manifold group should be Poisson
distributed with mean $\approx q/\num{\Out{Q}}$.  In the case where $q$ is
prime, this leads to the probability of a $Q$ cover being $1 -
e^{-q/2}$, which goes to $1$ as $\Out{Q} \to \infty$.  The last column of
Table~\ref{genus-2-simple} gives some data on this, using random curves
coming from our notion of a random Heegaard splitting.  It suggests
that, at least qualitatively, this last prediction of
Theorem~\ref{genus_2_image} really does hold.

To prove Theorem~\ref{genus_2_image}, we first show that $\tau_*$
pushes forward to the image group $\PSL{2}{\F_q}$.  
\begin{lemma}\label{conj-gens-psl2}
  Let $A$ and $B$ be elements of $\PSL{2}{\F_q}$ where $q$ is an odd
  prime power.  Suppose that $A$ and $B$ do not have have a common
  fixed point when acting on $P^1(\F_q)$.  Then there exists an element $T$
  in $\PGL{2}{\F_q}$ of order 2 such that
  \[
  T A T ^{-1} = A^{-1}  \mtext{and} T B T^{-1} = B^{-1}.
  \]
\end{lemma}
\begin{proof}
  In general, the trace of an element in $\PGL{2}{\F_q}$ depends on
  the lift to $\GL{2}{\F_q}$.  However, the elements $T$ in
  $\PGL{2}{\F_q}$ which have order 2 are exactly those where the trace
  of any lift is 0.  Now lift $A$ and $B$ to elements of
  $\SL{2}{\F_q}$, and consider the equations
  \[
  \tr(T) = 0, \tr(TA) = 0, \mtext{and} \tr(TB) = 0
  \]
  where $T$ is a $2 \times 2$ matrix over $\F_q$ (possibly singular).  As
  these are \emph{homogeneous} linear equations, there is a non-zero
  solution, call it $T$.  We claim that $T$ must be non-singular.  If
  not, change to a basis where the first vector spans the kernel of
  $T$, and so
  \[
  T = \begin{pmatrix}0 & t \\ 0 & 0\end{pmatrix}  \mtext{where $t \neq 0$.}
  \]
  But then $\tr(T A) = 0$ forces $A$ to be upper-triangular.  As the
  same is true for $B$, we have a contradiction in that $A$ and $B$
  have a common fixed point in $P^1(\F_q)$.  In $\PGL{2}{\F_q}$, the
  elements $T$, $TA$, and $TB$ have order 2, which proves the lemma.
\end{proof}

Now consider an epimorphism $\rho \maps \pi_1(H) \to \PSL{2}{\F_q}$.  By
the lemma, there exists a $T \in \PGL{2}{\F_q}$ so that the involution
of $\PSL{2}{\F_q}$ induced by $T$ is the push-forward of the
involution $\tau_*$ on $\pi_1(H)$.  So if $w$ is a word in standard form
coming from a non-separating curve, then there is a $U$ of order 2 in
$\PGL{2}{\F_q}$ such that $U \rho(w) U^{-1} = \rho(w)^{-1}$, where $U$ is
one of $T$, $A^{\pm 1}T$, or $B^{\pm 1} T$.  Thus just knowing that $w$
comes from a non-separating curve implies that $\rho(w)$ is sent to its
inverse by an involution of $\PSL{2}{\F_q}$ which is completely
determined by $\rho$ and the symmetry points of $w$.  The next lemma
shows that this restricts $\rho(w)$ to a proper subset of $\PSL{2}{\F_q}$, and
completes the proof of Theorem~\ref{genus_2_image}.
\begin{lemma}
  Let $q$ be an odd prime power, and $U$ an element of $\PGL{2}{\F_q}$
  of order 2.  Then the number of $W \in \PSL{2}{\F_q}$ such that
  $U W U^{-1} = W^{-1}$ is
  \[
  \frac{1}{2} \left( q^2 + (2 \epsilon_- - 1) q + 2 \epsilon_+ \right),
  \]
  where $\epsilon_{\pm}$ is $1$ if $\pm \det U$ is a square in
  $\F_q$ and $0$ otherwise.
\end{lemma}

\begin{proof}
  Requiring that $U W U^{-1} = W^{-1}$ is the same as saying that $U
  W$ has order 1 or 2.  In the former case, $W = U$ and this
  contributes to our count only when $U \in \PSL{2}{\F_q}$, that is,
  when $\epsilon_+ = 1$.  The latter case is the same as counting solutions
  to the equations,
  \[
  \tr( U W) = 0 \mtext{and} \det(W) = 1
  \]
  where here we have lifted everything to $\GL{2}{\F_q}$.  One can
  write out these equations explicitly (one is linear and the other
  quadratic), and it is not to hard to see that the number of
  solutions is equal to $q^2 + (2 \epsilon_- - 1) q$.  Passing to
  $\PSL{2}{}$ from $\SL{2}{}$ reduces the number of such $W$ by half.
  Combining, we get the count claimed for $W$.
\end{proof}

\subsection{Second to last relator}

Having worked out the simpler handlebody case, we return to our
original question about attaching the second to last relator.  Again,
let $\Sigma$ be a surface of genus $2$, and set 
\[
G = \pi_1(\Sigma) =\spandef{a, b, c, d }{ [a,b] = [c,d]},
\]
where our convention is that $[a,b] = aba^{-1}b^{-1}$.  If we choose
the base point for $\pi_1(\Sigma)$ to be the second leftmost fixed point of
$\tau$, then $\tau$ acts on $G$ in the following way
\begin{equation}\label{eq-tau-rules}
\tau(a) = a^{-1}, \tau(b) = b^{-1}, \tau(c) = x  c^{-1} x^{-1}, \mtext{and} \tau(d) = x d^{-1} x^{-1},
\end{equation}
where $x = a^{-1}b^{-1} d c$.  By taking arcs from the base point to
the other fixed points of $\tau$, we see that every non-separating
simple closed curve on $\Sigma$ can be represented by a element $w \in G$
such that $\tau_*(w) = s w s^{-1}$ where $s$ is in $S = \{ 1, a, b, x, x
d\}$.  Such a $w$ is said to be in \emph{standard form}.

Now consider a homomorphism $G \to \PSL{2}{\F_q}$.  The next
theorem gives a criterion for when $\tau$ pushes forward to an
automorphism of $\PSL{2}{\F_q}$.

\begin{theorem}
  Let $G$ be the fundamental group of a surface of genus 2. Let $\tau$
  be the automorphism of $G$ coming from the hyperelliptic involution.
  Consider an irreducible homomorphism $f: G \to \PSL{2}{\F_q}$, where
  $q$ is odd.  If $f$ lifts to a homomorphism into $\SL{2}{\F_q}$,
  then there is an element $T \in \PGL{2}{\F_q}$ such that $T f T^{-1} =
  f \circ \tau $.
\end{theorem}

Now let $M$ be a \3-manifold with boundary $\Sigma$, and suppose $f$ is
the restriction of a homomorphism $\pi_1(M) \to \PSL{2}{\F_q}$.  The
obstruction to lifting $f$ to $\SL{2}{\F_q}$ is an element of $H^2(G,
\Z/2\Z)$, which vanishes as $f$ extends over $M$.  So if $f$ is
irreducible, then it has the above symmetry, and this restricts the
image under $f$ of a non-separating simple closed curve similar to as
before.  If $w \in G$ is a standard form representative for such a
curve, then $f(w)$ must satisfy $T f(w) T^{-1} = f(s) f(w)^{-1}
f(s)^{-1}$, where $s$ is one of the five elements of $S$ and $T \in
\PGL{2}{\F_q}$ is the element inducing the symmetry. If instead $f$ is
reducible, its image lies in a proper subgroup of $\PSL{2}{\F_q}$.

Thus, in either case, the image of $f(w)$ is restricted. So as long as
there is some non-separating simple closed curve in the kernel of $f$,
one would expect that the probability that $f$ extends over the second to
last relator should be much higher than for a random word in a free
group.  However, as in the torus case, there could be situations where
there are no such curves in $\ker{f}$.  Unlike the torus case, we
do not know any examples where this actually occurs.  \comment{ or at
  least I do not --- Bill, do you remember any such examples?  There is
  an example of a epimorphism $f \maps G \to \PGL{2}{\F_8}$, which has
  this symmetry, and where $f$ does not extend over any attached
  handlebody.  However, there is still a non-separating simple closed
  curve in the kernel.  So the issue is really that once you have killed
  this curve the remaining torus maps to a non-cyclic subgroup of
  $\PGL{2}{\F_8}$, i.e. the issue occurs at the last, rather than next
  to last, relator.  (Note that $H_2$ of this group is 0.)  I also
  checked $\F_{11}$ and $\F_{13}$ and did not see any examples there.
}  Now let us prove the theorem.

\begin{proof} 
  Before beginning the proof itself, let us rephrase the question in
  order to make the algebra that follows seem more natural.  Let $X$
  be the hyperbolic orbifold $\Sigma/\tau$.  We have
  \[
  1 \to G \to \pi_1(X) \to \Z/2 \to 1,
  \]
  where $\pi_1(X)$ can be obtained by adding an element $t$ to $G$
  subject to the requirement that conjugating by $t$ induces the
  action of $\tau$ given in (\ref{eq-tau-rules}).  Thus the problem at
  hand is given $f \maps G \to \PSL{2}{\F_q}$, does $f$ extend to a
  homomorphism $\pi_1(X) \to \PGL{2}{\F_q}$?  The separating curve in
  $\Sigma$ representing $\left[ a, b \right]$ maps down to a curve
  $\gamma$ which separates $X$ into two regions containing $3$ cone
  points each; in $\pi_1(X)$ we have $\gamma^2 = \left[a,b\right]$.
  Roughly, Lemma~\ref{conj-gens-psl2} says that we can extend $f$ over
  each half of $X$, and the issue is whether these agree along
  $\gamma$.  It turns out that $f$ lifting to $\SL{2}{\F_q}$
  guarantees this.
  
 Let $\tilde{f}$ denote the lift of $f$ to
  $\SL{2}{\F_q}$, and set $A = \tilde{f}(a)$, $B = \tilde{f}(b)$, etc.
  To prove the theorem, it suffices to find a $T \in \GL{2}{\F_q}$ such
  that
  \begin{equation}\label{eq-traces-zero}
  \tr(T) = \tr(T A) = \tr(T B)= \tr(X^{-1}T) = \tr(X^{-1}T C) = \tr(X^{-1}TD) = 0
  \end{equation}
  as then all of the above elements have order $2$ in $\PGL{2}{\F_q}$
  and so $T f T^{-1} = f \circ \tau $.  The main case is when $f$ is
  irreducible when restricted to both of the subgroups $\left< a,
    b\right>$ and $\left<c, d\right>$, and we begin there.  Think of
  (\ref{eq-traces-zero}) as homogeneous linear equations in the
  entries of a $2 \times 2$ matrix $T$.  By the proof of
  Lemma~\ref{conj-gens-psl2}, any non-zero $T$ which satisfies the
  first $3$ trace conditions is necessarily non-singular.  Thus we
  just need to prove that the dimension of the solution space of
  (\ref{eq-traces-zero}) is positive dimensional.  To check this, we
  are free to enlarge our base field from $\F_q$ to its algebraic
  closure $k$.  Now over $k$, the equations (\ref{eq-traces-zero})
  have a non-zero solution if and only if there is one with $\det(T) =
  1$, since any non-zero determinant has a square root in $k$.  So
  henceforth, we try to solve equations (\ref{eq-traces-zero}) for $T
  \in \SL{2}{k}$.
  
  By Lemma~\ref{conj-gens-psl2}, we can choose $T, U \in \SL{2}{k}$ so
  that
  \[
    \tr(T) = \tr(T A) = \tr(T B)= 0 \mtext{and}  \tr(U) = \tr(U C) = \tr(U D)= 0.
  \]
  Our goal is show $U = X^{-1} T$.  As we are working in $\SL{2}$
  rather than $\PSL{2}$ we have $T^2 = -1$ not $1$, and the same for
  $TA$, $U$, etc.  Thus we have
  \[
  (B A T)^2 = - \left[B, A\right] = - \left[D, C\right] = (D C U)^2. 
  \]
  So $BAT$ and $DCU$ have the same square, which is not $-1$ because
  if $A$ and $B$ commuted we would be in the reducible
  case.  Any $S \neq -1$ in $\SL{2}{k}$ has at most $2$ square roots.
  Thus, after replacing $U$ with $-U$ if necessary, we have $BAT =
  DCU$ which implies $U = X^{-1}T$, as required.  This completes the
  proof of the theorem in the case when $f$ is irreducible when
  restricted to both $\pair{a,b}$ and $\pair{c,d}$.
  
  Now suppose instead $f$ is reducible when restricted to
  $\pair{a,b}$.  This is equivalent to $\tr(ABA^{-1}B^{-1}) = 2$, and
  the same reasoning implies that $f$ must be reducible when
  restricted to $\pair{c,d}$ as well.  Now $ABA^{-1}B^{-1}$ can not be
  parabolic, as then $f$ itself would be reducible.  Therefore, we are
  down to the case where $ABA^{-1}B^{-1}= 1$, i.e. both $\pair{A,B}$
  and $\pair{C,D}$ are abelian.  Now by changing generators of $G$ we
  can assume that $\pair{A, C}$ do not have a common fixed point in
  $P^1$.  By Lemma~\ref{conj-gens-psl2}, there is a $T \in GL_2(\F_q)$
  such that
  \[
 \tr(T) = \tr(T A) = \tr(T C B A) = 0.
  \]
  Now as $B$ commutes with $A$, the above also forces $\tr(T B) = \tr(T
  B A) = \tr(BAT) = 0$.  Rewriting $\tr(T C B A) = 0$ as $\tr(
  BAT \cdot C) = 0$, the commuting of $C$ and $D$ gives $\tr( BAT \cdot
  Y) = 0$ for all $Y \in \pair{C,D}$.  Expanding $X$ in
  (\ref{eq-traces-zero}), it follows that all those equations hold,
  completing the proof.
\end{proof}

\section{Covers of random Heegaard splittings}\label{heegaard_cover_intro}

In this section, we address the following question: Fix a finite group
$Q$ and genus $g$.  Consider a manifold $M$ obtained from a random
genus-$g$ Heegaard splitting.  What is the probability that $M$ has a
cover with covering group $Q$?  What is the distribution of the number
of covers?  In this section, we will begin looking at these by showing
that the answers to both these question are well-defined
(Prop.~\ref{makes_sense}).  We will then give three examples which
illustrate some of the key issues in computing these probabilities.
In later sections, we will compute these probabilities exactly for
abelian groups (Section~\ref{homology_Heegaard}), and we will give a
complete characterization of these probabilities for non-abelian
simple groups in the limit when $g$ is large
(Section~\ref{simple_quotients}).

Fix a genus-$g$ handlebody $H_g$, and denote $\partial H_g$ by $\Sigma$.  Let
$\M_g$ be the mapping class group of $\Sigma$.  Given $\phi \in \M_g$, let
$N_\phi$ be the closed \3-manifold obtained by gluing together two copies
of $H_g$ via $\phi$.  Our notion of a random Heegaard splitting of
genus $g$ is as follows.  Fix generators $T$ for $\M_g$.  A random
element $\phi$ of $\M_g$ of complexity $L$ is defined to be the result
of a random walk in the generators $T$ of length $L$.  Then we define
the \emph{manifold of a random Heegaard splitting} of genus $g$ and
complexity $L$ to be $N_\phi$, where $\phi$ is a random element of
$\M_g$ of complexity $L$.

We begin with the following proposition which shows that the questions we
are interested in make sense in the context of random Heegaard
splittings:
\begin{proposition}\label{makes_sense}
  Fix a Heegaard genus $g$ and a finite group $Q$.  Let $N$ be the
  manifold of a random Heegaard splitting of complexity $L$, and let
  $p(L)$ be the probability that $\pi_1(N)$ has an epimorphism onto
  $Q$.  Then $p(L)$ converges to a limit $p(Q,g)$ as $L$ goes to
  infinity.  Moreover, $p(Q,g)$ is independent of the choice of
  generators for $\M_g$, and the probability distribution of the
  number of epimorphisms also converges.
\end{proposition}

\begin{proof}
  Consider the collection $\A$ of all epimorphisms from $\pi_1(\Sigma)$ to
  $Q$, up to automorphisms of $Q$.  Let $\phi$ be in $\M_g$, and
  consider the associated \3-manifold $N_\phi$.  Identify $\Sigma$ with the
  boundary of the first copy of $H_g$ in $N_\phi$, so that a
  homomorphism $f$ in $\A$ factors through to one of $\pi_1(N_\phi)$ if
  and only if $f$ and $f \circ \phi_*^{-1}$ both extend over $H_g$.
  
  Let $\E \subset \A$ consist of those homomorphisms which do extend over
  $H_g$.  Before proving the full proposition, let us consider the
  simpler question: Given an $f$ in $\E$, what is the probability that
  it extends over the copy of $H_g$ in $N_\phi$?  That is, how often
  does $f \circ \phi_*^{-1}$ also lie in $\E$?  Consider the (left) action
  of $\M_g$ on $\A$ by precomposition: $\phi \cdot f = f \circ \phi_*^{-1}$.
  Thus we are interested in the probability that $\phi \cdot f$ lies in
  $\E$.  Now look at the Caley graph of this action of $\M_g$ on $\A$.
  That is, consider the graph whose vertices are $\A$ and whose edges
  correspond to the action of the chosen generators $T$ of $\M_g$.  If
  we do a random walk in $\M_g$, the corresponding sequence of $\phi_i
  \cdot f$ moves around this graph according to the labels on the edges.
  As in Lemma~\ref{uniform_measure}, the distribution of the $\phi_i \cdot
  f$ converges to the uniform measure on the orbit of $f$ in $\A$.
  Let $C$ be the orbit of $f$ in $\A$.  Since the $\phi_i \cdot f$ are
  nearly uniformly distributed in $C$ for large $i$, we see that the
  probability that $f$ extends to $N_\phi$ converges to $\num{C \cap
  \E}/\num{C}$.  Note that this limit depends only on the orbit $C$ and
  not on the choice of generators for $\M_g$.
  
  Returning to the question of the probability of $N_\phi$ having an
  epimorphism to $Q$, consider the action of $\M_g$ on subsets of
  $\A$.  Again, if we look at the images of $\E$ under $\phi_i$, these are
  nearly uniformly distributed in the orbit of $\E$ in the power set
  of $\A$; thus $p(L)$ converges to the proportion of subsets in the
  orbit of $\E$ which intersect $\E$.  Finally, the distribution of
  the number of epimorphisms also converges, to the corresponding
  finite averages over the orbit of $\E$ under $\M_g$.
\end{proof}

Next, we will give three detailed examples in genus 2 where we
calculate these probabilities exactly.  These illustrate the some
of the main issues and techniques that arise later in Sections
\ref{homology_Heegaard}--\ref{simple_quotients}.

\subsection{Example: $\Z/2$}\label{example-Z2}

For our first genus 2 example, let us begin with $Q = \Z/2$.  In this
case, $\A$ is isomorphic to the non-zero elements of $H^1(\Sigma ; \Z/2)$.
Thus $\num \A = 15$, and similarly $\num \E = 3$.  In order to compute
the probabilities, we need to understand the image of $\M_2$ in the
symmetric group of $\A$.  While the action is transitive, its image is
much smaller than all of $\mathrm{Sym}(\A)$: it is the 4-dimensional
symplectic group over $\F_2$, which we will call $G$.  As we do a
random walk in $\M_2$, the images under $\M_2 \to \mathrm{Sym}(\A_1)$
converge to the uniform distribution on $G$.  Thus $p(\Z/2, 2)$ is the
same as the probability that $(g \cdot \E) \cap \E \neq \emptyset$ for a random $g
\in G$.

First, let us compute the expected number of $\Z/2$-quotients.  To
begin, focus on whether a fixed $f \in \E$ extends to $\pi_1(N_\phi)$.
For $g \in G$, we have $f \in g \cdot \E$ if and only if $g^{-1} \cdot f \in
\E$.  Since the action of $G$ is transitive, $\big\{g^{-1} \cdot
f\big\}_{g \in G}$ is uniformly distributed in $\A$.  Therefore, the
probability that $f \in g \cdot \E$ is $\num \E / \num \A = 1/5$.  Thus
the expected number of $\Z/2$ quotients from $f$ is $1/5$.  As
expectations add, the overall expected number of quotients is
$\num{\E}^2/\num \A = 3/5$.  Computing $p(\Z/2, 2)$ is more
complicated because there may be correlations between different $f$ in
$\E$ extending (see the discussion in the next example).  In this
case, it turns out that averaging over all 720 elements of $G$ gives
that $p(\Z/2, 2) = 7/15$.

Of course, $p(\Z/2,2)$ is the same as the probability that $H^1(N_\phi;
\Z/2) \neq 0$.  In Section~\ref{homology_Heegaard}, we will work out the
distribution of the homology for general $\Z/p$.  Readers uninterested
in non-abelian groups can skip ahead to that section now.

\subsection{Example: $A_5$}

Next, let us consider the smallest non-abelian simple group $Q
= A_5$.  In this case, $\num \A = 2016$ and $\num \E = 19$.  The action of
$\M_2$ on $\A$ has two orbits $\A_1$ and $\A_2$ of size 1440 and 576
respectively, and $\E$ is completely contained in $\A_1$.  In order to
compute $p(A_5,2)$, we only need to understand the action on $\A_1$.
It turns out that this is the full alternating group
$\mathrm{Alt}(\A_1)$.  As before, we get that $p(A_5, 2)$ is the
probability that $(g \cdot \E) \cap \E \neq \emptyset$ for a random $g \in
\mathrm{Alt}(\A_1)$.  Now $\mathrm{Alt}(\A_1)$ acts transitively on
subsets of $\A_1$ of size 19, as $19 \leq \num{ \A_1} - 2$.  Thus $g \cdot \E$
is uniformly distributed over all subsets of $\A_1$ of size $19$.
Hence the probability that $(g \cdot \E) \cap \E \neq \emptyset$ is just the
probability that a randomly chosen subset of 19 elements of $\A_1$
intersects $\E$.  Thus
\[
p(A_5, 2) =  1 - \binom{1421}{19} \Big/  \binom{1440}{19} \approx 22.43\%.
\]
As in the previous example, the expected number of $A_5$ quotients is
$\num{\E}^2/\num{\A_1} \approx 0.2507$.

The case where $\M_g$ acts as the alternating group of an orbit is one
that we will find in general for non-abelian simple groups, so it is
worth discussing here the connection to the Poisson distribution.
Recall that the Poisson distribution with mean $\mu > 0$ is a
probability distribution on $\Z^+$ where $k$ has probability
$\frac{\mu^k}{ k!} e^{-\mu}$.  Roughly, the Poisson distribution
describes the number $k$ of occurrences of a preferred outcome in a
large ensemble of events where, individually, the outcome is rare and
independent, but in aggregate the expected number of occurrences is
$\mu > 0$.  In our context, we have a set $\A$ of size $n$ which
contains a marked subset $\E$ of size $a$; we then pick another subset
$\E'$ of size $a$ and want to know the size of $\E \cap \E'$.  If $n$ is
large, the distribution of $\num{\E \cap \E'}$ is essentially Poisson with
mean $\mu = a^2/n$.  In particular, the probability of at least one
intersection is $1 - e^{-\mu}$.  In the case of $A_5$, this
approximation gives the probability of an $A_5$ cover at 22.17\%.

It is worth mentioning that the alternating group action here makes
it very easy to compute $p(Q, g)$ from just the sizes of the orbits of
$\A$ and $\E$, and that this is not true in general.  For instance,
returning to the previous example, $p(\Z/2, 2) = 7/15 \approx 0.4667$ but
the probability that $\num \E = 3$ items chosen from $\num \A = 15$
things intersects a fixed set of size $3$ is larger: $47/91 \approx
0.5165$.  The reason for the difference is that the action of $\M_2$
is not 3-transitive, and there are positive correlations between
different elements of $\E$ factoring through to $\pi_1(N_\phi)$; in
particular, because the number of $\Z/2$-quotients is $\num {H^1(N_\phi,
\Z/2)} - 1 = 2^n - 1$, if we have more that one $\Z/2$-quotient then we
have 3 of them.

\subsection{Example: $\PSL{2}{\F_{13}}$}\label{example-PSL213}

Let us look at a more complicated example.  Consider $Q =
\PSL{2}{\F_{13}}$, a group of order 1092.  In this case, $\num \A =
623520$ and $\num \E = 495$.  In this case there are four orbits of
sizes 235680, 94080, 278400, and 15360.  Only the first two orbits
intersect $\E$ in subsets of size 307 and 188 respectively; set $\E_i
= \E \cap \A_i$.  The action on the first two orbits $\A_1$ and $\A_2$
are again by the full alternating groups $\mathrm{Alt}(\A_i)$.  Since
the two alternating groups have different orders, by
Lemma~\ref{lem-simple-surjects-to-product} the map $\M_2 \to
\mathrm{Alt}(\A_1) \times \mathrm{Alt}(\A_2)$ is surjective.  Therefore $g
\cdot \E_1$ and $g \cdot \E_2$ are independent of each other, and
\[
p(\PSL{2}{\F_{13}}, 2) = 1 - \prod_{i = 1}^2  \binom{\num{\A_i} - \num{\E_i}}{\num{\E_i}}
 \Big/ \binom{\num{\A_i} }{\num{ \E_i}} \approx 54.02\%.
\]
Further, the expected number of quotients is $\sum_{i=1}^2 \num{
  \E_i}^2/\num{\A_i} \approx 0.7756$.  Although $\E$ is contained in two
orbits, the overall distribution of quotients is still nearly Poisson
since the sum of two independent Poisson variables is also Poisson.

\subsection{Example: Small simple groups}\label{compare-simple}

As we will see in Lemma~\ref{size_of_A_g}, for simple groups the size of
$\A_g$ grows like $\num{Q}^{2g - 2}$.  Thus even for genus 2, it is
difficult to compute the action of $\M_g$ on $\A_g$ when $\num{Q}$ is
large.  However, the proof of Proposition~\ref{makes_sense} suggests a
way to approximate $p(Q,g)$ by looking only at $\num{Q}^{g-1}$
epimorphisms.  Namely, we first compute $\E_g$ and consider
$\phi \in \M_g$ which is evolving by a random walk in our fixed
generators.  The time average of $(\phi \cdot \E_g) \cap \E_g$ will then
converge to $p(Q, g)$.  For genus 2, we did this for the simple groups
of order less than $7000$.  The results are shown in
Table~\ref{genus-2-simple}, and we compare them to the results of our
earlier experiment \cite{DunfieldThurston:experiments}, as well as the
limit of $p(Q, g)$ as $g \to \infty$.

\begin{table}
\begin{tabular}{r|rcc|cc|cc}
$Q$  & Order &  $\num{\Out}$ & $\num{H_2}$& $p(Q,2)$ & Census & $p(Q,\infty)$ & Sing.~Quo \\
\hline
$A_5$     &    60 & 2 &2 & 22.4 & 26.9& 63.2  & 72.09 \\ 
$\PSL{2}{\F_7}$  &   168 &  2  &2 &30.8 & 28.2 & 63.2 & 88.95 \\
$A_6$     &   360 & 4 &6 & 28.9 & 31.4 & 77.7 &  85.04 \\ \hline
$\PSL{2}{\F_8}$  &   504 & 3 & 1& 22.1 &21.7  & 28.3& 89.37 \\
$\PSL{2}{\F_{11}}$ &   660 & 2 & 2& 41.7 &32.8 & 63.2& 98.63 \\ 
$\PSL{2}{\F_{13}}$ &  1092 & 2 &2 & 54.0 &41.1 & 63.2& 99.85 \\ \hline
$\PSL{2}{\F_{17}}$ &  2448 & 2 & 2&56.3 &43.1 & 63.2&  99.97 \\ 
$A_7$     &  2520 & 2 & 6 &  60.1 & 45.8 & 95.0 & 99.86 \\
$\PSL{2}{\F_{19}}$ &3420 &   2  &2 &55.9 & 44.4& 63.2 & 99.99   \\ \hline
$\PSL{2}{\F_{16}}$ &  4080 &  4 & 1&19.3 &18.3 &22.0 & 97.36 \\ 
$\PSL{3}{\F_3}$  &  5616 &  2 & 1&40.5 &28.0 & 39.3&  99.76 \\
$U_3(\F_3)$  &  6048 &  2 & 1& 31.5 & 18.0 & 39.3 & 99.57\\ \hline
$\PSL{2}{\F_{23}}$ &  6072 & 2 & 2 & 58.9  & 47.6 & 63.2 & 99.99 \\
\end{tabular}
\vspace{1cm}
\caption{
This table gives the percentage of genus 2 manifolds which have
particular finite simple quotients, and compares this data to other
samples.  The first 4 columns list the quotient group $Q$ and some of
its basic properties.  The number $p(Q,2)$ is the probability that a
manifold coming from a random Heegaard splitting of genus 2 will have
a cover with covering group $Q$; as discussed in
Section~\ref{compare-simple}, these numbers are actually an
approximation coming from looking at a random walk in $\M_2$ of length
$10^6$.  The Census column is the corresponding probability for those
genus 2 manifolds in the Hodgson-Weeks census that we studied in
\cite[\S5]{DunfieldThurston:experiments}.  The number $p(Q,\infty)$ is
the limit of $p(Q, g)$ as $g \to \infty$;  by
Theorem~\ref{simple_probs} we have  $p(Q, \infty) = 1 -
e^{-\num{H_2}/\num{\Out}}$ .  Finally, we looked at attaching
a single 2-handle to a genus 2 handlebody along a non-separating curve
to give a tunnel number one manifold with torus boundary; as before, 
the attaching curve is chosen by a random walk in $\M_2$.   
The last column records the probability that the resulting manifold 
has a $Q$ cover.  This data
is interesting to compare with the discussion in 
Section~\ref{genus-2-relator}.
}
\label{genus-2-simple}
\end{table}

\subsection{The general picture}

The rest of this section is devoted to what we can say in general
about $p(Q, g)$ for an arbitrary group, with particular emphasis on
the limiting picture as the genus $g$ goes to infinity.  We would like
to say that the probability distributions on the number of covers for
each $g$ converge to a limiting probability distribution as $g \to \infty$;
in particular, this would suggest some robustness in our notion of
random Heegaard splitting and that, perhaps, we should expect to find
this same limiting distribution with other notions of random
\3-manifold.  While we will build up quite a bit of information about
$\A_g, \E_g$ and $\A_g/\M_g$ for large $g$, we do not know how to prove
the existence of a limiting distribution in general.  Instead, we will
show here that the expected number of covers does have a limit as $g
\to \infty$ (Theorem~\ref{expectations_converge}).  In the special
cases of abelian and simple groups, we are able to obtain a complete
asymptotic picture (see Sections~\ref{homology_Heegaard} and
\ref{simple_quotients}).  The case of abelian groups is essentially
independent of the rest of this section, so if that is your primary
interest you can skip ahead to Section~\ref{homology_Heegaard}.

Examples~\ref{example-Z2}--\ref{example-PSL213} illustrate the key
issues that we encounter here.  As usual, let $H_g$ be a handlebody of
genus $g$, $\Sigma_g = \partial H_g$, and $\M_g$ be the mapping class group of
$\Sigma_g$.  We want to compute the probability $p(Q, g)$ that a
\3-manifold associated to a random genus-$g$ Heegaard splitting has a
cover with group $Q$.  As we saw in the examples, what we need to know
is how many epimorphisms $\pi_1(\Sigma_g) \to Q$ there are, and how $\M_g$
acts on them.  As we will see, the answer to the first question is easy
(Section \ref{counting_A_and_E}), and it is the second question that
requires more work.  For the latter question, for a general $Q$ we
will only be able to classify the orbit set $\A_g / \M_g$ (Sections
\ref{homology_classes}--\ref{subsec_orbit_classification}).  The basic
idea is this: given a homomorphism $f \maps \pi_1(\Sigma_g) \to Q$ in $\A_g$
we get a map from $\Sigma_g$ to the classifying space $BQ$, and thus an
element of $H_2(Q, \Z)$.  This homology class is invariant under the
action of $\M_g$ on $\A_g$.  The key is to show that for large $g$
this class is the only invariant of the action of $\M_g$ on $\A_g$
(Theorem~\ref{hom_classify_orbits}).  Finally, this section concludes
with a more group-theoretic point of view on some aspects of this
section as they relate to $\E_g$ (Section~\ref{orbit_of_E}).

\subsection{Counting $\A$ and $\E$}\label{counting_A_and_E}

In this subsection, we determine the number of elements of $\A$ and
$\E$ for a fixed group $Q$ as the genus gets large.  For results where
the genus is fixed but the (typically simple) group $Q$ gets large see
\cite{LiebeckShalev2005a} and \cite{LiebeckShalev2005b}.  If $Q$ is a
group, $Q'$ will denote the commutator subgroup.  The following lemma
computes for us the probability that a random homomorphism $F_{2 g} \to
Q$ factors through $\pi_1(\Sigma_g)$.

\begin{lemma}  
  Let $Q$ be a finite group.  If $(a_1,b_1,\ldots,a_g, b_g) \in Q^{2 g}$ is
  chosen uniformly at random, the probability that $\prod [a_i, b_i] =
  1$ converges to $1/\num{Q'}$ as $g \to \infty$.
\end{lemma}

\begin{proof}
  The set $T = \big\{[a,b]\big\}_{a,b \in Q}$ generates $Q'$.  Thus
  choosing $(a_i,b_i) \in Q^{2 g}$ at random is the same as choosing a
  string of $g$ elements of $T$.  That is, $\prod [a_i, b_i]$ is the
  result of a random walk in $Q'$ with respect to $T$.  As such, it
  converges to the uniform distribution on $Q'$.  Therefore, the
  probability of it being $1$ converges to $1/\num{Q'}$ as $g \to \infty$.
\end{proof}

For non-abelian simple groups, the set $T$ above is very large,
conjecturally all of $Q = Q'$.  In this case, the random walk is on a
(weighted) complete graph with vertices $Q$.  Thus it will converge to
the uniform distribution very quickly.   

Let $\A_g$ be the set of all epimorphisms of $\pi_1(\Sigma_g)$ to $Q$,
modulo automorphisms of $Q$.  

\begin{lemma}\label{size_of_A_g}
  Let $Q$ be a finite group.  Then $\num{\A_g} \sim \num{Q}^{2 g}/
  (\num{Q'} \num{\Aut(Q)}) $ as $g \to \infty$.
\end{lemma}

\begin{proof}
  Consider $\pi_1(\Sigma_g)$ in its standard presentation with $2 g$
  generators.  A $2 g$-tuple $(a_i, b_i) \in Q^{2 g}$ gives a possible
  homomorphism $\pi_1(\Sigma_g) \to Q$.  If the $2 g$-tuple is randomly
  chosen, the probability that the entries generate $Q$ goes to 1 as
  $g \to \infty$.  By the preceding lemma, the probability that it induces
  a homomorphism of $\pi_1(\Sigma_g)$ is $1/\num{Q'}$.  Combining the above, this
  says that the number of elements of $Q^{2g}$ which give epimorphisms
  is asymptotic to $\num{Q}^{2g}/\num{Q'}$.  Since the action of
  $\Aut(Q)$ on epimorphisms is free, we get the claimed formula for
  $\num \A$.
\end{proof}

We also record here the corresponding result about $\E$.

\begin{lemma}\label{size_of_E_g}
  Let $Q$ be a finite group.  Then $\num{\E_g} \sim
  \num{Q}^{g}/\num{\Aut(Q)}$.
\end{lemma}

\subsection{Associated homology classes}\label{homology_classes}

For us, an important invariant of an epimorphism $f \maps \pi(\Sigma_g) \to
Q$ is the associated homology class $c_f$ in $H_2(Q, \Z)$ coming from
the induced map on classifying spaces $\Sigma_g \to BQ$.  In the context of
finite groups, $H_2(Q,\Z)$ is usually called the Schur multiplier, and
it is an important invariant there, especially in the study of simple
groups (see e.g.~\cite{Wiegold82}).  In this section all homology
groups will have $\Z$ coefficients, so we will stop including this in
the notation.  It is not quite true that this homology class is
well-defined for elements of $\A_g$, but we can associate an element
of $H_2(Q)\big/ \Out(Q)$.  The issue here is that $\Aut(Q)$ may act
non-trivially on $H_2(Q)$; after all, elements in $\A_g$ are
equivalence classes of epimorphisms modulo the action of $\Aut(Q)$.
Now inner automorphisms of $Q$ induce self-maps of $BQ$ which are
homotopic to the identity, and so such automorphisms act trivially on
$H_2(Q)$.  Thus associated to an $[f] \in \A_g$ we get a well-defined
homology class $c_{[f]}$ in $H_2(Q) \big/ \Out(Q)$.  For simple
groups, it is often the case that $H_2(Q)$ is trivial or $\Z/2$; the
action of $\Out(Q)$ must be trivial in this case.  The first simple
group where the action is non-trivial is $A_6$, where $H_2(Q) = \Z/6$
and $\Out(Q) \to \Aut(H_2(Q)) = \Z/2$ is surjective.  In general, the
map $\Out(Q) \to \Aut(H_2(Q))$ need not be surjective, as the example
of $\PSL{3}{\F_4}$ shows.

Now the map $\A_g \to H_2(Q)\big/\Out(Q)$ is invariant under the
action of $\M_g$, since acting by a mapping class just re-marks the
surface and so does not change the image of the map $\Sigma_g \to BQ$.  Thus we
get a well-defined map $\A_g/\M_g \to H_2(Q)\big/ \Out(Q)$.  Later,
we will show that this map is a bijection for large $g$.  For this
reason, we are interested in the number of epimorphisms in each
homology class:
\begin{lemma}\label{class-probability}
  Let $Q$ be a finite group.  Fix $[c] \in H_2(Q, \Z)\big/ \Out(Q)$.
  Let $k$ be the number of elements of $H_2(Q, \Z)$ which are in the
  equivalence class $[c]$.   Then the ratio
  \[
  \num{ \setdef{ f \in \A_g}{ c_f = [c]} } \Big/ \num{\A_g}
  \]
  converges to $k \big/\num{H_2(Q,\Z)}$ as $g \to \infty$.  
\end{lemma}
\begin{proof}
  First, we explain how to compute $c_f$ directly using Hopf's
  description of $H_2$ of a group (see e.g.~\cite[\S II.5]{Brown82}).
  Let $G$ be a finitely generated group, and express it as a quotient
  of a free group $G = F/R$.  Then $H_2(G)$ can be naturally
  identified with $(F' \cap R)/[F, R]$.  From this point of
  view, $H_2(G)$ is a subgroup of leftmost term of the exact
  sequence
  \begin{equation}\label{hopfseq}
  1 \to R/[F, R] \to F/[F, R] \to G  \to 1;
  \end{equation}
  note that the leftmost term is central in the middle term.  When $G$
  is finite, $R/[F, R]$ is the direct sum of $(F' \cap R)/[F, R]$ and a
  free abelian group $A$ (see e.g.~\cite{Wiegold82}).  Taking the
  quotient by $A$ we get a exact sequence of finite groups
  \[
  0 \to H_2(G) \to S \to G \to 1,
  \]
  where $S$ is called a Schur cover of $G$.  Here, $H_2(G)$ is
  central in $S$ and contained in $S'$.  The Schur cover need not be
  unique unless $G$ is perfect; in that case $S \cong F'/[F,R]$ and $S$
  is the universal central extension of $G$.
  
  Now, let $f \maps \pi_1(\Sigma_g) \to Q$ be a homomorphism.  Consider
  standard generators $a_i, b_i$ of $\pi_1(\Sigma_g)$; in the Hopf picture,
  $H_2(\Sigma)$ is generated by the standard relator $\prod[a_i, b_i]$.  Fix
  a Schur cover $S$ of $Q$.  Consider an induced map of the
  sequences (\ref{hopfseq}) for $\pi_1(\Sigma_g)$ and $Q$.   Thus
  $c_f$ is given by the following proceedure: Pick lifts $s_i, t_i \in
  S$ of the images $f(a_i), f(b_i)$, and then $c_f = \prod[s_i, t_i]$.
  
  Now we determine the proportion of elements of $\A_g$ which have a
  fixed homology class.  It is easier here to work directly with
  epimorphisms $f \maps \pi_1(\Sigma_g) \to Q$ before quotienting out by
  $\Aut(Q)$.  As in Lemma~\ref{size_of_A_g}, consider a random
  $2g$-tuple $(s_i, t_i)$ of elements of $S$ and let $c = \prod [s_i,
  t_i]$.  The images of $(s_i, t_i)$ in $Q$ are also uniformly
  distributed, and therefore the image of $c$ in $Q$ is nearly
  uniformly distributed in $Q'$ if $g$ is large.  Recall that $H_2(Q)$
  lies in $S'$.  As we know $c$ is nearly uniformly distributed in
  $S'$, if we restrict to those $(s_i,t_i)$ which induce a
  homomorphism $\pi_1(\Sigma_g) \to Q$, the probability that $c$ is some
  particular element of $H_2(Q)$ is essentially $1/\num{H_2(Q)}$.
  Modding out by $\Aut{Q}$ gives the probability claimed in the
  original statement.
\end{proof}

\subsection{Stabilization}\label{sec_stabilization}

Now, we want to discuss our main topological tool for understanding
$\A_g$ when $g$ is large.  Let $f \maps \pi_1(\Sigma_g) \to Q$ be in $\A_g$.
A stabilization of $f$ is an element $f'$ of $\A_{g +h}$ obtained by
viewing $\Sigma_{g + h}$ as $\Sigma_g \# \Sigma_h$ and setting $f'$ to be $f$ on
$\Sigma_g$ and the trivial homomorphism on $\Sigma_h$; that is, $f'$ is the
composition of $f$ with the map on $\pi_1$ induced by the quotient map
$\Sigma_g \# \Sigma_h \to \Sigma_g$.  Note here we are \emph{not} fixing a
particular identification of $\Sigma_{g+h}$ with $\Sigma_g \# \Sigma_h$, so each $f
\in \A_g$ usually gives rise to many elements in $\A_{g + h}$.  So that
stabilization respects the associated classes in $H_2(Q)/\Out{Q}$, we do
require the identification of $\Sigma_{g+h}$ with $\Sigma_g \# \Sigma_h$ to be
orientation preserving.  Looking at it another way, an $f'$ in
$\A_{g+h}$ is the result of an $h$-fold stabilization if there is an
essential subsurface $S$ of $\Sigma_{g + h}$ which is a once-punctured
$\Sigma_h$, and where $f'$ is the trivial homomorphism on $\pi_1(S)$.

First, we will show that for large enough genus every element of
$\A_g$ is a stabilization:
\begin{proposition}\label{every-f-is-a-stabilzation}
  Let $Q$ be a finite group.  If $g > \num{Q}$ then every $f \maps
  \pi_1(\Sigma_g) \to Q$ is a stabilization.
\end{proposition}
Before giving the proof, let us point out an important consequence.  If
$f_1, f_2 \in A_g$ are in the same orbit under $\M_g$, then so are their
stabilizations in $\A_{g+1}$ under the action of $\M_{g+1}$.  Thus we
get a map of orbit sets $\A_g/\M_g \to \A_{g+1}/\M_{g+1}$.  The
proposition says that for $g > \num{Q}$ these maps are surjective; thus
once some threshold is crossed no new orbits appear, though existing
orbits can merge.  But since there are only finitely many orbits, this
merging process must eventually stop.  Thus we have
\begin{corollary}\label{orbits-stabilize}
  Let $Q$ be a finite group.  Then for large enough $g$, the number
  of $\M_g$ orbits in $\A_g$ is constant.  
\end{corollary}
It is worth noting that this argument gives no indication of when
$\num{\A_g/\M_g}$ stabilizes, even though
Proposition~\ref{every-f-is-a-stabilzation} gives an explicit bound.
Later, we will show the stable number of orbits is just
$\num{H_2(Q,\Z)/\Out(Q)}$.  We now prove the proposition.

\begin{proof}
  
  View $\Sigma_g$ as $g$ punctured-torus spokes attached to a central hub
  which is a $g$-times punctured sphere.  Consider standard $2g$
  generators $\{ a_i, b_i\}$ of $\pi_1(\Sigma_g)$, where each $(a_i, b_i)$
  pair are generators of the $i^{\mathrm{th}}$ spoke.  If we orient
  everything symmetrically, any element $w_i = a_1 \cdot a_2 \cdots a_i$ of
  $\pi_1(\Sigma_g)$, with $1 \leq i \leq g$, can be represented by an embedded
  non-separating simple closed curve in $\Sigma_g$.  
  
  Now fix $f \in A_g$.  We need to find an essential punctured torus on
  which $f$ is the trivial homomorphism.  By the pigeon hole
  principle, there exists $i < j$ such that $f(w_i) = f(w_j)$.
  Therefore, $f(a_{i+1} \cdots a_j) = 1$, and so we can find one
  non-separating simple closed curve in the kernel of $f$.  To turn
  this into an entire handle where $f$ is trivial, consider a maximal
  disjoint collection $c_1, \ldots ,c_k$ of such non-separating curves in
  the kernel of $f$.  By changing the basis of $\pi_1(\Sigma)$, we can make
  these $c_i$ be the curves $b_1, \ldots, b_k$ in our preferred basis.
  Then as before, there exists some $w = a_{i+1} \cdots a_j$ in the kernel
  of $f$.  By maximality, the curve $w$ must intersect one of our
  $c_i$, say $c_1$.  As $w$ and $c_1$ intersect in a single point,
  they have a regular neighborhood which is a punctured torus whose
  fundamental group is mapped trivially under $f$.  Thus $f$ is a
  stabilization.
\end{proof}

To complete our understanding of stabilization, we characterize the
stable equivalence classes of epimorphisms. Two epimorphisms $\alpha
\maps \pi_1(\Sigma_g) \to Q$ and $\beta \maps \pi_1(\Sigma_h) \to Q$ are called \emph{stably
  equivalent} if they have a common stabilization.  In particular,
$c_\alpha = c_\beta$ in $H_2(Q)$.  The following theorem of Livingston
\cite{Livingston85} shows that this homological condition is
sufficient as well as necessary.
\begin{theorem}[\cite{Livingston85}]\label{same-class-are-stably-same}
  Let $Q$ be a finite group, and consider two epimorphisms $\alpha \maps
  \pi_1(\Sigma_g) \to Q$ and $\beta \maps \pi_1(\Sigma_h) \to Q$.  Then $\alpha$ and $\beta$
  are stably equivalent if and only if $c_\alpha = c_\beta$ in $H_2(Q, \Z)$.
\end{theorem}
Since this result is crucial to our characterization of $\A_g/\M_g$
and is quick to prove, we include a complete proof, following
\cite{Livingston85}.  Unfortunately, the proof gives no control over
the amount of stabilization required.  One example where
stabilization is needed is $\PSL{2}{\F_{13}}$ where there are multiple
orbits of $\A_2$ which correspond to $0$ in $H_2$ (see
Example~\ref{example-PSL213}). 
For some meta-cyclic groups, Edmonds
showed that no stabilization is required \cite{Edmonds83}.  Zimmermann
\cite{Zimmermann87} gave a quite different, purely algebraic, proof of
this theorem which might be useful for an analysis of the degree of
stabilization required for particular $Q$.  We now give the proof of
the theorem.
\begin{proof}
  Suppose $\alpha$ and $\beta$ are as above, with $c_\alpha = c_\beta$ in $H_2(Q,
  \Z)$.  Now bordism is the same as homology in dimension $2$, so, in
  particular, the two induced maps $\Sigma_g \to BQ$ and $\Sigma_h \to BQ$ can
  be extended over some cobordism between $\Sigma_g$ and $\Sigma_h$
  \cite{ConnerBook} \cite{ConnerFloydBook}.  Thus there exists a
  \3-manifold $M$ with two boundary components $\partial_1M = \Sigma_g$ and
  $\partial_2M = \Sigma_h$, and an epimorphism $f \maps \pi_1(M) \to Q$ which
  restricts to $\pi_1(\partial_iM)$ as $\alpha$ and $\beta$ respectively.  Take a
  relative handle decomposition of $M$ with all the 2-handles added
  after the 1-handles; that is $M$ is $(\partial_1M) \times I$, plus some
  1-handles, plus some 2-handles, plus $(\partial_2M) \times I$.  Let $M_1$ be
  $(\partial_1M) \times I$ and the 1-handles, and let $M_2$ be the 2-handles and
  $(\partial_2 M) \times I$.  Thus $M$ is the union of $M_1$ and $M_2$ along a
  surface $\Sigma$.  To conclude, we will show that the restriction of $f$
  to $\pi_1(\Sigma)$ is a stabilization of both $\alpha$ and $\beta$.
  
  First consider $\alpha$ and $M_1$.  Because $\alpha$ is surjective, we
  can slide the attaching maps of the 1-handles of $M_1$ around so
  that all their cores map trivially under $f$.  This shows $(\Sigma,f)$
  is a stabilization of $\alpha$.  If we flip the handles over, the same
  reasoning applies to $M_2$ and $\beta$.  Thus $\alpha$ and $\beta$ have
  a common stabilization.  
\end{proof}

\subsection{Characterization of $\A_g/\M_g$ for large $g$}\label{subsec_orbit_classification}
Recall from Section~\ref{homology_classes} that we have a natural map
$\A_g/\M_g \to H_2(Q)/\Out(Q)$.  We will now show:
\begin{theorem}\label{hom_classify_orbits}
Let $Q$ be a finite group.   For all large $g$, the map 
\[
\A_g/\M_g \to H_2(Q, \Z)/\Out(Q)
\]
is a bijection.  
\end{theorem}
\begin{proof}
  First, by Corollary~\ref{orbits-stabilize} the size of $\A_g/\M_g$
  is constant for large $g$, and the stabilization maps $\A_g/\M_g \to
  \A_{g+1} / \M_{g+1}$ are bijections.  Moreover, these stabilization
  maps are compatible with the maps to $H_2(Q)$.
  Theorem~\ref{same-class-are-stably-same} implies that any two
  elements of $\A_g$ which represent the same class in
  $H_2(Q)/\Out(Q)$ become the same in some $\A_{g+h}/\M_{g+h}$; thus
  $\A_g/\M_g \to H_2(Q)/\Out(Q)$ is injective for all large $g$.
  Lemma~\ref{class-probability} shows that every class is realized for
  large enough $g$, so the map $\A_g/\M_g \to H_2(Q)/\Out(Q)$ is
  surjective as well (alternatively, this follows because bordism is
  the same as homology in dimension 2).  Thus $\A_g/\M_g \to
  H_2(Q)/\Out(Q)$ is a bijection for all large $g$, as claimed.
\end{proof}

Unfortunately, this theorem and its constituent parts do not seem to
be enough to show the existence of a limiting distribution of the
number of quotients as the genus $g$ goes to infinity.  However, it is 
easy to show that the number of expected quotients does
converge:
\begin{theorem}\label{expectations_converge}
  Let $Q$ be a finite group.  Let $E(Q, g)$ be the expected number of
  covers with covering group $Q$ of a \3-manifold associated to a
  random Heegaard splitting of genus $g$.  Then as the genus $g$ goes
  to infinity
  \[
   E(Q,g)  \to  \frac{\num{Q'} \num{H_2(Q, \Z)}}{\num{\Aut(Q)}}.
  \]
\end{theorem}
\begin{proof}
  Let $g$ be large enough so that $\A_g/\M_g \to H_2(Q)/\Out(Q)$ is a
  bijection.  Let $A_g^0$ be those elements $f$ of $\A_g$ which have
  $c_f = 0$ in $H_2(Q)/\Out(Q)$; by choice of $g$, $\A_g^0$ is a
  single orbit of $\M_g$.  By Lemmas~\ref{size_of_A_g} and
  \ref{class-probability}, we have
  \[
  \num{\A_g^0} \sim \frac{\num{Q}^{2g}}{\num{Q'} \num{\Aut(Q)} \num{H_2(Q, \Z)}}.
  \]
  As usual, let $\E_g$ be those epimorphisms which extend over the
  handlebody.  By Lemma~\ref{size_of_E_g}, we have $\num{\E_g} \sim
  \num{Q}^g/\Aut(Q)$.  Since $\E_g \subset \A_g^0$ and the action of $\M_g$
  on $\A_g^0$ is transitive, as in Example~\ref{example-Z2} we have
  that
  \[
  E(Q, g) = \frac{\num{\E_g}^2}{\num{\A_g}} \sim \frac{\num{Q'}\num{H_2(Q)}}{\num{\Aut(Q)}},
    \]
    as desired.  
\end{proof}

\subsection{Characterization of $\A'_g/\M_g$ for large $g$}\label{subsec_orbit_classification_prime}
In Section~\ref{simple_quotients}, we will need a slight variant of
Theorem~\ref{hom_classify_orbits} and its precursors in
Section~\ref{sec_stabilization}.  For a finite group $Q$ and genus
$g$, define $\A_g'$ to be the set of epimorphisms $\pi_1(\Sigma)$ onto
$Q$; unlike $\A_g$ we do not mod out by the action of $\Aut(Q)$.  As
with $\A_g$, elements in $\A'_g$ have associated homology classes, but
now these are well-defined in $H_2(Q, \Z)$.  In this context, we have
the following analog of Theorem~\ref{hom_classify_orbits}:
\begin{theorem}\label{hom_classify_orbits_prime}
  Let $Q$ be a finite group.   For all large $g$, the map $\A'_g/\M_g \to
  H_2(Q, \Z)$ is a bijection.  
\end{theorem}
The proof of this is the same as for
Theorem~\ref{hom_classify_orbits}, after modifying
Corollary~\ref{orbits-stabilize} and Lemma~\ref{class-probability} by
replacing $\A_g$ by $\A'_g$ and $H_2(Q)/\Out(Q)$ by $H_2(Q)$.  In case
you are wondering why we do not work with $\A'_g$ throughout, the
reason is that the action of $\M_g$ on $\A'_g$ can not be highly
transitive in general, e.g.~the full alternating group of each
component; this is because $\num{\E'_g \cap (\phi \cdot \E'_g)}$ must always
be an integer multiple of $\Aut(Q)$.

\subsection{Orbit of $\E_g$ under $\M_g$}\label{orbit_of_E}

In the case of $\PSL{2}{\F_{13}}$ (Example~\ref{example-PSL213}) we
saw that the action of $\M_2$ on $\E_2$ was not transitive, and there
were two distinct orbits.  However, one consequence of
Theorem~\ref{hom_classify_orbits} is that for any $Q$ the action of
$\M_g$ on $\E_g$ is always transitive when $g$ is large enough.  In
fact, something more is true: the action of the stabilizer of $\E_g$
is transitive for large $g$.  Moreover, when $Q$ is a non-abelian
simple group, the image of the stabilizer of $\E_g$ in $\Sym(\E_g)$
contains $\Alt(\E_g)$.  Both of these facts follow immediately from
the well-developed theory of the action of $\Out(F_g)$ on epimorphisms
$F_g \to Q$.  We outline the needed results below, following the survey
\cite{Pak2001}, to which we refer the reader for further results along
these lines.  The material in this subsection is not strictly
necessary for the rest of this paper, but has two merits.  The first
is that it gives explicit bounds, in terms of $Q$, of when these
phenomena occur.  The second is that the strong result in the case of
simple groups allows us to avoid the use of the Classification of
Simple Groups for some of the results in
Section~\ref{simple_quotients}.

First, note that the stabilizer of $\E_g$ in $\M_g$ contains the
mapping class group $\M(H_g)$ of the handlebody $H_g$.  Let $F_g =
\pi_1(H_g)$ be the free group on $g$ generators.  The mapping class
group $\M(H_g)$ acts on $\E_g = \left\{ f \maps F_g \surjects Q
\right\}\big/\Aut(Q)$ via a homomorphism $\M(H_g) \to \Out(F_g)$; based
on the generators of $\Out(F_g)$ that we list below, it is easy to see
that $\M(H_g) \to \Out(F_g)$ is surjective.  Thus, understanding the
action of $\M(H_g)$ is equivalent to understanding the action of
$\Out(F_g)$ on $\E_g$.  In the language of \cite[\S2.4]{Pak2001}, elements
of $\E_g$ are called $T$-systems, and $\E_g$ is denoted by $\Sigma_g(Q)$.  

Before continuing, recall that $\Out(F_g)$ is generated by three
kinds of elements: replacing a generator by its product with another
generator (on either the left or the right), permuting two generators,
and inverting a generator.  If we think of elements of $\E_g$ as
$g$-tuples of elements of $Q$ (up to $\Aut(Q)$), then these generators
carry over to natural operations on such tuples.  Now let us state
the main fact of this subsection:
\begin{proposition}\label{act_on_E}
  Let $Q$ be a finite group.  Then for all large $g$, the action of
  $\Out(F_g)$ on $\E_g$ is transitive.  Hence the action of $\M_g$ on
  $\E_g  \subset \A_g$ is also transitive for large $g$.  
\end{proposition}
How big does $g$ need to be for the conclusion to hold?  For simple
groups, a conjecture of Wiegold posits that $g \geq 3$ suffices, and
this is known for some classes of such groups
\cite[Thm.~2.5.6]{Pak2001}.  In general, the bound in the proof we will
give is proportional to $\log\num{Q}$.  Now, let us prove the
proposition, following Prop.~2.2.2 of \cite{Pak2001}, which in turn
follows \cite{DiaconisSaloff-Coste1998}.

\begin{proof}
  Let $d$ be the minimal number of generators of $Q$.  Let $\bar{d}$
  be the maximal size of a minimal generating set for $Q$; that is, a
  generating set for which no proper subset generates.  We will show
  that the action of $\Out(F_g)$ is transitive provided that $g \geq d +
  \bar{d}$, which we now assume.  Let $q_1, \ldots, q_d$ be a minimal
  generating set for $Q$.  Consider $f_0 = (q_1,\ldots,q_d, 1, \ldots, 1)$ in
  $\E_g$.  Let $(r_1, \ldots r_g)$ be any element in $\E_g$, which we will
  move to $f_0$ by elements in $\Out(F_g)$.  Since there is some
  subset of the $r_i$ of size $\bar{d}$ which generate $Q$, permute
  the $r_i$ so that $r_{d + 1}, \ldots, r_g$ generates $Q$.  Now by using
  elements of $\Out(F_g)$ which multiply $(r_1,\ldots,r_d)$ by elements of
  $(r_{d+1}, \ldots r_g)$, we can turn $(r_1,\ldots, r_g)$ into $(q_1,\ldots,q_d,
  r_{d+1},\ldots,r_g)$.  Repeating this using the first $d$ places
  instead lets us get to $f_0$ by further elements of $\Out(Q)$.
  Thus the action of $\Out(F_g)$ on $\E_g$ is transitive.
\end{proof}

Before moving on, let us note that $\bar{d}$ is no more than the length
of a properly increasing sequence of subgroups
\[
\{1\} = Q_0 < Q_1 < \cdots < Q_k = Q.
\]
Since the index of one $Q_i$ in the next is at least $2$, we get that
$\bar{d} \leq \log_2\num{Q}$.  Thus in Proposition~\ref{act_on_E}, we
can take $g \geq 2 \log_2\num{Q}$ as $d + \bar{d} \leq 2 \bar{d}$.  Again,
see \cite{Pak2001} for sharper results about particular groups.

Since we are particularly interested in non-abelian simple groups, the
following theorem is of interest.
\begin{theorem}\label{act_on_E_simple}
  Let $Q$ be a non-abelian finite simple group and $g \geq 3$.  Then
  there is some orbit $\E'$ of $\E_g$ where $\Out(F_g)$ acts on $\E'$
  as either the full symmetric or alternating group.
\end{theorem}
In particular, when $\Out(F_g)$ acts transitively, the image of the
homomorphism $\Out(F_g) \to \Sym(\E_g)$ is either the whole thing or
has index two.  Theorem~\ref{act_on_E_simple} was proved by Gilman
\cite{Gilman77} for $g \geq 4$ and improved to $g = 3$ by Evans
\cite{Evans1993}.  The statements they give seem to put
restrictions on the simple group, but their conditions actually hold
for all finite simple groups; see the discussion in
\cite[Thm.~2.4.3]{Pak2001}.  The proof of Theorem~\ref{act_on_E_simple}
is fairly short and elementary but we will not reproduce it here; see
\cite{Gilman77} for more.  We should mention, though, that as we
stated it the result is more complex because it uses that every finite
simple group is generated by two elements, one of which has order two;
as such it uses part of the Classification of Finite Simple Groups.
However, for the way we will use Theorem~\ref{act_on_E_simple}, we only
need that the conclusion holds for all large $g$, and this is exactly
what Gilman's elementary argument shows.

\section{Covers where the quotient is simple}
\label{simple_quotients}

This section is devoted to the proof of the following theorem, which
gives a complete asymptotic picture in the case of a non-abelian
simple group
\begin{theorem}\label{simple_probs}
  Let $Q$ be a finite non-abelian simple group.  Then as the genus $g$ goes to
  infinity,
  \[
  p(Q, g) \to 1 - e^{-\mu} \mtext{where} \mu = \num{ H_2(Q, \Z)}
  /\num{\Out(Q)}.
  \]
  Moreover, the limiting distribution on the number of $Q$-covers
  converges to the Poisson distribution with mean $\mu$.  
\end{theorem}

As an example, if $Q = \PSL{2}{\F_p}$ where $p$ is an odd prime then
$\mu = 1$.  Thus for large genus the probability of a $\PSL{2}{\F_p}$
cover is about $1 - e^{-1} \approx 0.6321$.  Unfortunately, the proof does
not give information about the rate of convergence, and we will
discuss this issue in detail later.

It is interesting to compare this result with the corresponding result
for a general random group coming from a balanced presentation.
Theorem~\ref{thm-simple-balanced} shows that as the number of
generators goes to infinity, the limiting distribution is also
Poisson.  However, the mean is much smaller, namely $1/\num{\Aut(Q)}
\leq 1/\num{Q}$.  In particular, the mean goes to $0$ as the size of the
groups increases.  In contrast, we just saw that for 3-manifold groups
the limiting mean of the number of $\PSL{2}{\F_p}$ covers is
independent of the size of the group. 

To prove the theorem, the thing that we need beyond
Section~\ref{heegaard_cover_intro} is a better understanding of the
action of the mapping class group $\M_g$ on the set of epimorphisms
$\A_g$ (notion is as in Section~\ref{heegaard_cover_intro}).  In
particular, we will show that for a fixed non-abelian simple group $Q$
the action is eventually by the full alternating group of each
component (Theorem~\ref{action_is_alternating}).  Before stating the
first lemma, we need some definitions.  For a finite group $Q$ and an
element $c \in H_2(Q)/\Out(Q)$ set
\[
\A_g^c = \setdef{f \in \A_g}{c_f = c}.
\]
Recall that the action of a group on a set $\Omega$ is $k$-transitive if
the action on $k$-tuples of \emph{distinct} elements of $\Omega$ is
transitive.  The following lemma is the heart of the proof of
Theorem~\ref{simple_probs}:
\begin{lemma}\label{action-is-very-trans}
  Let $Q$ be a non-abelian finite simple group, with $c \in
  H_2(Q)/\Out(Q)$.  Let $k$ be a positive integer.  Then for all large
  $g$, the action of $\M_g$ on $\A_g^c$ is $k$-transitive.
\end{lemma}

\begin{proof}
  Consider the direct product $Q^k$ of $k$ copies of $Q$.  We will use
  $\A_g(Q)$ and $\A_g(Q^k)$ to denote epimorphisms to $Q$ and $Q^k$
  respectively.  As in Section~\ref{subsec_orbit_classification}, we
  will also consider the set $\A'_g(Q^k)$ of epimorphisms to $Q^k$
  where we do not mod out by $\Aut(Q^k)$.  By
  Theorem~\ref{hom_classify_orbits_prime}, for all large $g$, the map
  $\A'_g(Q^k) / \M_g \to H_2(Q^k)$ is a bijection; henceforth we will
  assume this holds.
  
  Now consider $k$-tuples $([f_1],\ldots,[f_k])$ and $([h_1],\ldots,[h_k])$ of distinct
  elements of $\A_g^c(Q)$.  Let $\tilde{c}$ be a lift of $c$ to
  $H_2(Q, \Z)$, and choose our representatives
  \[
  f_i, h_i \maps \pi_1(\Sigma_g) \to Q
  \]
  so that all the $c_{f_i}$ and $c_{h_i}$ are equal to $\tilde{c}$ in
  $H_2(Q)$.  Now consider the induced product maps $f, h \maps
  \pi_1(\Sigma_g) \to Q^k$.  Since $Q$ is simple and the $f_i$ represent
  \emph{distinct} classes in $\A_g$, the homomorphism $f$ is also
  surjective (Lemma~\ref{lem-simple-surjects-to-product}).  Similarly
  $h$ is also an epimorphism, and so both $f$ and $h$ are elements of
  $\A'_g(Q)$.
  
  As $Q$ is non-abelian, $H_1(Q, \Z) = 0$, and so $H_2(Q^k)$ is the
  direct sum $\left( H_2(Q) \right)^k$.  Thus $c_f = c_h =
  \tilde{c}^k$ in $H_2(Q^k)$.  By our choice of genus, there must be a
  $\phi \in \M_g$ such that $h = \phi \cdot f$.  Then $\phi$ carries
  $([f_1],\ldots,[f_k])$ to $([h_1],\ldots,[h_k])$ and so the action of $\M_g$
  is $k$-transitive.
\end{proof}

We will also need the following variant of the preceding lemma, whose
proof is essentially identical.
\begin{lemma} \label{lem-trans-alt}
  Consider a finite set $Q_i$ of non-abelian finite simple groups, and
  fix $c_i \in H_2(Q_i)/\Out(Q_i)$.  Suppose there are no isomorphisms
  from $Q_i$ to $Q_j$ which take $c_i$ to $c_j$.  Then the action
  of $\M_g$ on $\A_g^{c_1}(Q_1) \times \cdots \times \A_g^{c_n}(Q_n)$ is transitive.
\end{lemma}

Now we use Lemma~\ref{action-is-very-trans} to prove:
\begin{theorem}\label{action_is_alternating}
  Let $Q$ be a non-abelian finite simple group, with a fixed class $c$ in
  $H_2(Q)/\Out(Q)$.  Then for all large $g$, the action of $\M_g$ on
  $\A_g^c$ is by the full alternating group $\Alt(\A_g^c)$.  Moreover,
  the actions on different $\A^c_g$ are independent, in the sense that
  the map $\M_g \to \prod \Alt(\A_g^c)$ is surjective.
\end{theorem}

You can view this theorem as saying that the action of $\M_g$ on
$\A_g$ is nearly as mixing as possible, at least when the genus is
large.  As such, it is analogous to Goldman's theorem that the action
of $\M_g$ on the $\SU{2}$-character variety is ergodic for any genus
$\geq 2$ \cite{Goldman1997}.  While the proof here does not give an
explicit bound on when this mixing behavior occurs, we suspect that
it typically occurs as soon as $g \geq 3$.

\begin{proof}
  One consequence of the Classification of Simple Groups is that a
  group which acts $6$-transitively on a finite set $\Omega$ must contain
  $\Alt(\Omega)$ (see e.g. ~\cite[Thm.~7.3A]{DixonMortimerBook}).  By
  Theorem~\ref{action-is-very-trans}, the action of $\M_g$ on $\A_g^c$
  is 6-transitive for all large $g$; hence the action contains
  $\Alt(\A_g^c)$.  Since $\M_g$ is a perfect group for genus at least
  $2$, the image of $\M_g \to \Sym(\A_g^c)$ must be $\Alt(\A_g^c)$ and
  not $\Sym(\A_g^c)$, proving the first part of the theorem.
  
  For the independence, if $\M_g \to \prod \Alt(\A_g^c)$ is not surjective
  then there are distinct $c$ and $d$ such that $\M_g \to \Alt{ \A_g^c}
  \times \Alt{\A_g^d}$ is not surjective.  Then there is a bijection $\alpha
  \maps \A_g^c \to \A_g^d$ which is compatible with the $\M_g$ action
  on both of these sets; in particular, the action of $\M_g$ on pairs
  $(f, h) \in \A_g^c \times \A_g^d$ is not transitive.  But that
  contradicts Lemma~\ref{lem-trans-alt}, finishing the proof.
\end{proof}

You might wonder if we really need the Classification of Simple Groups
to prove this theorem.  In the case of the orbit $\A_g^0$, which is
what we need for Theorem~\ref{simple_probs}, we can replace the
Classification by some much less difficult results.  First choose $g$
large enough so that the action on $\A_g^0$ is 2-transitive.  Let $G$
be the image of $\M_g \to \Sym(\A_g^0)$.  Gilman showed (see
Theorem~\ref{act_on_E_simple}) that the action of the stabilizer of
$\E_g$ in $G$ is all of $\Alt(\E_g)$ for large $g$.  Thus $G$ has
subgroups $K \lhd H \leq G$ so that $H/K$ is isomorphic to $\Alt(\E_g)$;
in this context, $G$ is said to have a \emph{section} which is
$\Alt(\E_g)$.  Then Theorem 5.5B of \cite{DixonMortimerBook} says that
since $G$ is 2-transitive, either $G$ contains $\Alt(\A_g^0)$ or
$\num{\E_g} < 6 \log \num{\A_g^0}$.  The latter can not happen for large
$g$ by the results of Section~\ref{counting_A_and_E}, and so $M_g$
acts on $\A_g^0$ by the full alternating group.  Both Gilman's result
and the theorem about permutation groups are fairly elementary, and,
while certainly not trivial, they are orders of magnitude easier than
the Classification.

The main result of this section now follows easily:

\begin{proof}[Proof of Theorem~\ref{simple_probs}]
  We need to show that the distribution of the number of $Q$-covers
  converges to the Poisson distribution with mean
  $\num{H_2(Q)}/\num{\Out(Q)}$.  By the proceeding theorem, choose $g$
  large enough so that the action of $\M_g$ on $\A_g^0$ is by
  $\Alt(\A_g^0)$.  Moreover, by Lemmas~\ref{size_of_A_g},
  \ref{size_of_E_g}, and \ref{class-probability} we have that
  \[
  \num{\A_g^0} \sim \frac{\num{Q}^{2g - 2}}{\num{\Out(Q)}
    \num{H_2(Q)}} \mtext{and}
  \num{\E_g} \sim \frac{\num{Q}^{g-1}}{\num{\Out(Q)}}.
 \]
 Since the action of $\Alt(\A_g^0)$ on $\A_g^0$ is
 $(\num{\A_g^0}-2)$-transitive, the action is in particular
 $\num{\E_g}$-transitive.  Thus the number of $Q$-covers is
 distributed exactly as the number of intersections of $\E_g$ with a
 randomly chosen subset of $\A_g^0$ of size $\num{\E_g}$.  By the
 above, the expected number of covers is $\num{\E_g}^2/\num{\A_g^0} \sim
 \num{H_2(Q)}/\num{\Out(Q)}$.  Thus the distribution of $Q$-covers
 converges to a Poisson distribution with mean
 $\num{H_2(Q)}/\num{\Out(Q)}$, as desired.
\end{proof}

There is another proof of Theorem~\ref{simple_probs} that is worth
mentioning, which bypasses Theorem~\ref{action_is_alternating} and
instead relies directly on Lemma~\ref{action-is-very-trans} and some
elementary probability theory.  Let $X_g$ denote the random variable
whose distribution is the number of $Q$-covers for genus $g$.  For
comparison purposes, set $Y_g$ to be the random variable which is the
number of intersections of a randomly chosen subset of $\A_g^0$ of
size $\num{\E_g}$ with $\E_g$.  Let $\mu = \num{H_2(Q)}/\num{\Out(Q)}$.
We know that the $Y_g$ limit to the Poisson distribution with mean
$\mu$; we want to show the same for $X_g$.

Fix a positive integer $k$.  Lemma~\ref{action-is-very-trans} says
that the action of $\M_g$ on $\A_g^0$ is $k$-transitive for all large
$g$.  Think of $X_g$ as a sum of variables $I_f$ for $f \in \E_g$,
which are indicator functions for whether a particular $f$ extends.
The $k$-transitivity of the action implies that if we look at $k$
variables $I_{f_1}, \ldots, I_{f_k}$, then these are distributed exactly as
if we picked a $\phi \cdot \E_g$ as a random subset of $\A_g^0$ of size
$\num{E_g}$.  In other words, $I_{f_1}, \ldots, I_{f_k}$ have the same
distribution if we decomposed $Y_g$ instead of $X_g$.  Thus the
$k^{\mathrm{th}}$-moment of $X_g$ is equal to that of $Y_g$.
Therefore the $k^{\mathrm{th}}$-moments of the $X_g$ converge to the
$k^{\mathrm{th}}$-moments of the Poisson distribution with mean $\mu$.
This is enough to show that the $X_g$ converge to that Poisson
distribution (see e.g.~\cite[Thm.~19]{BollobasBook1985}; see also
\cite[Ch.~8]{AlonSpencerBook} for related statements).

\subsection{Sequences of simple covers}

In trying to understand the set of all finite covers of a \3-manifold
$M$, a natural question is whether $\pi_1(M)$ has a quotient
which is in some particular class of groups, such as alternating
groups, or groups of the form $\PSL{2}{\F_p}$.  If $M$ is hyperbolic,
then the congruence quotients give an infinite number of
$\PSL{2}{\F_p}$-covers.  Moreover, Lubotzky has shown that congruence
subgroups have density zero among all finite-index subgroups of
$\pi_1(M)$ \cite{Lubotzky95}.  However, even the answer to the
following question is unknown:
\begin{question}\label{has-alternating-quotient}
  Does every hyperbolic \3-manifold have a cover with group $A_n$ for some
  $n \geq 5$?
\end{question}
While this question for hyperbolic surfaces is trivially yes (since
those groups surject onto a free group of rank 2), the 2-dimensional
case becomes interesting if we generalize to orbifolds.  In this case,
Everitt \cite{Everitt2000} has shown that the fundamental group of any
hyperbolic 2-orbifold surjects onto $A_n$ for all large $n$; the hard
cases here are the $S^2(p,q,r)$ orbifolds, and there Everitt manages to
build explicit covers which realize these groups.  Recently, Liebeck
and Shalev have given another proof which studies the question from a
very different probabilistic point of view
\cite{LiebeckShalev2005a}.

Theorem~\ref{simple_probs} suggests that we might expect the
answer to Question~\ref{has-alternating-quotient} to be yes because
$\num{H_2(A_n)}/\num{\Out(A_n)} = 1$ for all large $n$.  Thus if
having covers with different groups are independent random events, we
would naively expect a random \3-manifold to have at least one, indeed
infinitely many, $A_n$ covers.  In this direction, we can show:
\begin{theorem}\label{random-alternating-quotient}
  Let $\epsilon > 0$.  Then for all sufficiently large $g$, the probability
  that the \3-manifold obtained from a random genus-$g$ Heegaard
  splitting has an $A_n$-cover for some $n \geq 5$ is at least $1-\epsilon$.
  Moreover, the same is true if we require that there be at least some
  fixed number $k$ of such covers.
\end{theorem}

What about similar questions where we consider a collection
$\{Q_i\}_{i=0}^\infty$ of non-abelian simple groups?  For this we will show:
\begin{theorem}\label{random-quotient-sequence}
  Let $\{Q_i\}_{i=0}^\infty$ be a sequence of distinct non-abelian
  finite simple groups; set $\mu_i = \num{H_2(Q_i)}/\num{\Out(Q_i)}$.
  Suppose further that $\sum \mu_i = \infty$.  Fix $\epsilon > 0$.  Then for
  all sufficiently large $g$, the probability that the \3-manifold
  obtained from a random genus-$g$ Heegaard splitting has an
  $Q_i$-cover for some $i$ is at least $1-\epsilon$.  Moreover, the same is
  true if we require that there be at least some fixed number $k$ of
  such covers.
\end{theorem}

Since $\num{H_2(A_n)}/\num{\Out(A_n)} = 1$ for all large $n$, this
theorem immediately implies Theorem~\ref{random-alternating-quotient}.
Another class where it applies is $Q_i = \PSL{2}{\F_{p_i}}$ where
$p_i$ is prime; an example of where it does not is $Q_i =
\PSL{2}{\F_{p^{i^2}}}$ where $p$ is a fixed prime.  In fact, it is
not hard to work out exactly which sequences $\{Q_i\}$ satisfy the
hypothesis from the Classification of Finite Simple Groups (see
\cite[\S3]{ATLAS}).  In particular, any sequence $\{Q_i\}$ must either
contain an infinite number of $A_n$ or an infinite number of Chevalley
groups.  In the former case, Theorem~\ref{random-quotient-sequence}
applies.  For the Chevalley case, these are groups of Lie type, and
let $q_i = p_i^{f_i}$ be the field of definition of $Q_i$.  Then
$\mu_i = c_i/f_i$ where $c_i$ is a non-zero rational number which is
universally bounded above and below (see \cite[\S3]{ATLAS}, especially
Table 5).  Thus Theorem~\ref{random-quotient-sequence} applies in this
case if and only if $\sum 1/f_i$ diverges.  Thus
Theorem~\ref{random-quotient-sequence} applies to $\PSL{2}{\F_{p^i}}$
but not $\PSL{2}{\F_{p^{i^2}}}$.  Now we will give the proof of the
theorem.
  
\begin{proof}[Proof of Theorem~\ref{random-quotient-sequence}]
  The key lemma is the following, which says that covers with
  different $Q_i$ are essentially independent.
  \begin{lemma}\label{lemma-prod-surj}
    Let $Q_1, \ldots, Q_n$ be non-abelian finite simple groups.  Then for
    all large $g$, the map $\M_g \to \prod \Alt(\A_g^0(Q_i))$ is surjective. 
  \end{lemma}
  First, let us see why this lemma implies the theorem.  For a fixed
  genus $g$, let $X_i$ denote the random variable corresponding to the
  number of $Q_i$ covers.  Fix a positive integer $n$.  Choose $g$
  large enough so that the lemma holds for $X_1,\ldots,X_n$, and so that
  the expectations $E(X_i)$ are very nearly $\mu_i$ for $1\leq i \leq n$.
  From the lemma, the $X_i$ are independent random variables which are
  nearly Poisson with mean $\mu_i$.  Thus their sum $X_1 + \cdots + X_n$ is
  nearly Poisson with mean $M_n = \mu_1 + \cdots + \mu_n$.  In particular,
  the probability of having a cover with group one of $Q_1, \ldots ,
  Q_n$ is about $1 - e^{-M_n}$.  Since $\sum \mu_i$ diverges, by
  increasing $n$ we can make this probability as close to $1$ as we
  like (of course, $g$ may well have to increase as we change $n$).
  Similarly, we can require some fixed number $k$ of $Q_i$ covers.
  This proves the theorem modulo the lemma.
  
  As for the lemma, by Theorem~\ref{action_is_alternating} gives a $g$
  so that $\M_g \to \Alt \big( \A_g^0(Q_i) \big)$ is surjective for $i
  \leq n$.  If the product $\M_g \to \prod \Alt(\A_g^0(Q_i))$ is not
  surjective, then because the factors are simple, there must be a
  pair $i, j$ such that $\num{\A_g^0(Q_i)} = \num{\A_g^0(Q_j)}$ and a
  bijection $\A_g^0(Q_i) \to \A_g^0(Q_j)$ which is compatible with the
  $\M_g$ action on both of these sets.  In particular, the action of
  $\M_g$ on pairs $(f, h) \in \A_g^0(Q_i) \times \A_g^0(Q_j)$ is not
  transitive.  But that would contradict Lemma~\ref{lem-trans-alt},
  thus proving Lemma~\ref{lemma-prod-surj} and hence the theorem.
\end{proof}

\section{Homology of random Heegaard splittings}
\label{homology_Heegaard}

In this section, we work out the distribution of the homology of a
\3-manifold coming from a random Heegaard splitting of genus $g$.
First, let us set up the point of view we will take in this section.
As always, let $H_g$ be a genus-$g$ handlebody, $\Sigma_g$ its boundary,
and $\M_g$ the mapping class group of $\Sigma_g$.  Given $\phi \in \M_g$ let
$M_\phi$ be the resulting \3-manifold.  Let $J$ be the kernel of the map
$H_1(\Sigma_g; \Z) \to H_1(H_g; \Z)$.  Then $H_1(M_\phi; \Z)$ is the quotient
of $H_1(\Sigma_g; \Z)$ by the subgroup $\pair{ J, \phi_*^{-1}(J)}$.  Note
that counting algebraic intersections between two cycles gives $H_1(\partial
H_g; \Z)$ a natural symplectic form.  With respect to this symplectic
structure $J$ is a Lagrangian subspace.  The natural map $\M_g \to
\Sp{2g}{\Z}$ coming from the action of mapping classes on
$H_1(\Sigma_g;\Z) \cong \Z^{2g}$ is surjective.  Thus understanding
the distribution of the homology of $M_\phi$ boils down to considering the
distribution of Lagrangians under the action of $\Sp{2g}{\Z}$.

\subsection{Symplectic groups over $\F_p$}   

As with the case of random finitely presented groups, we will work one
prime at time.  We will work in the following framework, where $p$ is
a fixed prime (or prime power).  Let $J$ be a vector space over $\F_p$
of dimension $g$, and let $K$ be its dual.  A nice way of thinking
about the symplectic vector space over $\F_p$ of dimension $2g$ is to
consider $V = J \oplus K$ with symplectic form given by $\pair{ (j_1, k_1)
  , (j_2, k_2) } = k_1(j_2) - k_2(j_1)$.  Notice that both $J$ and $K$
are Lagrangian subspaces of $V$.

Now set $G = \Sp{}{V} =\Sp{2g}{\F_p}$.  Identify $V$ with $H_1(\Sigma_g;
\F_p)$ so that the kernel of the map to $H_1(H_g; \F_p)$ is the
Lagrangian $J$.  The action of $\M_g$ on homology gives
a surjection $\M_g \to G$.  If $\phi \in \M_g$ is the result of a long
random walk then its image in $G$ is nearly uniformly distributed.
Now by the analog of the Gram-Schmidt process for a symplectic vector
space, it is easy to show that the action of $G$ on the set of
Lagrangians is transitive.  Thus $\phi^{-1}_*(J)$ is nearly uniformly
distributed among all Lagrangians.  Hence the asymptotic probability
that $\dim( H_1(M_\phi ; \F_p) ) = k$ is simply the ratio of the number
of Lagrangians $L$ in $V$ such that $\dim(J \cap L) = k$ to the total
number of Lagrangians.

\subsection{Counting Lagrangians}

First, let us count the total number of Lagrangians in $V$; we will do
this by computing the size of the stabilizer $S \leq G$ of the fixed
Lagrangian $J$.  Notice we have a natural homomorphism $S \to
\GL{}{(J)}$.  This map is surjective, for if $A \in \GL{}{(J)}$ then
recalling that $K = J^*$, we have that the map $A \oplus (A^*)^{-1}$ of $V$
respects the symplectic form and hence is in $S$.  An element $s$ of the
kernel of $S \to \GL{}{(J)}$ is determined by its restriction to $K$.
Write the restriction as $B_1 \oplus B_2$ where $B_1$ maps $K$ into $J$ and
$B_2$ maps $K$ into $K$.  For each $j \in J$ and $k \in K$ we have $\pair{
  k, j } = \pair{ s(k), s(j) } = \pair{ B_1(k) + B_2(k) , j } =
\pair{B_2(k), j}$ and hence $B_2 = \mathrm{Id}$ on $K$.  Moreover, the
requirement that $s$ preserve the symplectic form is then equivalent to
$\pair{ B_1(k_1), k_2 } = \pair{B_1(k_2), k_1}$ for each $k_1, k_2 \in
K$.  Such a $B_1$ is called \emph{symmetric}.  If we write $B_1$ with
respect to a pair of dual basis for $J$ and $K$ then this is
equivalent to the resulting matrix being symmetric.  Thus
\begin{equation*}
\begin{split}
\num{S} =  \num{ \{ \text{$g \times g$ symmetric matrices
    }  \} } \cdot \num{\GL{}{(J)}}&=  p^{g(g+1)/2}  \prod_{k=1}^g \left( p^g - p^{k - 1}\right) \\ 
 &= p^{g^2}\prod_{k=1}^g \left( p^k - 1\right) .
\end{split}
\end{equation*}

To complete the count of Lagrangians, we also need to know the order
of $G$ itself.  Fix a basis for $V$ consisting of a basis $\{e_i\}$
for $J$ and the dual basis $\{e^i\}$ for $K$.  To build an element of
$G$, there are $p^{2g} - 1$ choices for the image $v$ of $e_1$, and
as the image $w$ of $e^1$ must satisfy $\pair{w , v} = 1$, there are
$p^{2g - 1}$ choices for $w$.  The other basis vectors must be sent
into $\vecspan(v,w)^\perp$.  Thus inductively one has
\[
\num{G} = \prod_{k = 1}^g (p^{2 k} - 1)p^{2 k - 1} = p^{g^2}\prod_{k = 1}^g  \left(p^{2 k} - 1\right).
\]

Thus, the number of Lagrangians in $V$ is 
\[
\frac{\num{G}}{\num{S}} = \prod_{k = 1}^g \frac{p^{2 k} - 1}{p^k - 1} = \prod_{k = 1}^g \left( p^k + 1 \right).
\]

\subsection{Counting transverse Lagrangians}  

Now we will calculate the probability that $H_1(M_\phi ; \F_p) = 0$.  To
do this, we need to count the number of Lagrangians $L$ which are
transverse to our base Lagrangian $J$.  Given such an $L$, projection
of $L$ onto $K$ is surjective; thus we can view
it as the graph of a map $B \maps K \to J$.  The
requirement that the graph $L$ is Lagrangian is equivalent to
\[ 
\pair{ k_1 +  B(k_1) , k_2 + B(k_2)} = 0 \iff \pair{B(k_1), k_2} = \pair{B(k_2),  k_1} \mtext{for all $k_1, k_2 \in K$.}
\]
That is, $B$ is symmetric in the same sense as before.  Thus the
number of Lagrangians transverse to $J$ is equal to the number of
symmetric $g \times g$ matrices over $\F_p$, that is $p^{g(g+1)/2}$.
Combining this with our count of Lagrangians gives the following:
\begin{theorem}\label{thm-homology-Heegaard}
  Fix a Heegaard genus $g$.  Then the asymptotic probability that
   \[
 H_1(M_\phi ; \F_p) = 0 \mtext{is} \prod_{k = 1}^g \frac{1}{1 + p^{-k}}
  \]
  as we let the length of the random walk generating $\phi \in \M_g$ go
  to infinity.
\end{theorem}
Notice that for a fixed genus $g$ the larger $p$ is the closer this
probability is to $1$.  Thus the asymptotic probability that
$H_1(M_\phi; \Q)$ vanishes is $1$, as you would expect.  On the other
hand, if we consider a pair of primes $p$ and $q$ then
Lemma~\ref{lem-simple-surjects-to-product} implies that the induced
map
\[
\M_g \to \Sp{2g}{\F_p} \times \Sp{2g}{\F_q}
\]
is surjective.  Thus homology with $\F_p$ versus $\F_q$ coefficients
are asymptotically independent variables.  Moreover, the sum of the
probabilities that $H_1(M_\phi; \F_p) \neq 0$ diverges (albeit slowly, like
the harmonic series).  Thus the probability that some $H_1(M_\phi; \F_p)
\neq 0$, where $p$ is larger than some fixed $C$, goes to $1$ as the
complexity of $\phi$ goes to $\infty$.  Hence the expected size of $H_1(M_\phi
; \Z)$ goes to infinity as the complexity of $\phi$ increases.
Summarizing:
\begin{corollary}\label{cor-things-dont-repeat}
  Fix a Heegaard genus $g$.  Then with asymptotic probability $1$ the
  rational homology $H_1(M_\phi; \Q)$ vanishes.  However, for a fixed $C$
  the probability that $\num{H_1(M_\phi ; \Z)} > C$ goes to 1 as the
  complexity of $\phi \to \infty$.
\end{corollary}
Notice that the second conclusion implies that while a fixed manifold
$N$ occurs with many different $\phi$, as the complexity of $\phi$
goes to infinity the probability that $N = M_\phi$ goes to zero.  Thus
asymptotic probabilities are, in particular, not highly disguised
statements about some finite set of manifolds.

It is also interesting to compare Theorem~\ref{thm-homology-Heegaard}
to the corresponding results for the group $\Gamma$ of random $g$-generator
balanced presentations.  In that case, by
Proposition~\ref{prop-homology-random-group}, the probability that
$H_1(\Gamma; \F_p) = 0$ is
\[
\prod_{k=1}^g\left( 1 - p^{-k} \right).
\]
Comparing term by term shows that this is less than the same
probability for $\M_\phi$; thus the homology of the random group is
\emph{more} likely to be non-zero than for the random \3-manifold.
Note also that in both cases the probabilities have the same limit
(namely $0$) as $p$ goes to infinity.  Both of these facts are the
reverse of what we saw for non-abelian simple groups.

\subsection{General distribution} 

To complete the picture, we need to find the probability that
$\dim( H_1(M_\phi; \F_p) ) = d$ for general $d$.  To do this, we start
by parameterizing those Lagrangians $L$ where $\dim( L \cap J) = d$.
Let us fix $A = L \cap J$.  Any Lagrangian $L$ intersecting $J$ in $A$ is
contained in $A^\perp = J \oplus \mathrm{Ann}(A)$, where $\mathrm{Ann}(A)$ is
the subset of $K = J^*$ which annihilates $A$.  Note that $A^\perp/A = (
J/A ) \oplus \mathrm{Ann}(A)$ inherits a natural symplectic form, and $L$
projects to a Lagrangian in $A^\perp/A$ which is transverse to $J/A$.
Indeed, Lagrangians in $A^\perp$ which are transverse to $J/A$ exactly
parameterize Lagrangians in $V$ which intersect $J$ in $A$.  Thus the
latter are parameterized by $(g-d) \times (g-d)$ symmetric matrices.  Hence
the number of Lagrangians intersecting $J$ in a $d$ dimensional
subspace is:
\begin{multline*}
 \num{ \left\{ \text{$(g - d) \times (g - d)$ symmetric matrices } \right\} } \cdot \num{ \left\{ \text{dim $d$ subspaces of $J$} \right\} } \\ 
= p^{ (g - d + 1)(g - d)/2} \prod_{k = 1}^d \frac{p^{g - k + 1} - 1}{p^k - 1}.
\end{multline*}
If we set $c_d$ to be the probability that $H_1( M_\phi ; \F_p)$ has
dimension $d$, then a concise way of summarizing this distribution is 
\begin{equation*}
\begin{split}
c_0 &= \prod_{k=1}^g \left(1 + p^{-k} \right)^{-1} \mtext{and} \\
\frac{c_d}{c_0} &= \prod_{k = 1}^d  \frac{1 - p^{-g + k - 1}}{p^k - 1} \mtext{for $1 \leq d \leq g$.}
\end{split}
\end{equation*}

\subsection{Large genus limit}  

As with the case of simple groups, the distributions of $\dim
H_1(M_\phi; \F_p)$ have a well-defined limit distribution as the genus
$g$ goes to infinity.  Namely, it is the probability distribution
given by
\begin{equation*}
\begin{split}
\widetilde{c}_0 &= \prod_{k=1}^\infty \left(1 + p^{-k} \right)^{-1} \mtext{and} \\
\frac{\widetilde{c}_d}{\widetilde{c}_0} &= \prod_{k = 1}^d  \left( p^k - 1 \right)^{-1} \mtext{for $1 \leq d$.}
\end{split}
\end{equation*}
The reader can check that this is really a probability distribution,
i.e.~that it has unit mass.  One way to do this is to use that the
finite approximates have unit mass, and then observe that $c_d/c_0$ is
an increasing function of $g$.

\subsection{$p$-adic point of view}  

It is possible to work out the $p$-adic distribution of homology,
analogous to Section~\ref{sec-balanced-abelian}.  However, we will not
go into this here.  The key technical tool needed is the ability to
count \emph{non-singular} $g \times g$ matrices over $\F_p$, which is done
in \cite{Carlitz54} and \cite{MacWilliams69}; see also
\cite[\S4]{Stanley98}.

\comment{ The number of non-singular symmetric $g \times g$ matrices over $\F_p$ is
$
  \prod_{k \text{ \rm even }}  p^k \prod_{k \text{ \rm odd }}  (p^k - 1)  \text{ \rm \quad where $1 \leq k \leq g$. }$
   double-check about if used as I have it right as I reformatted it
  a bit for readability.}

\section{Homology of finite-sheeted covers}\label{cover-homology}

In this section, we try to determine the probability that a cover has
$\beta_1 > 0$, where the covering group $Q$ is fixed.  As before we work
in the context of random Heegaard splittings.  Our original goal was
to find groups where this probability is positive, and so demonstrate
that the Virtual Haken Conjecture is true in many cases, perhaps with
probability 1.  Unfortunately, our results are in the other
direction; in particular the main result in this section is:
\begin{theorem}\label{thm-no-homology-abelian-cover}
  Let $Q$ be a finite abelian group.  The probability that the
  3-manifold obtained from a random Heegaard splitting of genus 2 has
  a $Q$-cover with $\beta_1 > 0$ is $0$.
\end{theorem}
The restriction to genus 2 is almost certainly artificial and simply
due to a technical difficulty in the general case that we were unable
to overcome.  

What about when $Q$ is non-abelian?  We did some computer experiments,
and could not find any groups $Q$ for which $Q$-covers seemed to have
$\beta_1 > 0$ with positive probability.  However, these experiments are
not completely convincing, even for the groups examined; as you will
see in the proof of Theorem~\ref{thm-no-homology-abelian-cover}, one
needs to consider certain subgroups of the mapping class group of
truly staggering index, even for fairly small $Q$.  Thus one can not
rule out experimentally some very small probability that a $Q$-cover
has $\beta_1 > 0$.  We still suspect that the probability of $\beta_1 > 0$
is $0$ for all $Q$, and would be surprised if there was a $Q$ for
which the probability of $\beta_1 > 0$ was greater than, say, $10^{-4}$.

Let us now contrast these results with our earlier experimental
results in \cite{DunfieldThurston:experiments}.  There, we examined
about $11{,}000$ small volume hyperbolic 3-manifolds.  Of these, more
than $8{,}000$ had 2-generator fundamental groups, and therefore
almost certainly have Heegaard genus 2.  About $1.7\%$ of these genus
2 manifolds had a $\Z/2$ cover with $\beta_1 > 0$, and at least $1.7\%$
had a $\Z/3$ cover with $\beta_1 > 0$.  Overall, $7.3\%$ of the genus 2
manifolds have an abelian cover with $\beta_1 > 0$ (of the full sample of
$11{,}000$ manifolds, this rises to $9.6\%$).  In comparing this with
Theorem~\ref{thm-no-homology-abelian-cover}, however, it is important
to keep in mind that the covering group $Q$ is fixed in the theorem.
In particular, we do not know the answer to the following question,
even experimentally:
\begin{question}
  Let $M$ be obtained from a random Heegaard splitting of genus $g$.
  What is the probability that the maximal abelian cover of $M$ has
  $\beta_1 > 0$?
\end{question}

Compared to the experiments in \S5 of
\cite{DunfieldThurston:experiments} which dealt with non-abelian
simple covers, the contrast becomes much more marked.  In
\cite[\S5]{DunfieldThurston:experiments}, we found for a fixed such
group $Q$ that a very large proportion of the covers ($> 50\%$) had
$\beta_1 > 0$; indeed, the expectation for $\beta_1$ grew linearly with
$\num{Q}$.  However, in our experiments for random genus 2 Heegaard
splittings we found that the probability seems to be 0 for the first
few group ($A_5$, $\PSL{2}{\F_7}$, $A_6$, and $\PSL{2}{\F_8}$), and we
expect that this pattern continues for all simple groups.

\subsection{Homology and subgroups of the mapping class group}  

The goal of this section is to prove
Theorem~\ref{thm-no-homology-abelian-cover}.  While
Theorem~\ref{thm-no-homology-abelian-cover} restricts to genus 2 and
abelian covering groups, with the exception of
Lemma~\ref{lemma-ratio-to-zero} this section will be done without
these restrictions.  So from now on, we fix a Heegaard genus $g \geq 2$
and an arbitrary finite group $Q$.  Let $H$ be a handlebody of that
genus, and $\Sigma$ be $\partial H$.  Let $\M$ denote the mapping class group of
$\Sigma$, with a fixed generating set $T$.  We focus on a single
epimorphism $f \maps \pi_1(\Sigma) \to Q$ which extends over $H$.
Proposition~\ref{makes_sense} shows that for $\phi$ the result of a
random walk in $\M$, the probability that $f$ extends to the
\3-manifold $N_\phi$ converges to a positive number as the length of the
walk goes to infinity.  We want to show that the probability that the
corresponding cover of $N_\phi$ has $\beta_1 > 0$ is $0$; as there are only
finitely many choices for $f$ this will suffice to prove
Theorem~\ref{thm-no-homology-abelian-cover}.

Recall from the proof of Proposition~\ref{makes_sense} that $f$
extends over $N_\phi$ if and only if $\phi \cdot f = f \circ \phi^{-1}_*$
extends over $H$.  We will consider each of the possibilities for $\phi
\cdot f$ individually.  We start with the case that $\phi \cdot f = f$ in the
set $\A$ of all such epimorphisms, up to $\Aut(Q)$.  This case is a
little simpler, and it forms the core for the argument in general.
Let $K$ be the kernel of $f \maps \pi_1(\Sigma) \to Q$.  We are focusing on
the subgroup
\[
\M_f = \setdef{ \phi \in \M}{ \phi \cdot f = f \ \mbox{in $\A$, or equivalently $\phi_*(K) = K$} },
\]
which has finite index in $\M$.  The mapping classes in $\M_f$ are
exactly those which lift to the cover $\tilde{\Sigma} \to \Sigma$ corresponding
to $K$, and let $\tilde{\M}_f$ be the group of all such lifts, modulo
isotopy within this restricted class of homeomorphisms of
$\tilde{\Sigma}$.  Regarding $Q$ as the group of covering translations
of $\tilde{\Sigma}$, we have the exact sequence
\[
1 \to Q \to \tilde{\M}_f \to \M_f \to 1.
\]
Here, we are using that $g \geq 2$, as in the torus case $Q$ need not
inject into $\tilde{\M}_f$.  Note that $Q$ is not central unless
$\M_f$ acts trivially on $\pi_1(\Sigma) / K$.

If $\phi \in \M_f$ then the $Q$-cover corresponding to $f$ will be
denoted $\tilde{N}_\phi$.  If $\tilde{\phi} \in \tilde{\M}_f$ is a lift of
$\phi$, then $\tilde{N}_\phi = N_{\tilde{\phi}}$.  The covering group $Q$
acts on the homology $H_1(\tilde{N}_\phi; \Q)$, and we can consider the
decomposition of this action into irreducible (over $\Q$)
representations.  This is the same decomposition as considering $L =
H_1(\tilde{N}_\phi; \Z)/(\mbox{torsion})$, which can be regarded as a
lattice in $H_1(\tilde{N}_\phi; \Q)$, and looking at the sublattices on
which $Q$ acts irreducibly.  (The direct sum of these sublattices has
finite index in $L$, but is not necessarily all of $L$.)  We will
examine $L$ by looking at the submodules individually.  Equivalently,
for each irreducible representation $\rho \maps Q \to \GL{}(V_i)$ where
$V_i$ is a lattice, we consider the homology $H_1(N_\phi, V_i)$ with
coefficients twisted by $\rho$.

To understand the homology of $\tilde{N}_\phi$, we first consider the
action of the covering group $Q$ on $V = H_1(\tilde{\Sigma}; \Z)$.  Again,
the action decomposes into a sum of (rationally) irreducible
representations.  We group these representations by isomorphism type,
and so express a finite index sublattice of $V$ as $V_0 \oplus V_1 \oplus \cdots
\oplus V_n$, where the action on $V_i$ is a direct sum of copies of a
single representation, which differs for distinct indices $i$.  The
action of $\tilde{\M}_f$ on $V$ preserves the decomposition into $V_i$
because $Q$ is a normal subgroup of $\tilde{\M}_f$.

Now consider the cover $\tilde{H}$ of the handlebody $H$ corresponding
to $f$, and set $W$ to be the kernel of $V \to H_1(\tilde{H}; \Z)$.  Then
we have $H_1(\tilde{N}_\phi) = V \big/ \left< W, \tilde{\phi}_*^{-1}(W)
\right> $.  As the action of $Q$ on $\tilde{\Sigma}$ extends over
$\tilde{H}$, the kernel $W$ is a $Q$-invariant subset of $V$.  Hence
it also decomposes into $W_0 \oplus W_1 \oplus \cdots \oplus W_n$ where $W_i = W \cap
V_i$.

\begin{lemma}\label{lemma-W-half-dim-V}
Each $W_i$ is half-dimensional in $V_i$.
\end{lemma}

\begin{proof}  
  
  Basically, we can compute both dimensions explicitly using the
  Euler characteristic.  Consider an irreducible rational $Q$-module
  $T \cong \Q^d$ corresponding to some $V_i$.  Suppose $X$ is a finite CW
  complex with a epimorphism $\pi_1(X) \to Q$.  Then the Euler
  characteristic of the twisted homology $H_*(X, T)$ is just $\dim (T) \chi
  (X)$, as is clear from the chain complex used to compute the former.
  Moreover, $H^0(X, T)$ and $H_0(X, T)$ both vanish because $T$ is
  irreducible; in general, $H^0(X,T)$ is the submodule of $T$
  consisting of invariant vectors, and $H_0(X, T)$ is the module of
  co-invariants ($= V/ \spandef{ q v - v}{q \in Q, v \in T}$)
  \cite[\S III.1]{Brown82}.
 
  Consider the homology $H_*( H, T)$ of the handlebody $H$.  As
  $H$ is homotopy equivalent to a bouquet of $g$ circles, we have that
  the only non-zero homology is in dimension $1$.  Thus
  \[
  \dim H_1(H, T) = - \chi ( H_*(H, T) ) = (g - 1) \dim T 
  \]
  and so $H_1( \tilde{H}, \Q )$ contains exactly $(g - 1)$ copies of
  $T$.
  
  Next, consider $H_*(\Sigma, T)$.  Poincar\'e duality implies $H_2(\Sigma, T)
  \cong H^0(\Sigma, T) = 0$   \cite[VIII.10]{Brown82}.  Thus again we have
  \[
  \dim H_1(\Sigma, T) = - \chi ( H_*(\Sigma, T) ) =  (2 g  -  2) \dim T.
  \]
  and so $V_i$ consists of $(2g - 2)$ copies of $T$.  Counting
  dimensions now shows that $W_i$ must consist of $(g - 1)$ copies of
  $T$, proving the lemma.
\end{proof}

We now focus attention on one summand $V_i$, looking at the homology
piece by piece.  As we only care about rational homology, set
\[
 U_i = \left( V_i \big/ \left< W_i , \tilde{\phi}_*^{-1}(W_i)\right>\right) \otimes \Q .
\]  Consider the homomorphism $\rho \maps \tilde{\M}_f \to
\GL{}(V_i)$ induced by the projection $V \to V_i$.  We want to show
that $U_i$ almost surely vanishes, or equivalently that the image of
$\rho$ almost surely takes $W_i$ to a complementary subspace.  More
precisely we need:
\begin{claim}\label{claim-prob-hom-small}
  Given $\epsilon > 0$ there exists a $C_0$ such that the following holds.
  If $\phi$ is the result of a random walk in $\M$ of length $C \geq C_0$
  then the probability
  \[
  P \left\{ \phi \in \M_f \mtext{and} U_i \neq 0 \right\} < \epsilon.
  \]
\end{claim}
We will deduce the above claim with the help of the next lemma, but
first some notation.  Let $P$ be the orbit of $V_i$ in the
half-dimensional Grassmannian $\Gr(W_i)$.  Let $B$ be the $X \in P$
such that $X$ is not transverse to $V_i$.  We will want to work mod a
fixed large prime $p$.  We denote reduction mod $p$ by a bar, e.g.
$\overline{V}_i= V_i/pV_i$.  We are interested in $\overline{\rho} \maps
\tilde{\M}_f \to \GL{}{(\overline{V}_i)}$, whose image we call
$\overline{G}$.  We want to understand the action of $\overline{G}$ on
the half-dimensional Grassmannian $\Gr(\overline{W}_i)$.  Let
$\overline{P} \subset \Gr(\overline{V}_i)$ be the orbit of $\overline{W}_i$
under $\overline{G}$; equivalently, $\overline{P}$ is the image of $P$
under $V_i \to \overline{V}_i$.  We set $\overline{B}$ to be the
reduction of $B$ mod $p$. Thus for $\tilde{\phi} \in \tilde{\M}_f$ we
have $U_i \neq 0$ implies that
$\overline{\rho}(\tilde{\phi})(\overline{V}_i) \in \overline{B}$.  Now we
can state:
\begin{lemma}\label{lemma-ratio-to-zero}
  Assume that $Q$ is cyclic and $g=2$.  Then the ratio
  $\num{\overline{B}}/\num{\overline{P}}$ goes to $0$ as $p \to \infty$. 
\end{lemma}

For a finite abelian group, any irreducible representation over $\C$
is \1-dimensional and the image consists of roots of unity; hence it
factors through some cyclic quotient.  It easily follows that all
$\Q$-irreducible representations also factor through cyclic quotients.
Thus Lemma~\ref{lemma-ratio-to-zero} is sufficiently general for its
role in proving Theorem~\ref{thm-no-homology-abelian-cover}.  Before
proving the lemma, let us deduce Claim~\ref{claim-prob-hom-small} from
it.  Let $\tilde{\Gamma}$ be the stabilizer of $\overline{W}_i$ in
$\overline{P}$, which is a subgroup of $\tilde{M}_f$ of index
$\num{\overline{P}}$.  Let $\Gamma$ be the image of $\tilde{\Gamma}$ in
$\M_f$.  As $Q$ fixes $W_i$ it is contained in $\tilde{\Gamma}$ and so $[
\M_f : \Gamma ] = [\tilde{\M}_f : \tilde{\Gamma}]$.  Now, if $\phi \in \M_f$ then
$U \neq 0$ implies that $\phi$ is in one of $\num{\overline{B}}$ cosets of
$\Gamma$ in $\M_f$.  As always, if $\phi$ is the result of a long random
walk in $\M$ then the location of $\phi$ in the finite coset space
$\M/\Gamma$ is nearly uniform.  Thus for long walks
\begin{align*}
P \left\{ \phi \in \M_f \mtext{and} U_i \neq 0 \right\} \lessapprox
\frac{\num{\overline{B} }}{[\M:\Gamma]} &= \frac{1}{[\M: \M_f]} \cdot \frac{\num{\overline{B}}}{[\M_f : \Gamma]}\\
   &= \frac{1}{[\M: \M_f]} \cdot \frac{\num{\overline{B}}}{\num{\overline{P}}}
\end{align*}
By the lemma, we can make the rightmost term as small as we want,
which proves Claim~\ref{claim-prob-hom-small}.  We now turn to the
proof of the lemma.

\begin{proof}[Proof of Lemma~\ref{lemma-ratio-to-zero}]
  We think of the handlebody $H$ as the union of two solid tori $H_1$
  and $H_2$ glued along a disc.  Then $\pi_1(H)$ is freely generated by
  $\{x_1, x_2\}$ where each $x_i$ generates $\pi_1(H_i)$.  Given our $f
  \maps \pi_1(\Sigma) \to Q$ which extends over $H$, we can choose this
  decomposition so that the induced map $f \maps H \to Q$ is such that
  $f(x_1)$ generates $Q$ and $x_2$ generates the kernel of $f$.  The
  corresponding $Q$-cover of $H$ is easy to describe; it consists of
  the $\num{Q}$-fold cyclic cover of $H_1$ with $\num{Q}$ copies of
  $H_2$ glued on like spokes around a central hub.
  
  Now we construct an element $\phi \in \M_f$ which allows us to
  estimate the size of $\overline{P}$ from below.  Let $T_j$ be the
  torus with one hole that is $\Sigma \cap H_j$.  Choose simple closed
  curves $\alpha_j$ and $\beta_j$ in $T_j$ which meet in one point, and where
  $\alpha_j$ bound a disc in $H_j$.  Consider the element $\phi \in \M_f$
  which is a Dehn twist in $\beta_2$ followed by $\num{Q}$ Dehn twists
  in $\beta_1$.  A lift $\widetilde{\phi}$ of $\phi$ is easy to describe:
  the preimage of $\beta_2$ consists of $\num{Q}$ curves, one on each
  $T_2$ spoke, while the preimage of $\beta_1$ is a single curve running
  once round the $\widetilde{T}_1$ hub; the lift $\widetilde{\phi}$ is
  simply the Dehn twist along this disjoint collection of curves.  The
  important thing for us is that $\widetilde{\phi}$ takes the
  Lagrangian kernel of
  \[
  H_1(\widetilde{\Sigma}; \Z) \to H_1(\widetilde{H}; \Z)
  \]
  to something transverse with itself; moreover, the images of this
  Lagrangian under proper powers $\widetilde{\phi}^n$ are all
  distinct, mutually transverse subspaces.

  Now consider the particular summand $V_i$ of $H_1(\widetilde{\Sigma};
  \Z)$ in the case at hand.  For $\widetilde{\phi}$, we know that the
  orbit of the Lagrangian $W_i$ inside $V_i$ is infinite.  Therefore, by
  choosing $p$ large we can make $\overline{P}$ as large as we want.
  Turning now to understanding $\overline{B}$, by the Euler
  characteristic calculation of Lemma~\ref{lemma-W-half-dim-V}, $W_i$
  consists of a single irreducible $Q$-module.  Therefore, for
  $\widetilde{\psi} \in \widetilde{M}_f$ the intersection
  $\widetilde{\psi}(W_i) \cap W_i$ is either $\{ 0 \}$ or all of $W_i$ since
  this intersection is $Q$-invariant.  Thus $B$ consists solely of $W_i$ and
  $\num{\overline{B}} = 1$.  Since we can make $\num{\overline{P}}$ as
  large as we want, the Lemma follows.  
\end{proof}

This completes the proof of
Theorem~\ref{thm-no-homology-abelian-cover} in the case that $\phi \cdot f
= f$.  We now turn to the general case, where $\phi \cdot f = g$ for an
arbitrary $g \maps \pi_1(\Sigma) \to Q$ extending over $H$.  Consider
\[
\M_{f,g} = \setdef{ \phi \in \M }{ \mbox{$\phi \cdot f = g$ in  $\A$, or equivalently $\phi_*(K) = \ker g$} },
\]
which is no longer a subgroup but is a coset of $\M_f$.  Fix $\phi_0$
such that $\phi_0 \cdot g = f$.  Given $\phi \in \M_{f,g}$ we have $\phi_0 \circ \phi
\in \M_f$.  Schematically, we have:
\[\begin{CD}
  \tilde{\Sigma}_f @>\tilde{\phi_0 \circ \phi}>> \tilde{\Sigma}_f @>\tilde{\phi_0}^{-1}>> \tilde{\Sigma}_g \\ 
  @VVV @VVV @VVV \\ 
  \Sigma @>\phi_0\circ\phi >> \Sigma @>\phi_0^{-1}>> \Sigma  \\
\end{CD}
\]
where we have distinguished the $Q$-covers of $\Sigma$ corresponding to
$f$ and $g$ by subscripts.  Let $V_f = H_1(\tilde{\Sigma}_f ; \Z)$ and let
$W_f$ be the kernel of $V_f \to H_1(\tilde{H}_f; \Z)$; similarly, let $V_g$
and $W_g$ be the corresponding lattices for $\tilde{\Sigma}_f$.  Then
moving from the right hand column of the diagram to the middle we get
\[
H_1( \tilde{N}_\phi ) = V_g \big/ \left<W_g , \tilde{\phi}_*(W_f) \right> \cong  V_f \big/ \left< {\phi_0}_*(W_g), (\tilde{\phi}_0 \circ \tilde{\phi})_* (W_f)\right>.
\] 
Now the subspace ${\phi_0}_*(W_g)$ is $Q$-invariant, so as before
we can break this question into separate questions for each summand
$V_{f,i}$ of $H_1( \tilde{\Sigma}_f ; \Z)$ corresponding to an
irreducible representation.  Again, we are interested in the
orbit $P$ of the Lagrangian $W_{f,i} \subset V_{f,i}$ under the elements of
$\M_f$.  The only difference is that the set $B$ should now be taken
to be those $X \in P$ which are not transverse to
${\phi_0}_*(W_{g,i})$ (rather than to $W_{f,i}$).  If we now assume
that the genus is 2, then, as in the proof of
Lemma~\ref{lemma-ratio-to-zero} the $Q$ actions on
${\phi_0}_*(W_{g,i})$ and elements of $P$ are irreducible.  Thus
$B$ has either $0$ or $1$ element depending on whether
${\phi_0}_*(W_{g,i}) \in P$.  The general case can now be completed in
exactly the same way as before.  This finishes the proof of
Theorem~\ref{thm-no-homology-abelian-cover}.

\subsection{Possible generalizations}

It would be nice to remove the restrictions that $Q$ is abelian and
that the genus is $2$ from the hypotheses of
Theorem~\ref{thm-no-homology-abelian-cover}.  These two conditions are
actually quite separate from the point of view of the proof, so we
discuss them each in turn.  

First, the fact that $Q$ is abelian (or really, cyclic) ensured that
the cover $\tilde{\Sigma} \to \Sigma$ was concrete enough to exhibit an
element $\phi \in \M_f$ showing that each $P$ is infinite.  (Note that
the $\phi$ given in Lemma~\ref{lemma-ratio-to-zero} easily generalizes
to any genus.)  Such $\phi$ could probably also be found for other
groups, especially small cases like $S_3$.  More ambitiously, if one
wanted to do a whole class of simple groups, e.g. ~alternating groups,
one would be badly hampered if one could not remove the genus
restriction --- after all, Section~\ref{simple_quotients} only applies
in that case.

The genus 2 hypothesis is used at only one point --- to show
$\num{\overline{B}} = 1$ and thus allowing us to make
$\num{\overline{B}}/\num{\overline{P}}$ small simply by knowing that
$P$ is infinite.  For higher genus, we would need a more detailed
picture of the image
\[
\M_f \to \Sp{}(V)
\]
where $V = H_1(\tilde{\Sigma}; \Z)$ in order to compare the relative
sizes of $\overline{B}$ and $\overline{P}$.   

However, a more abstract point of view might also work to circumvent
this issue.  Consider a finitely generated subgroup $\Gamma$ of
$\Sp{2n}(\Z)$, which we think of as sitting inside $G = \Sp{2n}(\R)$.
Fix a standard integral Lagrangian $W$ in $\R^{2n}$ and set
\[
D = \setdef{g \in G}{g(W) \cap W \neq \{ 0 \}}
\]
which is a proper real-algebraic subvariety (but not subgroup) of $G$.
Now fix generators of $\Gamma$, and consider the probability that a
random walk of length $N$ in $\Gamma$ lies in $D$.  It seems very
reasonable that, as long as $\Gamma$ does not have a finite-index
subgroup which is contained in $D$, then this probability goes to $0$
as $N \to \infty$.

Indeed one can prove this with some additional hypothesis on the
Zariski closure $\overline{\Gamma}$ of $\Gamma$.  What one needs is the
following.  Consider the mod $p$ reduction $\overline{\Gamma}(\F_p)$ of
$\overline{\Gamma}$, which has a natural structure as an algebraic variety
over $\F_p$.   By hypothesis, $D(\F_p)$ is a proper subvariety.  It
follows that
\[
\frac{\num{D(\F_p)}}{\num{\overline{\Gamma}(\F_p)}} \to 0 \mtext{as} p \to \infty
\]
as each has roughly as many points as the projective space over $\F_p$
of the appropriate dimension (this is a weak form of the Weil
Conjecture).  Thus if one has that the mod $p$ reduction map
\begin{equation}\label{eq-approximation}
  \Gamma \to \overline{\Gamma}(\F_p)
\end{equation}
is surjective for arbitrarily large $p$, this would prove the desired claim.  

It is known that (\ref{eq-approximation}) holds when $\Gamma$ is a
lattice in $\overline{\Gamma}(\F_p)$, or if $\overline{\Gamma}(\F_p)$
is simple and simply connected.  (The latter is the celebrated
Nori-Weisfeiler Strong Approximation Theorem \cite{Weisfeiler1984},
see also \cite[\S Windows]{LubotzkySegal} for a discussion.)  It is unclear if either of these hypotheses hold in our setting.

\end{document}